\documentclass[12pt,a4paper]{article}
\usepackage[latin1]{inputenc}
\usepackage[english]{babel}
\usepackage{graphicx}
\usepackage{subfigure}
\usepackage{amsmath}
\usepackage{amsfonts}
\usepackage{amssymb}
\usepackage{eucal}

\begin{document}

\title{\vspace*{2cm} {\bf High order semilagrangian methods \\for the
BGK equation}}

\author{
M.\ Groppi$^{1}$, G.\ Russo$^{2}$, G.\ Stracquadanio$^{1}$ \bigskip \\
{\footnotesize  {\it $1$  Dept. of Mathematics and Computer Science, Univ. of
Parma, Italy}}\\
{\footnotesize {\it $2$  Dept. of Mathematics and Computer Science, Univ. of
Catania, Italy}}\vspace{.2cm} }
\date{}

\maketitle
\begin{abstract}
A new class of high-order accuracy numerical methods for the BGK model of the Boltzmann equation is presented. The  schemes are based on a semi-lagrangian formulation of the BGK equation; time integration is dealt with DIRK (Diagonally Implicit Runge Kutta) and BDF methods; the latter turn out to be accurate and computationally less expensive than the former. 
Numerical results
and examples show that the schemes are reliable and efficient for
the investigation of both rarefied and
fluid regimes in gasdynamics.

\end{abstract}
\section{Introduction}
In the kinetic theory of gases, the dynamics of a monoatomic rarefied gas system is described by the Boltzmann equation \cite{Cerc}. The numerical approximation of this equation is not trivial due to the complex structure of the collision operator. The BGK equation, introduced by Bhatnagar, Gross and Krook \cite{bgk} and independently by Welander \cite{wel} is a simplified model of the Boltzmann equation. In the BGK model the collision operator is substituted by a relaxation operator;  the initial value problem reads as
\begin{equation}
\begin{array}{l}
\dfrac{\partial f}{\partial t}+{v}\cdot\nabla_{{x}}f=Q_{BGK}[f]\equiv\dfrac{1}{\varepsilon}(M[f]-f), \;\;\;({x},{v},t)\in \mathbb{R}^{d}\times\mathbb{R}^{N}\times\mathbb{R}^{+} \label{BGK}\\\\
f({x},{v},0)=f_{0}({x},{v}),
\end{array}
\end{equation}
where $d$ and $N$ denote the dimension of the physical and velocity spaces respectively, and $\varepsilon^{-1}$ is the collision frequency, that, throughout this paper, is assumed to be a fixed constant for simplicity. $M[f]$ denotes the local Maxwellian with the same macroscopic moments of the distribution function $f({x},{v},t)$, and is given by
\begin{equation}
M[f]({x},{v},t)=\dfrac{\rho({x},t)}{\big[2\pi RT({x},t)\big]^{N/2}}\exp\bigg(-\dfrac{({v}-{u}({x},t))^{2}}{2RT(x,t)}\bigg),\label{Maxwellian}
\end{equation}
where $R$ is the  ideal gas constant and $\rho({x},t) \in \mathbb{R}^{+}$, $u({x},t) \in \mathbb{R}^{N}$ and $T({x},t)  \in \mathbb{R}^{+}$ denote the macroscopic moments of the distribution function $f$, that is: density, mean velocity and temperature, respectively. They are obtained in the following way
\begin{equation}
(\rho,\rho u,E)^{T}=\langle f\phi({v})\rangle,\;\;\text{where}\;\; \phi({v})=\Big(1,{v},\frac{1}{2}|v|^{2}\Big)^{T}, \,\, \langle g\rangle=\int_{\mathbb{R}^{N}}g({v})\,d_N{v}.\label{exact moments}
\end{equation}
The physical quantity $E({x},t)$ is the total energy that is related to the temperature $T(x,t)$ by the underlying relation:
$$E({x},t)=\frac{1}{2}\rho({x},t) u({x},t)^{2}+\frac{N}{2}\rho(x,t)RT({x},t).$$
\newline
\noindent
The BGK model (\ref{BGK}) satisfies the main properties of the Boltzmann equation \cite{bgk,wel}, such as conservation of mass, momentum and energy, as well as the dissipation of entropy. In details,
\begin{equation}
\langle M[f] \phi(v) \rangle= \langle f \phi(v)\rangle, \qquad \quad \int_{\mathbb{R}^{N}} Q_{BGK}[f] \log f d_Nv \leq 0\,.
\end{equation}

 The equilibrium solutions are clearly Maxwellians, indeed the collision operator vanishes for $f=M[f]$.
 The BGK model is computationally less expensive than the Boltzmann equation, due mainly to the simpler form of the collision operator, but it still provides qualitatively correct solutions for the macroscopic moments near the fluid regime\footnote{More precisely, from the BGK model, to zero-th order in $\varepsilon$, one obtains the compres\-sible Euler equations in the fluid-dynamic limit, while to first order in $\varepsilon$, the moments satisfy equations of compressible Navier-Stokes type, but with the wrong value for the Prandtl number. This problem can be fixed by resorting to the so-called ES-BGK model \cite{Ho}, but in the present paper we shall restrict to the classical BGK model.}. These two aspects, the lower computational complexity and the correct description of the hydrodynamic limit, explain the interest in the BGK model over the last years.  Without expecting to be exhaustive, we refer for instance to \cite{Perthame,LSR,PP07,Yun,Miess,AAP,PR} and the references therein for a more in-depth analysis of the various aspects (theoretical and numerical) of BGK models. In particular, in the last years a lot of numerical schemes have been proposed to solve the BGK equation in an efficient way; just to mention a few, the very recent papers \cite{Dimarco131,Dimarco132} concern  methods based on splitting techniques, while the scheme proposed in \cite{PP12} takes advantage from the explicit advancing in time of the macroscopic fields involved in the BGK operator.

The aim of this paper is to develop high order semilagrangian numerical schemes for the BGK equation.
  Semilagrangian methods for BGK models  have recently received increasing interest \cite{FilRus,RusSanYun}, since they  well describe either a rarefied or a fluid regime. The relaxation operator is treated implicitly and the semilagrangian treatment of the convective part avoids the classical CFL stability restriction. Moreover, in this work time integration is dealt with  BDF methods along characteristics, which turn out to be accurate but computationally less expensive than Diagonally Implicit Runge Kutta (DIRK) methods introduced in \cite{santpre} and analyzed in \cite{RusSanYun}.

The paper is organized as follows. In Section 2 the semilagrangian method is introduced and the first order method is described; in Section 3 higher order methods are presented, based on  BDF or DIRK schemes for time integrations; the possibility to avoid interpolation is also investigated in Section 4. For simplicity, all schemes are described for the 1+1D BGK model. In Section 5 we describe how to extend the methods to 1D in space and 3D in velocity in slab geometry (Chu reduction). Numerical results are shown in Section 6, with the aim of showing the performance and the accuracy of the proposed methods in various examples.

\section{Lagrangian formulation and first order\\ scheme}
\begin{figure}[h]
{\includegraphics[trim=3.2cm 21.3cm 3.8cm 4cm,scale=1]{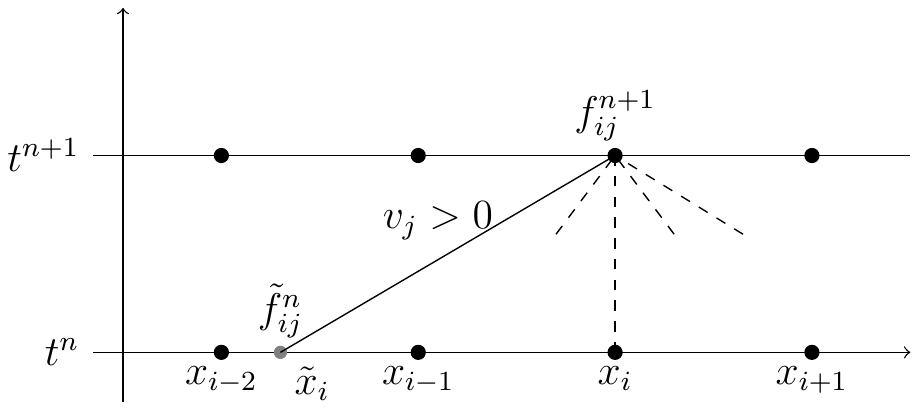}}
\caption{\footnotesize{Representation of the implicit first order scheme. The foot of the characteristic does not lie on the grid, and some interpolation is needed to compute $\tilde{f}^{n}_{ij}$.}}\label{first_order_fig}
\end{figure}

We shall restrict to the BGK equation in one space and velocity dimension (namely $d=N=1$ in (\ref{BGK}),(\ref{Maxwellian}).
In the Lagrangian formulation, the time evolution of $f(x,v,t)$ along the characteristic lines is given by the following system:
\begin{equation} \label{lagrBGK}
\begin{array}{l}
\dfrac{df}{dt}=\dfrac{1}{\varepsilon}(M[f]-f),\\\\
\dfrac{dx}{dt}=v,\\\\
x(0)=\tilde{x},\;\;f(x,v,0)=f_{0}(x,v)\;\;t\geq 0,\;\;x,\,v\,\in\, \mathbb{R}.
\end{array}
\end{equation}
For simplicity, we assume constant time step $\Delta t$ and uniform grid in physical and velocity space, with mesh spacing $\Delta x$ and $\Delta v$ respectively, and denote the grid points by $t^{n}=n\Delta t$, $x_{i}=x_{0}+i\Delta x,\, i=0,\cdots, N_{x}$, $v_{j}=j\Delta v,\, j=-N_{v},\cdots, N_{v}$, where $N_{x}+1$ and $2N_{v}+1$ are the number of grid nodes in space and velocity respectively, so that $[x_{0},x_{N_{x}}]$ is the space domain. We also denote the approximate solution $f(x_{i},v_{j},t^{n})$ by $f^{n}_{ij}$.\\
Relaxation time $\varepsilon$ is typically of the order of the Knudsen number, defined as the ratio between the molecular mean free path length and a representative macroscopic length; thus, the Knudsen number can vary in a wide range, from order greater than one (in rarefied regimes)  to very small values (in fluid dynamic regimes).\\
For this reason, if we want to capture the fluid-dynamic limit, we have to use an L-stable scheme in time. An implicit first order L-stable semilagrangian scheme (Fig. \ref{first_order_fig}) can be achieved in this simple way
\begin{equation}
f^{n+1}_{ij}=\tilde{f}^{n}_{ij}+\dfrac{\Delta t}{\varepsilon}(M[f]^{n+1}_{ij}-f^{n+1}_{ij}).\label{First order}
\end{equation}
The quantity $\tilde{f}^{n}_{ij}\backsimeq f(x_{i}-v_{j}\Delta t,v_{j},t^{n})$ can be computed by suitable reconstruction from $\{f^{n}_{\cdot j}\}$; linear reconstruction will be sufficient for first order scheme, while higher order reconstructions, such as ENO or WENO  \cite{WENO}, may be used to achieve high order avoiding oscillations. The convergence of this first order scheme has been studied in \cite{RusSanYun}.\\$M[f]^{n+1}_{ij}$ is the Maxwellian constructed with the macroscopic moments of $f^{n+1}$:
$$M[f]^{n+1}_{ij}=M[f](x_{i},v_{j},t^{n+1})=\dfrac{\rho^{n+1}_{i}}{\sqrt{2\pi RT^{n+1}_{i}}}exp\bigg(-\dfrac{(v_{j}-u^{n+1}_{i})^{2}}{2RT^{n+1}_{i}}\bigg).$$
This formula requires the computation of the discrete moments of $f^{n+1}$, through a numerical approximation of the integrals  in (\ref{exact moments}). This is obtained in the following standard way\footnote{Computing the moments using this approximation of the integrals has the consequence that the discrete Maxwellian $M^{n+1}_{ij} = \frac{\rho^{n+1}_{i}}{\sqrt{2\pi RT^{n+1}_{i}}}exp(-\frac{(v_{j}-u^{n+1}_{i})^{2}}{2RT^{n+1}_{i}})$ does not have the same discrete moments as $f^{n+1}_{ij}$. The discrepancy is very small if the distribution function is smooth and the number of points in velocity space is large enough, because midpoint rule is spectrally accurate for smooth functions having (numerically) compact support. However, for small values of $N_{v}$, such discrepancy can be noticeable. To avoid this drawback, Mieussens introduced a discrete Maxwellian \cite{Miess,Miess2}. The computation of the parameters of such Maxwellian requires the solution of a non linear system. A comparison between the continuous and discrete Maxwellian can be found, for example, in \cite{PupAl}. Here we shall neglect this effect, and assume that, using eq. (\ref{approximate moments}), $M^{n+1}_{ij}$ and $f^{n+1}_{ij}$ have the same moments with sufficient approximation.}:
\begin{equation}
\begin{array}{l}
\rho^{n+1}_{i}=\sum_{j=-N_{v}}^{N_{v}}f^{n+1}_{ij}\Delta v,\\\\
u^{n+1}_{i}=\dfrac{1}{\rho^{n+1}_{i}}\sum_{j=-N_{v}}^{N_{v}}v_{j}f^{n+1}_{ij}\Delta v,\\\\
E^{n+1}_{i}=\sum_{j=-N_{v}}^{N_{v}}\dfrac{1}{2}v_{j}^{2}f^{n+1}_{ij}\Delta v.\label{approximate moments}
\end{array}
\end{equation}
From now on, we will denote formulas in (\ref{approximate moments}) with the more compact notation: $(\rho_{i}^{n+1},(\rho u)^{n+1}_{i},E_{i}^{n+1})=m[f_{i\cdot}^{n+1}]$, where, in general, $m[f]$ will indicate the approximated macroscopic moments related to the distribution function $f$.

Now it is  evident that Equation (\ref{First order}) cannot be immediately solved for $f^{n+1}_{ij}$. It is a non linear implicit equation because the Maxwellian depends on $f^{n+1}$ itself through its moments. To solve this implicit step one can act as follow. Let us take the moments of equation (\ref{First order}); this is obtained at the discrete level multiplying both sides by $\phi_{j}\Delta v$, where $\phi_{j}=\{1,v_{j},v_{j}^2\}$ and summing over $j$ as in (\ref{approximate moments}). Then we have
$$\sum_{j}(f^{n+1}_{ij}-\tilde{f}^{n}_{ij})\phi_{j}=\dfrac{\Delta t}{\varepsilon}\sum_{j}(M[f]^{n+1}_{ij}-f^{n+1}_{ij})\phi_{j},$$
which implies that
$$\sum_{j}f^{n+1}_{ij}\phi_{j}\simeq\sum_{j}\tilde{f}^{n}_{ij}\phi_{j},$$
because, by definition, the Maxwellian at time $t^{n+1}$ has the same moments as $f^{n+1}$ and we assume that equations (\ref{approximate moments}) is accurate enough. This in turn gives
\begin{equation}
m[f^{n+1}_{i\cdot}]\simeq m[\tilde{f}^{n}_{i\cdot}].\label{Maxw n+1 approx Maxw tilde}
\end{equation}
Once the Maxwellian at time $t^{n+1}$ is known using the approximated macroscopic moments $m[\tilde{f}^{n}_{i\cdot}]$, the distribution function $f^{n+1}_{ij}$ can be explicitly computed
\begin{equation}
f^{n+1}_{ij}=\dfrac{\varepsilon\tilde{f}^{n}_{ij}+\Delta tM^{n+1}_{ij}}{\varepsilon+\Delta t}.\label{first order solution implicit step}
\end{equation}
This approach has already been used in \cite{RusSanYun}, \cite{santpre}, \cite{PhdSant} and in \cite{PP07} in the context of Eulerian schemes.

\section{High order methods}
\subsection{Runge-Kutta methods}
The scheme of the previous section corresponds to implicit Euler applied to the BGK model in characteristic form. High order discretization in time can be obtained by Runge-Kutta or BDF methods.\\
In \cite{santpre}, \cite{PhdSant}, the relaxation operator has been dealt with an L-stable diagonally implicit Runge-Kutta scheme \cite{HaiWa}. DIRK schemes are completely characterized by the triangular $\nu\times\nu$ matrix $A=(a_{lk})$, and the coefficient vectors $c=(1,\cdots,c_{\nu})^{T}$ and $b=(b_{1},\cdots,b_{\nu})^{T}$, which are derived by imposing accuracy and stability constraints \cite{HaiWa}.\\
 DIRK schemes can be represented through the Butcher's table
$$\begin{tabular}{c|c}
$c$ & $A$ \\\hline
 & $b^{T}$ \\
\end{tabular}.$$
Here we consider the following DIRK schemes
$$\text{RK2}=\begin{tabular}{c|cc}
$\alpha$ & $\alpha$ & 0\\
1 & 1-$\alpha$ & $\alpha$\\\hline
  & 1-$\alpha$ & $\alpha$
\end{tabular},\;\;\;\;\;
\text{RK3}=\begin{tabular}{c|ccc}
$\frac{1}{2}$ & $\gamma$ & 0 & 0\\
 $(1+\gamma)/2$ & $(1-\gamma)/2$ & $\gamma$ & 0\\
 1 & $1-\delta-\gamma$&$\delta$ &$\gamma$\\\hline
 &$1-\delta-\gamma$&$\delta$&$\gamma$\\
\end{tabular}
$$
which are a second and third order L-implicit schemes, respectively \cite{ARS97}. The coefficient $\alpha$ is  $$\alpha=1-\frac{\sqrt{2}}{2},$$ 
while $\gamma$ is the middle root of $6x^{3}-18x^{2}+9x-1,$ $\gamma\simeq 0.4358665215,$ and $\delta=3/2\gamma^{2}-5\gamma+5/4\simeq-0.644363171.$
Both RK schemes have the property that the last row of the matrix $A$ equals $b^{T},$ therefore the numerical solution is equal to the last stage value. Such schemes are called ``stiffly accurate''. An A-stable scheme which is stiffly accurate is also L-stable \cite{HaiWa}.
Applying the DIRK schemes to the characteristic formulation of the BGK equation (\ref{lagrBGK}), the numerical solution is obtained as
\begin{equation}
f^{n+1}_{ij}=f^{(\nu,n)}_{ij}+\Delta t\sum_{\ell=1}^{\nu}b_{\ell}K^{(\nu,\ell)}_{ij},\label{general form RK}
\end{equation}
where
$$K^{(\nu,\ell)}_{ij}=\dfrac{1}{\varepsilon}(M[F^{(\nu,\ell)}_{ij}]-F^{(\nu,\ell)}_{ij})$$
denote the RK fluxes on the characteristics $x=x_{i}+v_{j}(t-t^{n+1}),$ and
$$F^{(\nu,\ell)}_{ij}=f^{(\nu,n)}_{ij}+\Delta t\sum_{k=1}^{\ell}a_{\ell k}K^{(\ell,k)}_{ij}$$
are the stage values;  the first index of the pair  $(\nu,\ell)$  indicates that we are along the $\nu$-th characteristic and the second one denotes that we are computing the $\ell$-th stage value. Moreover $f^{(\nu,n)}_{ij}\equiv f(t^{n},x_{i}-c_{\nu}\Delta tv_{j},v_{j}).$

In a standard DIRK method, the $\ell$-th stage value, say $F_{ij}^{(\nu,\ell)},$ is evaluated by solving an implicit equation involving only $F_{ij}^{(\nu,\ell)},$ since the previous stage values have already been computed, due to the triangular structure of the matrix $A$. In our case this is not so easy, because if the point corresponding to stage $\ell$ along the characteristics is not a grid point, it is no possible to compute the moments of the Maxwellian at that point in space-time; indeed, after multiplying by $\phi_{j}$ and summing on $j$, the elements of the sum are computed in variable space points, so we cannot take advantage of the useful properties of the collision invariants. For this reason, we need two kinds of stage values: the stage value along the characteristics, $F_{ij}^{(\nu,\ell)}$, and the stage values on the grid, $F_{ij}^{(\ell,\ell)}$ (see Figure \ref{RK2_fig} and \ref{RK3_fig}).\\
Second and third order RK schemes are described below.
\subsubsection{RK2}
\begin{figure}[h]
{\includegraphics[trim=3cm 20.2cm 4cm 4cm,scale=1]{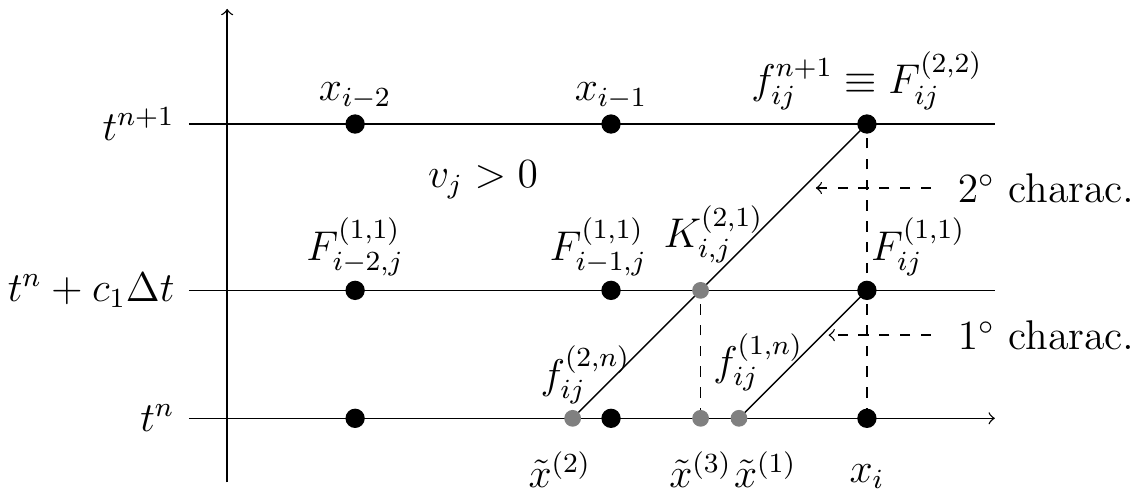}}
\caption{\footnotesize{Representation of the RK2 scheme. The black circles denote grid nodes, the gray ones the points where interpolation is needed.}}\label{RK2_fig}
\end{figure}

The general form of RK2 is (see Fig. \ref{RK2_fig})
\begin{equation}
f^{n+1}_{ij}=f^{2,n}_{ij}+\Delta t(b_{1}K^{(2,1)}_{ij}+b_{2}K^{(2,2)}_{ij}).\label{general form RK2}
\end{equation}
First we compute $F_{ij}^{(1,1)}$ in the grid node by $$F_{ij}^{(1,1)}=\frac{\varepsilon f_{ij}^{(1,n)}+\Delta ta_{11}M_{ij}^{(1,1)}}{\varepsilon+\Delta ta_{11}}.$$ The Maxwellian $M_{ij}^{(1,1)}=M[F_{ij}^{(1,1)}]$ can be evaluated  using the macroscopic moments $m[f_{i\cdot}^{(1,n)}],$ using an argument similar to the one adopted in (\ref{Maxw n+1 approx Maxw tilde}). $f_{i,j}^{(1,n)}=f(t^{n},x_{i}-a_{11} v_{j}\Delta t,v_{j})$ can be computed by a suitable WENO space reconstruction at time $t^{n}$ \cite{WENO}; in Appendix we report the WENO reconstructions adopted in the paper.

\noindent

Once the implicit step is solved, the Runge-Kutta fluxes $K_{ij}^{(2,1)} =\\ \dfrac{1}{\varepsilon}(M[F^{(1,1)}_{ij}] -
F^{(1,1)}_{ij})$ are computed by high order interpolation on the intermediate nodes $\tilde{x}^{(3)}$ along the characteristics. Then the second stage value can be computed by
\begin{equation}
F_{ij}^{(2,2)}=f_{ij}^{(2,n)}+\Delta t\bigg(a_{21}K_{ij}^{(2,1)}+a_{22}\frac{1}{\varepsilon}(M[F_{ij}^{(2,2)}]-F_{ij}^{(2,2)})\bigg).     \label{sv2}
\end{equation}
Equation (\ref{sv2}) cannot be immediately solved because the Maxwellian depends on $F_{ij}^{(2,2)}$ itself. However, if we take the moments of both sides of equation (\ref{sv2}), we can compute the moments of $F_{ij}^{(2,2)}$ since the elements of the sum on $j$ containing the Maxwellian $M[F_{ij}^{(2,2)}]$ are now on  fixed space points. Indeed
$$\sum_{j}\bigg(F_{ij}^{(2,2)}-f_{ij}^{(2,n)}-\Delta ta_{21}K_{ij}^{(2,1)}\bigg)\phi_{j}=a_{22}\frac{\Delta t}{\varepsilon}\sum_{j}\bigg(M[F_{ij}^{(2,2)}]-F_{ij}^{(2,2)}\bigg)\phi_{j}=0,$$
thus the moments are given by $m[F_{i\cdot}^{(2,2)}]=m[f_{i\cdot}^{(2,n)}+\Delta ta_{21}K_{ij}^{(2,1)}],$ so we can compute $M[F_{ij}^{(2,2)}],$ and solve the implicit step for $F_{ij}^{(2,2)}.$\\
Notice that $f_{ij}^{n+1}=F_{ij}^{(2,2)},$ because the scheme is stiffly accurate, i.e the last row of the matrix $A$ is equal to the vector of weights.
\subsubsection{RK3}
\begin{figure}[h]
{\includegraphics[trim=4.5cm 16.35cm 2.5cm 4cm,scale=1]{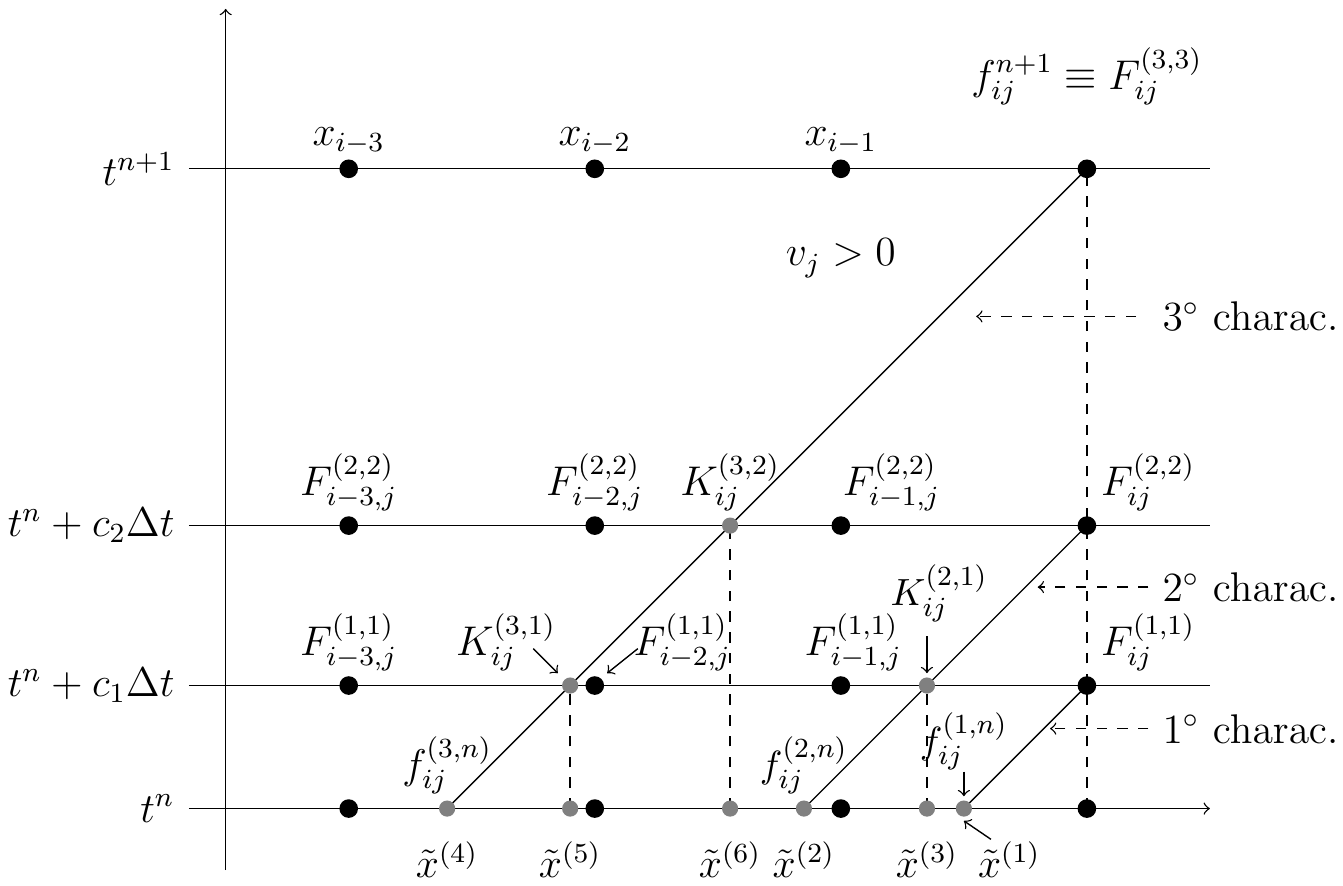}}
\caption{\footnotesize{Representation of the RK3 scheme. The black circles denote grid nodes, the gray ones the points where interpolation is needed.}}\label{RK3_fig}
\end{figure}

The RK3 scheme works in a similar way, and Fig. \ref{RK3_fig} shows procedure. \smallskip

\textbf{Algorithm (RK3)}
\begin{itemize}
\item[-] Calculate $f^{(1,n)}_{i,j}=f(t^{n},\tilde{x}^{(1)}=x_{i}-c_{1}v_{j}\Delta t,v_{j}),$ $f^{(2,n)}_{i,j}=f(t^{n},\tilde{x}^{(2)}=x_{i}-c_{2}v_{j}\Delta t,v_{j}),$
$\tilde{f}^{3,n}_{i,j}=f(t^{n},\tilde{x}^{(4)}=x_{i}-v_{j}\Delta t,v_{j})$ by interpolation from $f^{n}_{\cdot j}$;
\item[-] Calculate $F^{(1,1)}_{ij}$ in the grid node using the technique (\ref{Maxw n+1 approx Maxw tilde}), (\ref{first order solution implicit step}), with $\Delta t$ replaced by $c_{1}\Delta t$. Given $F^{(1,1)}_{ij}$, one can evaluate the Runge-Kutta fluxes  $K^{(1,1)}_{ij}=\dfrac{1}{\varepsilon}\bigg(M[F^{(1,1)}_{ij}]-F^{(1,1)}_{ij}\bigg)$ in the grid nodes and then calculate $K^{(2,1)}_{ij}$ and $K^{(3,1)}_{ij}$ by interpolation from $K^{(1,1)}_{\cdot j}$ in $\tilde{x}^{(3)}=x_{i}-(c_{2}-c_{1})v_{j}\Delta t$ and $\tilde{x}^{(5)}=x_{i}-(1-c_{1})v_{j}\Delta t$, respectively;
\item[-] Calculate $F^{(2,2)}_{ij}$ in the grid node using RK2 scheme described in the previous section with time step $c_{2}\Delta t$. Given $F^{(2,2)}_{ij}$, one can evaluate $K^{(2,2)}_{ij}=\dfrac{1}{\varepsilon}\bigg(M[F^{(2,2)}_{ij}]-F^{(2,2)}_{ij}\bigg)$ in the grid nodes and then calculate $K^{(3,2)}_{ij}$  by interpolation from $K^{(2,2)}_{\cdot j}$ in $\tilde{x}^{(6)}=x_{i}-(1-c_{2})v_{j}\Delta t$;
\item[-] Now one can update $f^{n+1}_{ij}$ using (\ref{general form RK}), taking into account that the method is stiffly accurate and using the properties of the collision invariants to solve the implicit step.
\end{itemize}

\subsubsection{Summary of the Runge-Kutta schemes}
Three schemes based on RK are tested in the paper:
\begin{itemize}
\item[-]scheme RK2W23: uses WENO23 for the interpolation and RK2, as described above, for time integration;
\item[-]scheme RK3W23: uses WENO23 for the interpolation and RK3, as described above, for time integration;
\item[-]scheme RK3W35: uses WENO35 for interpolation and RK3 for time integration.
\end{itemize}

\textbf{Remark}\medskip \\
In practice, the Runge-Kutta fluxes can be computed from the internal stages. For example, using RK2, we have
$$K_{ij}^{(1,1)}\dfrac{1}{\varepsilon}(M[F^{(1,1)}_{ij}]-F^{(1,1)}_{ij})=\dfrac{F^{(1,1)}_{ij}-\tilde{f}^{2,n}_{ij}}{\Delta ta_{11}}.$$
The latter expression can be used in the limit $\varepsilon\rightarrow 0,$ with no constraint on the time step.

\subsection{BDF methods}
In this section we present a new family of high order semilagrangian schemes, based on BDF. The backward differentiation formula are implicit linear multistep methods for the numerical integration of ordinary differential equations $y' = g(t, y)$ \cite{HaiWa}.
Using the linear polynomial interpolating $y^{n}$ and $y^{n-1}$ one obtains the simplest BDF method (BDF1) that correspond to backward Euler, used  in Section 1.

Here the characteristic formulation of the BGK model, that leads to ordinary differential equations, is approximated  by using  BDF2 and BDF3 methods, in order to obtain high order approximation. The relevant expressions, under the hypothesis that the time step $\Delta t$ is fixed, are:
\begin{equation}
\text{BDF2}:=\;\;\;\;y^{n+1}=\dfrac{4}{3}y^{n}-\dfrac{1}{3}y^{n-1}+\dfrac{2}{3}\Delta t\,g(y^{n+1},t^{n+1}),\hspace{1.8cm}\label{BDF2}
\end{equation}
\begin{equation}
\text{BDF3}:=\;\;\;\;y^{n+1}=\dfrac{18}{11}y^{n}-\dfrac{9}{11}y^{n-1}+\dfrac{2}{11}y^{n-2}+\dfrac{6}{11}\Delta t\,g(y^{n+1},t^{n}).\label{BDF3}
\end{equation}
Here we apply the BDF methods along the characteristics.
\subsubsection{BDF2}
The numerical approximation of the first equation in (\ref{lagrBGK}) is obtained as
\begin{equation}
\text{BDF2}:=\;\;\;\;f^{n+1}_{i,j}=\dfrac{4}{3}f^{n,1}_{ij}-\dfrac{1}{3}f^{n-1,2}_{ij}+\dfrac{2}{3}\frac{\Delta t}{\epsilon}(M[f]^{n+1}_{ij}-f^{n+1}_{ij}),\label{BDF2 approx}
\end{equation}
where $f^{n-(s-1),s}_{i,j}\simeq f(t^{n-(s-1)},x_{i}-sv_{j}\Delta t,v_{j})$, can be computed by suitable reconstruction from $\{f^{n-(s-1)}_{\cdot j}\}$; high order reconstruction will be needed for BDF2 and BDF3 schemes, and again we make use of  WENO techniques \cite{WENO} for accurate non oscillatory reconstruction.

\begin{figure}[h]
{\includegraphics[trim=3.5cm 19cm 2cm 3cm,scale=1]{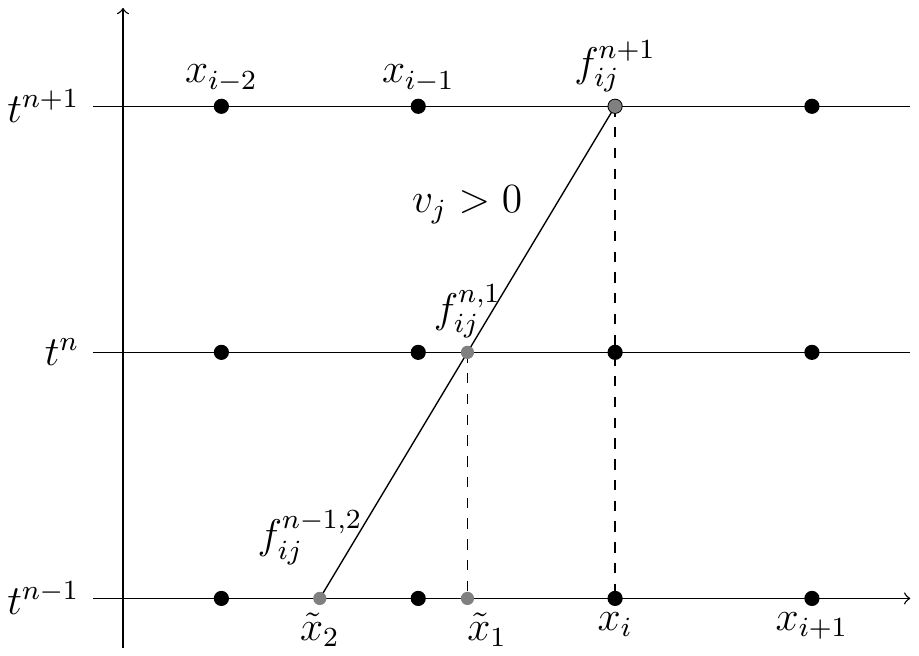}}
\caption{\footnotesize{Representation of the BDF2 scheme. The black circles denote grid nodes, the gray ones the points where interpolation is needed.}}\label{BDF2_fig}
\end{figure}

To compute the solution $f^{n+1}_{ij}$ from equations (\ref{BDF2 approx}), also in this case one has to solve a non linear implicit equation. We can act as previously  done for the backward Euler method, by taking advantage of the properties of the collision invariants. Thus we multiply both sides of the equation (\ref{BDF2 approx}) by $\phi_{j}$ and sum over $j$, getting
$$\sum_{j}\bigg(f^{n+1}_{ij}-\frac{4}{3}f^{n,1}_{ij}+\frac{1}{3}f^{n-1,2}_{ij}\bigg)\phi_{j}=\frac{2\Delta t}{3\varepsilon}\sum_{j}(M[f]^{n+1}_{ij}-f^{n+1}_{ij})\phi_{j},$$\\
which implies that
$$\sum_{j}(f^{n+1}_{ij})\phi_{j}=\sum_{j}\bigg(\frac{4}{3}f^{n,1}_{ij}-\frac{1}{3}f^{n-1,2}_{ij}\bigg)\phi_{j},$$
so in Equation (\ref{BDF2 approx}) we can compute $M[f^{n+1}_{ij}]$ with the usual procedure adopting the approximated macroscopic moments
\begin{equation}
(\rho_{i}^{n+1},(\rho u)^{n+1}_{i},E^{n+1}_{i})=m\bigg[\frac{4}{3}f^{n,1}_{i\cdot}-\frac{1}{3}f^{n-1,2}_{i\cdot}\bigg].\label{approx moments for BDF2}
\end{equation}
Once the Maxwellian $M[f^{n+1}_{ij}]$ is computed, the distribution function value $f^{n+1}_{ij}$ can be easily obtained from schemes (\ref{BDF2 approx}) for BDF2. The procedure for BDF2 is sketched in Fig. \ref{BDF2_fig} and described the following algorithm.

\textbf{Algorithm (BDF2)}\\
\begin{itemize}
\item[-]  Calculate $f^{n-1,2}_{ij}=f(t^{n-1},\tilde{x}_{2}=x_{i}-2v_{j}\Delta t,v_{j})$, $f^{n,1}_{ij}=f(t^{n},\tilde{x}_{1}=x_{i}-v_{j}\Delta t,v_{j})$ by interpolation from $f^{n-1}_{\cdot j}$ and $f^{n}_{\cdot j}$ respectively;
\item[-] Compute the Maxwellian $M[f^{n+1}_{ij}]$ using  (\ref{approx moments for BDF2}) and upgrade the numerical solution $f_{ij}^{n+1}$.
\end{itemize}

\noindent
A similar algorithm is obtained using BDF3, as we will see in the next subsection.

\subsubsection{BDF3}

The numerical solution of the BGK equation in (\ref{lagrBGK}) is obtained as
\begin{equation}
f^{n+1}_{i,j}=\dfrac{18}{11}f^{n,1}_{ij}-\dfrac{9}{11}f^{n-1,2}_{ij}+\dfrac{2}{11}f^{n-2,3}_{ij}+\dfrac{6}{11}\frac{\Delta t}{\epsilon}(M[f]^{n+1}_{ij}-f^{n+1}_{ij}),\label{BDF3 approx}
\end{equation}
where $f^{n-(s-1),s}_{i,j}$ can be computed by suitable reconstruction from $\{f^{n-(s-1)}_{\cdot j}\}$.

\begin{figure}[h]
{\includegraphics[trim=3.5cm 18.5cm 3.5cm 4.5cm,scale=1]{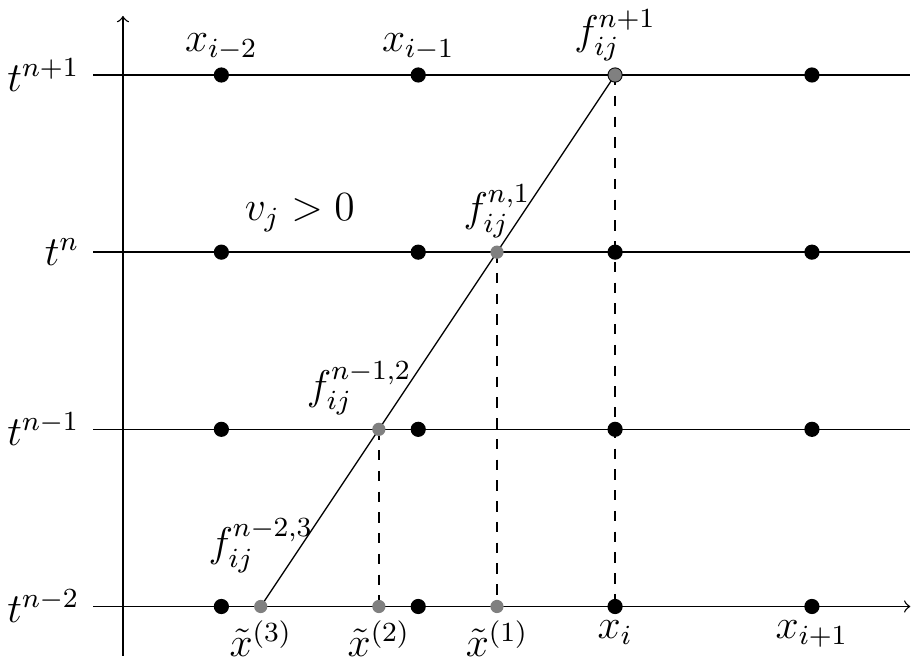}}
\caption{\footnotesize{Representation of the BDF3 scheme. The black circles denote grid nodes, the gray ones the points where interpolation is needed.}}\label{BDF3_fig}
\end{figure}

To compute the solution $f^{n+1}_{ij}$ from equation (\ref{BDF3 approx}) we need again to take moments of such equation
$$\sum_{j}\bigg(f^{n+1}_{ij}-\frac{18}{11}f^{n,1}_{ij}+\frac{9}{11}f^{n-1,2}_{ij}-\frac{2}{11}f^{n-2,3}_{ij}\bigg)\phi_{j}=\frac{6\Delta t}{11\varepsilon}\sum_{j}(M[f]^{n+1}_{ij}-f^{n+1}_{ij})\phi_{j},$$
which implies that
$$\sum_{j}(f^{n+1}_{ij})\phi_{j}=\sum_{j}\bigg(\frac{18}{11}f^{n,1}_{ij}-\frac{9}{11}f^{n-1,2}_{ij}+\frac{2}{11}f^{n-2,3}_{ij}\bigg)\phi_{j},$$
so in Equation (\ref{BDF3 approx}) we can compute $M[f^{n+1}_{ij}]$ with the usual procedure, adopting the approximated macroscopic moments $$(\rho_{i}^{n+1},(\rho u)^{n+1}_{i},E^{n+1}_{i})=m\bigg[\frac{18}{11}f^{n,1}_{i\cdot}-\frac{9}{11}f^{n-1,2}_{i\cdot}+\frac{2}{11}f^{n-2,3}_{i\cdot}\bigg]$$
Once the Maxwellian $M[f^{n+1}_{ij}]$ is computed, the distribution function value $f^{n+1}_{ij}$ can be easily obtained from schemes (\ref{BDF3 approx}) for BDF3. This procedure is sketched in Fig. \ref{BDF3_fig}.

To compute the starting values $f^{1}_{ij}$ for  BDF2  and $f^{1}_{ij},\;f^{2}_{ij}$ for  BDF3  we have used, as predictor, Runge Kutta methods of order 2 and 3, respectively.\\

\subsubsection{Summary of the BDF schemes}
Three schemes based on BDF are tested in the paper:
\begin{itemize}
\item[-]scheme BDF2W23: uses WENO23 for the interpolation and BDF2, as described above, for time integration;
\item[-]scheme BDF3W23: uses WENO23 for the interpolation and BDF2, as described above, for time integration;
\item[-]scheme BDF3W35: uses WENO35 for interpolation and BDF3 for time integration.
\end{itemize}

At variance with Runge-Kutta methods, BDF methods do not need to compute intermediate stage values, and the implicit step for the Maxwellian is solved only once during a time step. Moreover, we have to interpolate in less  out-of-grid points (for instance, we have to perform only 3 interpolations in a time step using BDF3, versus 6 interpolations needed to advance one time step using a DIRK method of order 3). This makes BDF methods very efficient from a computational point of view.

\section{Semi-lagrangian schemes without interpolation}
As we can observe, the cost of the schemes presented above, is mainly due  to the interpolation, especially when we  use  high order interpolation techniques.\\
In order to reduce the  computational cost we look for schemes that avoid interpolation. The key idea is to choose a discretization parameters in such a way that all the characteristics connect grid points in space. This is obtained, for example, by choosing $\Delta v\Delta t=\Delta x$ (See Fig. \ref{fig_without_int}).

\begin{figure}[h]
{\includegraphics[trim=1cm 21.1cm 2cm 6.7cm,scale=1]{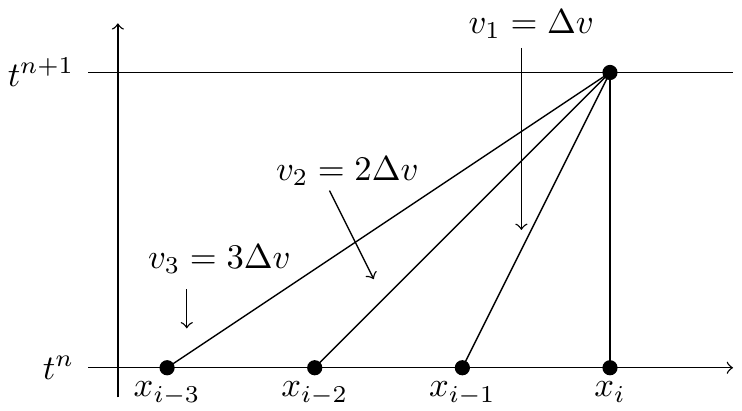}}
\caption{\footnotesize{Implicit first order scheme without interpolation.}} \label{without_interpolation_fig}
\end{figure}

This choice corresponds to solving the equation at each characteristics by implicit Euler, thus resulting in a first order method in time. In order to increase the order of accuracy one can resort to BDF or RK time discretization. BDF2 and BDF3 can be easily applied in this setting.\\ The use of higher order RK schemes requires that the stage value lie on the grid as well. This is obtained by imposing $\Delta v\Delta t=s\Delta x,$ $s\in \mathbb{N}.$ In this case the coefficients of the vector $\underline{c}$ must be multiples of  $1/s$. Moreover we need a L-stable scheme. Imposing accuracy and stability constraints on the coefficients of the Butcher's table together with the fact that coefficients of the vector $\underline{c}$ must be multiples of  $1/s$, we obtain some DIRK methods of second and third order, using respectively $s=3$ and $s=4.$\\ Using schemes that avoid the interpolation, the choice of the time step is determined by the other discretization steps, the CFL number is fixed to $sN_{v}.$ This means that we have very large time step. A Runge-Kutta method that avoids interpolation is the following
$$\begin{tabular}{c|cc}
1/3 & 1/3 & 0\\
1 & 3/4 & 1/4\\\hline
  & 3/4 & 1/4.
\end{tabular}$$
This is a second order method, diagonally implicit and L-stable, because it is stiffly-accurate and A-stable, and allows us to avoid interpolation using $s=3.$\\
These schemes are much simpler to implement and therefore each time step can be advanced very efficiently. However they require a very fine grid in space. A comparison with more standard semilagrangian methods that make use of interpolation will be presented im the section on numerical results.

\section{Chu reduction model}
The methods have been extended  to treat problem in 3D in velocity, 1D in space, in slab geometry. The technique used is the Chu reduction \cite{Chu}, which, under suitable symmetry assumption, allows to transform a 3D equation (in velocity) in a system of two equations 1D (in velocity), where the schemes  previously introduced can be applied.\\
We consider the application of BGK equation to problems with axial symmetry with respect to an axis (say, $x_{1}\equiv x$), in the sense that all transverse spatial gradients vanish, and the gas is drifting only in the axial direction. In such cases, distribution functions $f(t,x,\vec{v})$ depend on the full velocity vector $\vec{v}$ (i.e molecular trajectories are three-dimensional) but dependence on the azimuthal direction around the symmetry axis is such that all transverse components of the macroscopic velocity $\vec{u}$ vanish (i.e. $u_{2}=u_{3}=0$). \\
Let us introduce the new unknowns
\begin{equation}
g_{1}(t,x,v)=\int_{\mathbb{R}^{2}}f(t,x,(\vec{v}))\,dv_{2}dv_{3},\quad g_{2}(t,x,v)=\int_{\mathbb{R}^{2}}(v_{2}^{2}+v_{3}^{2})f(t,x,(\vec{v}))\,dv_{2}dv_{3},\label{marginal function}
\end{equation}
each depending only on one space and one velocity variable. Multiplication of (\ref{BGK}) by $1$ and $(v_{2}^{2}+v_{3}^{2})$ and integration with respect to $(v_{2},v_{3})\in\mathbb{R}^{2}$ yields then the following system of BGK equations for the unknown vector $g=(g_{1},g_{2})$, coupled with initial conditions
\begin{equation}
\dfrac{\partial g_{i}}{\partial t}+v\dfrac{\partial g_{i}}{\partial x}=\dfrac{1}{\varepsilon}(M[f]_{i}-g_{i}), \;\;\;(t,x,v)\in\mathbb{R}_{+}\times\mathbb{R}\times\mathbb{R}, \label{BGK Chu}\
\end{equation}
$$g_{i}(0,x,v)=g_{i,0}(x,v),\;\;\; i=1,2.$$
The BGK system (\ref{BGK Chu}) describes a relaxation process towards the vector function $(M[f]_{1},M[f]_{2})$, which is obtained by Chu transform of (\ref{Maxwellian}) with $N=3$ and has the form
$$(M[f]_{1},M[f]_{2})=(M[f]_{1},2RT\,M[f]_{1}),$$
where
$$M[f]_{1}=\dfrac{\rho(t,x)}{\sqrt{(2\pi RT(t,x)}}\exp\bigg(-\dfrac{(v-u(t,x))^{2}}{2RT(t,x)}\bigg).$$
The macroscopic moments of the distribution function $f$, needed to evaluate $M[f]_{1}$ are given in terms of $g_{1}$ and $g_{2}$ as:
$$\rho=\int_{\mathbb{R}}g_{1}\,dv_{1},\quad\quad\quad u=\dfrac{1}{\rho}\int_{\mathbb{R}}vg_{1}\,dv_{1},$$
$$3RT=\dfrac{1}{\rho}\bigg[\int_{\mathbb{R}}(v-u)^{2}g_{1}\,dv_{1}+\int_{\mathbb{R}}g_{2}\,dv_{1}\bigg].$$
The following relation will be useful to solve the implicit step:
\begin{equation}
\int_{\mathbb{R}}(v_{1}-u_{1})^{2}(M[f]_{1}-g_{1})\,dv_{1} + \int_{\mathbb{R}} (M[f]_{2}-g_{2})\,dv_{1}=0.\label{chu_implicit_step}
\end{equation}
Indeed
$$3RT\rho = \int_{\mathbb{R}}(v_{1}-u_{1})^{2}M[f]_{1}\,dv_{1} + 2RT\int_{\mathbb{R}} M[f]_{1}\,dv_{1},$$
and
$$3RT\rho = \int_{\mathbb{R}^{3}}f((v_{1}-u_{1})^{2}+v_{2}^{2}+v_{3}^{2})\,dv_{1}dv_{2}dv_{3}$$
$$= \int_{\mathbb{R}}(v_{1}-u_{1})^{2}g_{1}\,dv_{1} + \int_{\mathbb{R}} g_{2}\,dv_{1}.$$
Taking the difference we obtain (\ref{chu_implicit_step}).\\
The discrete version of the first order implicit scheme (in a similar way one can extend to high order schemes) of (\ref{BGK Chu}) is
\begin{equation}
g_{s,ij}^{n+1}=\tilde{g}_{s,ij}^{n}+\dfrac{\Delta t}{\varepsilon}(M_{s,ij}^{n+1}-g_{s,ij}^{n+1})\,\,\,\,s=1,2.\label{BGK_Chu_discr}
\end{equation}
To solve the implicit step we have to compute $m[g_{1,i\cdot}^{n+1}].$ The density $\rho_{i}^{n+1}$ and the momentum $(\rho u)_{i}^{n+1}$ can be easily computed multiplying the first equation of (\ref{BGK_Chu_discr}) by $1$ and $v_{j}$ and summing over $j.$ In this way we get
$$\rho_{i}^{n+1}=\Delta v\sum_{j}\tilde{g}_{1,ij}^{n},\quad\quad\quad(\rho u)_{i}^{n+1}=\Delta v\sum_{j}v_{j}\tilde{g}_{1,ij}^{n}.$$
To obtain the temperature $T_{i}^{n+1},$ instead we have to multiply by $(v_{1,j}-u_{1,j})^{2}$ and by $1$ respectively the first and the second equation of (\ref{BGK_Chu_discr}), and than summing over $j.$\\
Now, using the discrete analogue of (\ref{chu_implicit_step}):
$$\Delta v\sum_{j}(v_{1,j}-u_{1,j})^{2}(M_{1,ij}^{n+1}-g_{1,ij}^{n+1})+\Delta v\sum_{j}(M_{2,ij}^{n+1}-g_{1,ij}^{n+1})=0,$$
one can compute the temperature $T_{i}^{n+1}$ in this way:
$$3R\rho_{i}^{n+1}T_{i}^{n+1}=\Delta v\sum_{j}(v_{j}-u_{j})^{2}\tilde{g}_{1,ij}^{n}+\Delta v\sum_{j}\tilde{g}_{2,ij}^{n}.$$
Once the new moments $\rho_{i}^{n+1},$ $(\rho u)_{i}^{n+1}$ and $T_{i}^{n+1},$ are computed, we can solve the implicit step and to upgrade the numerical solution.

\section{Numerical tests}

We have considered two types of numerical tests with the purpose of verifying the accuracy (test 1) and the shock capturing properties (test 2) of the schemes. Different values of the Knudsen number have been investigated in order to observe the behavior of the methods varying from the rarefied ( $\varepsilon\simeq 1$) to the fluid ($\varepsilon\simeq 10^{-6}$) regime. We use units for temperature such that $R=1.$\\
In the first part of the section we consider the single 1D model and we explore the choice of the optimal CFL. A comparison between semilagrangian schemes with and without interpolation is also presented. The second part of the section is devoted to the results on the method applied to the 1D space-3D velocity case in slab geometry (Chu reduction).

\subsection{Regular velocity perturbation}
\begin{figure}[h]
\subfigure
{\includegraphics[trim=3cm 7cm 2cm 8.5cm,scale=0.4]{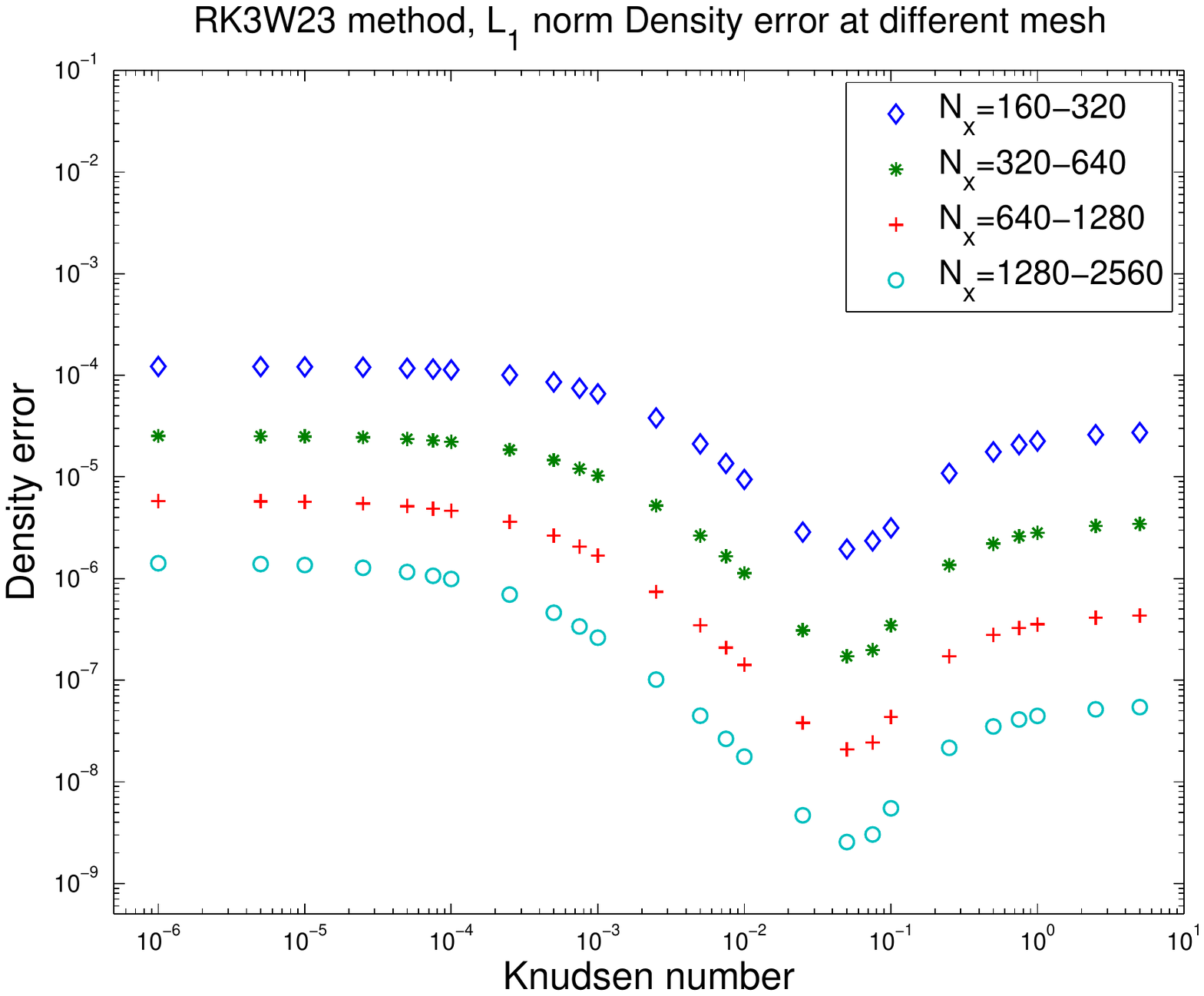}}
\subfigure
{\includegraphics[trim=1cm 7cm 2cm 8.5cm,scale=0.4]{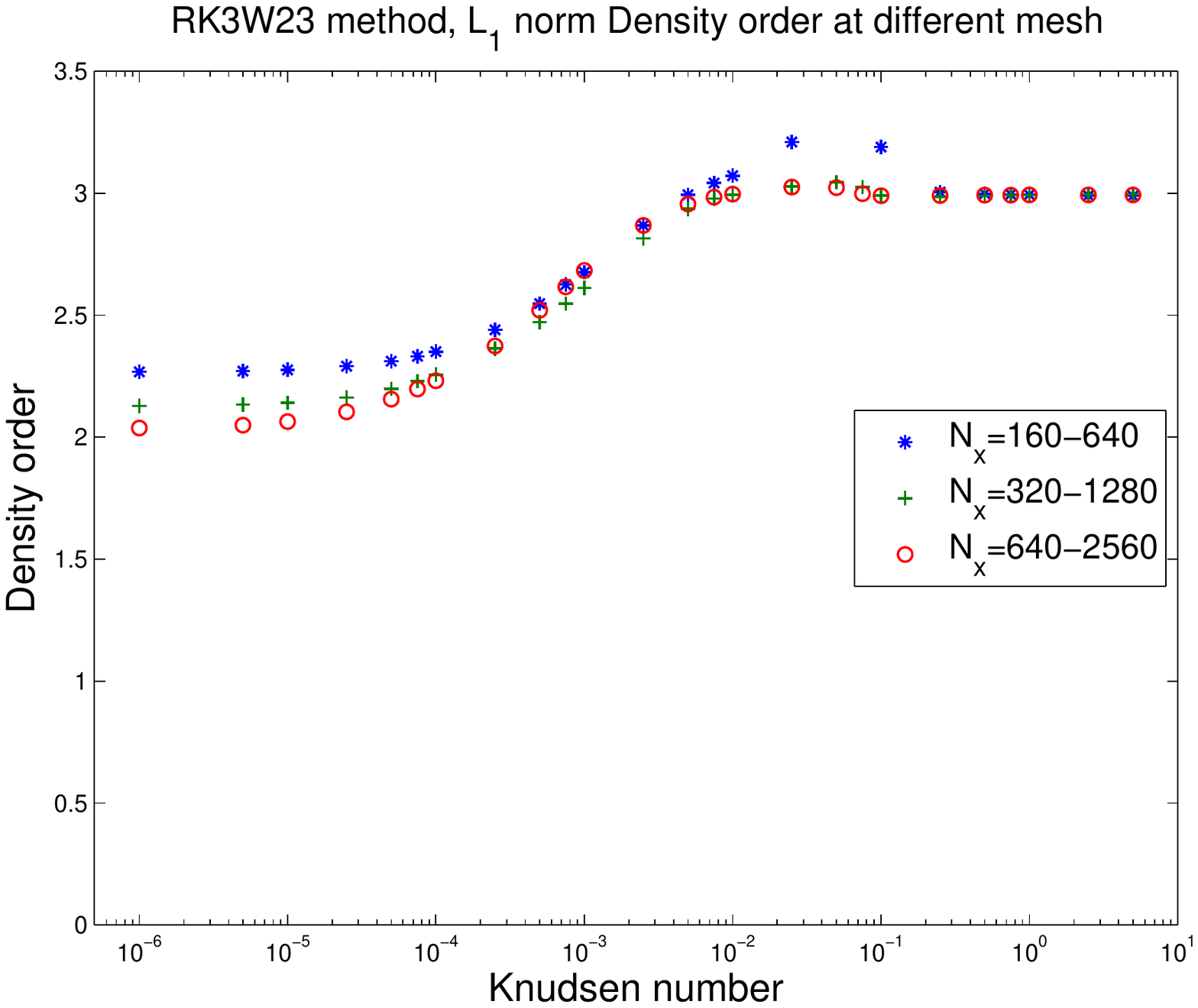}}
\subfigure
{\includegraphics[trim=3cm 7cm 2cm 7cm,scale=0.4]{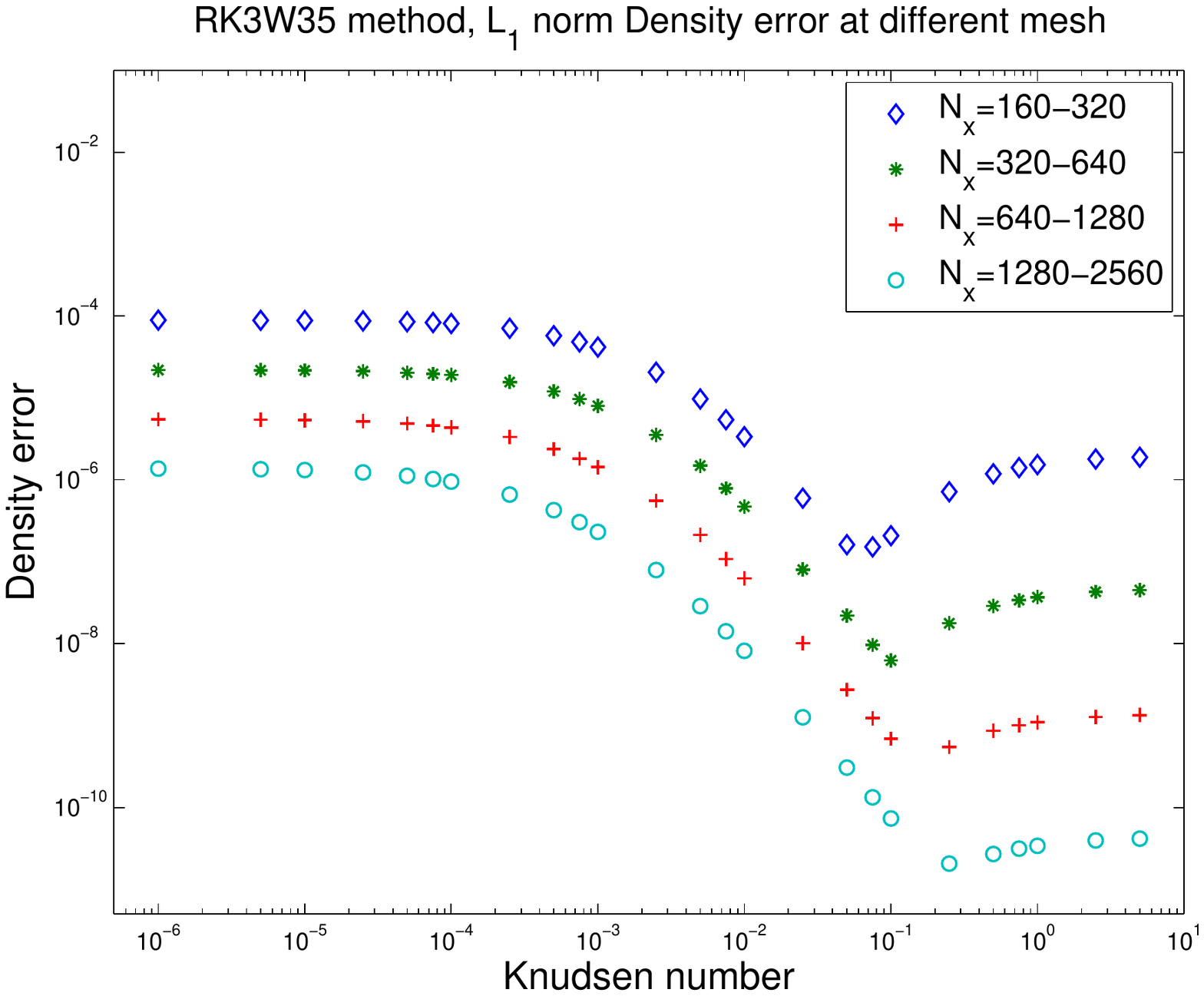}}
\subfigure
{\includegraphics[trim=1cm 7cm 2cm 7cm,scale=0.4]{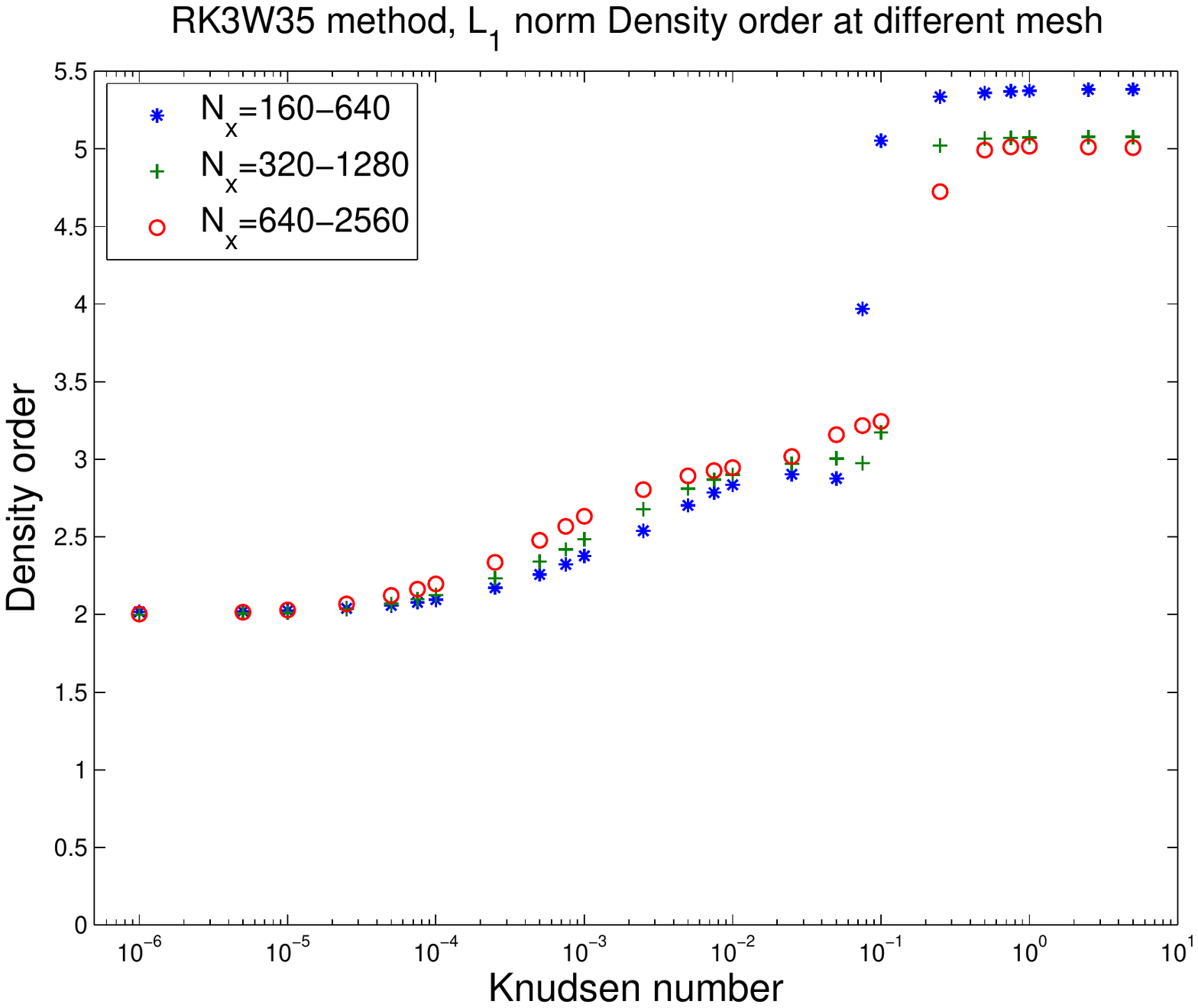}}
\caption{\footnotesize{$L_{1}$ error and accuracy order of RK3W23 and RK3W35, varying $\varepsilon$, using periodic boundary condition.}}\label{fig_RK3W23_order_periodic}
\end{figure}

\begin{figure}[h]
\subfigure
{\includegraphics[trim=3cm 7cm 2cm 8.5cm,scale=0.4]{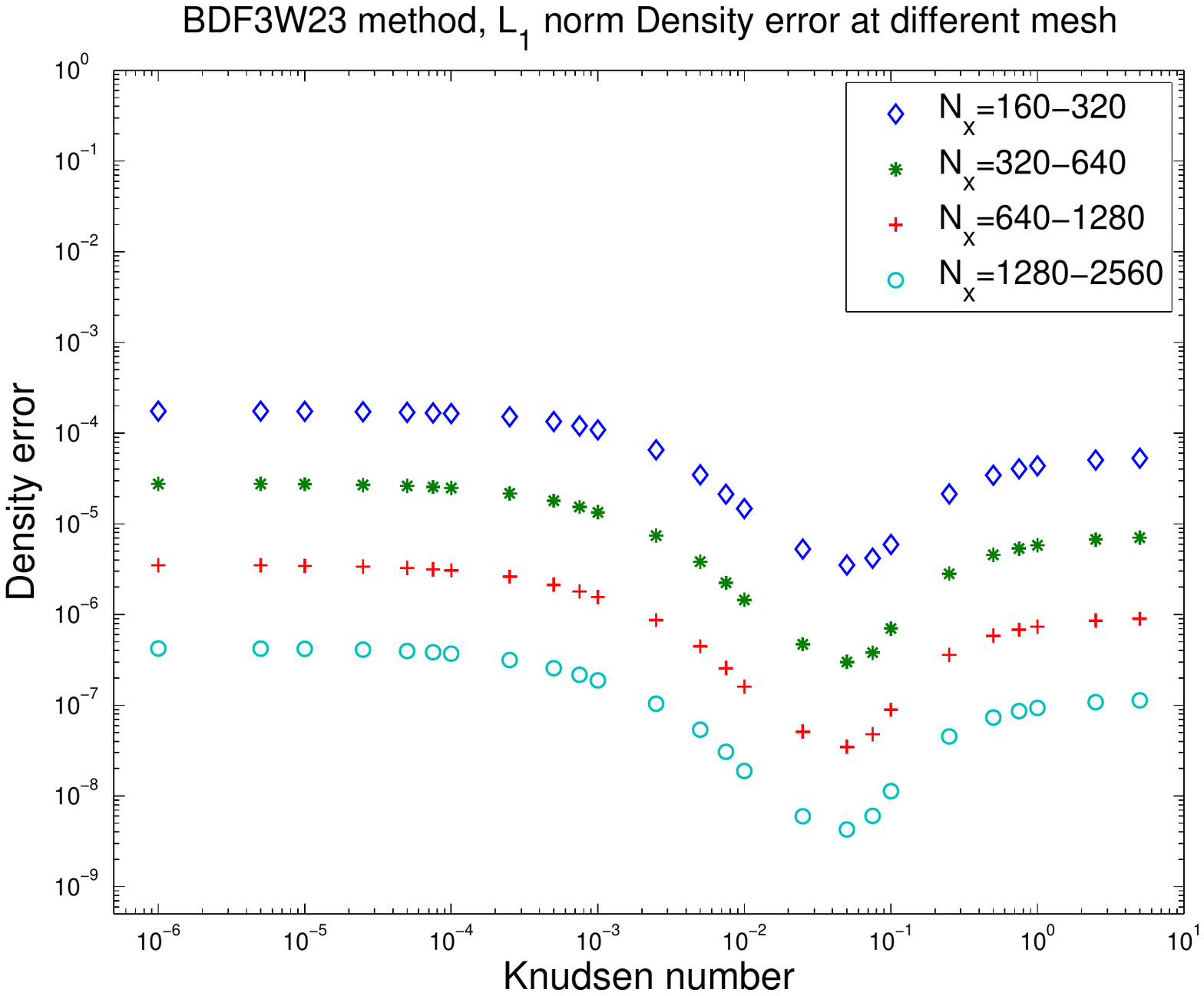}}
\subfigure
{\includegraphics[trim=1cm 7cm 2cm 8.5cm,scale=0.4]{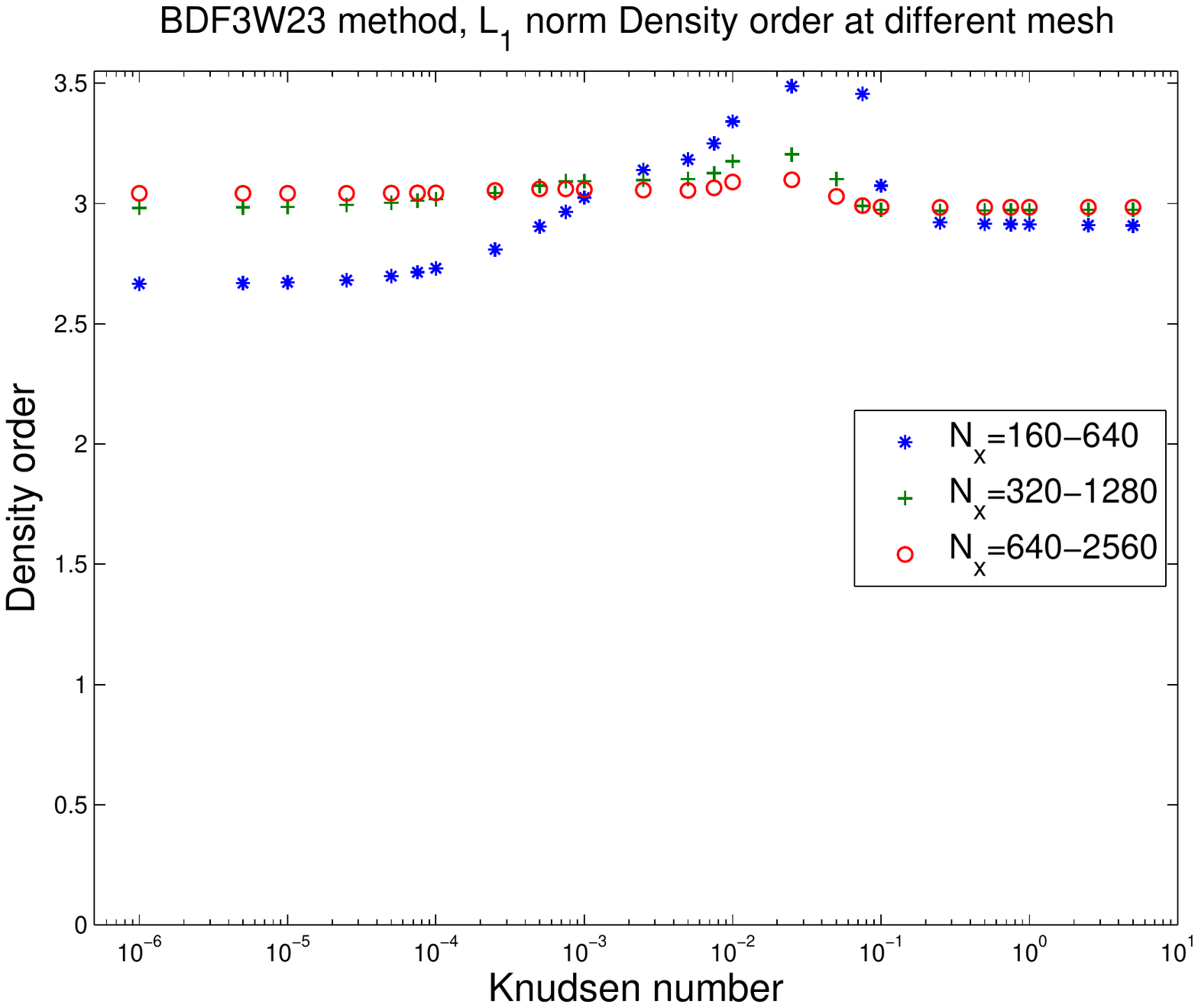}}
\subfigure
{\includegraphics[trim=3cm 7cm 2cm 7cm,scale=0.4]{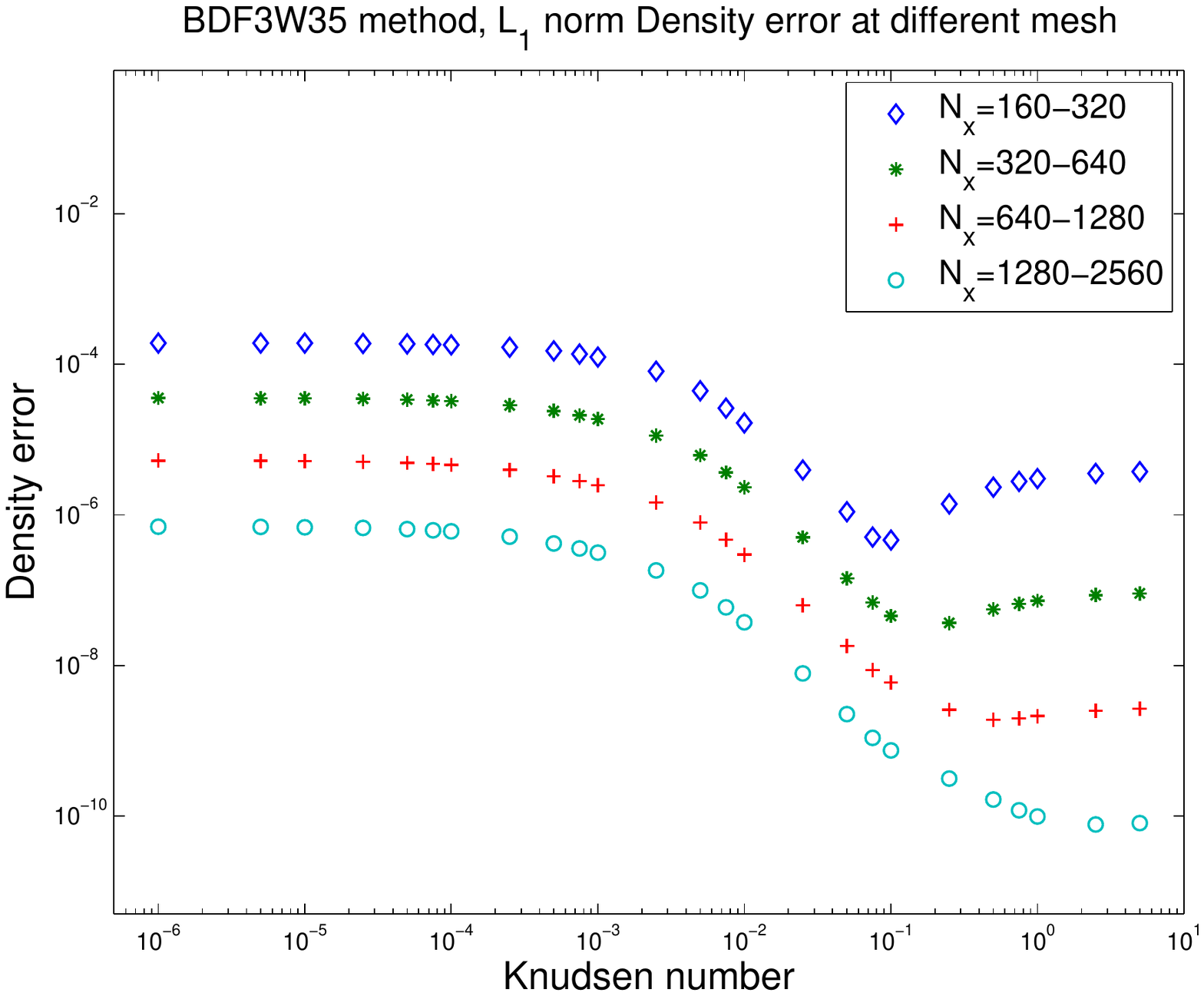}}
\subfigure
{\includegraphics[trim=1cm 7cm 2cm 7cm,scale=0.4]{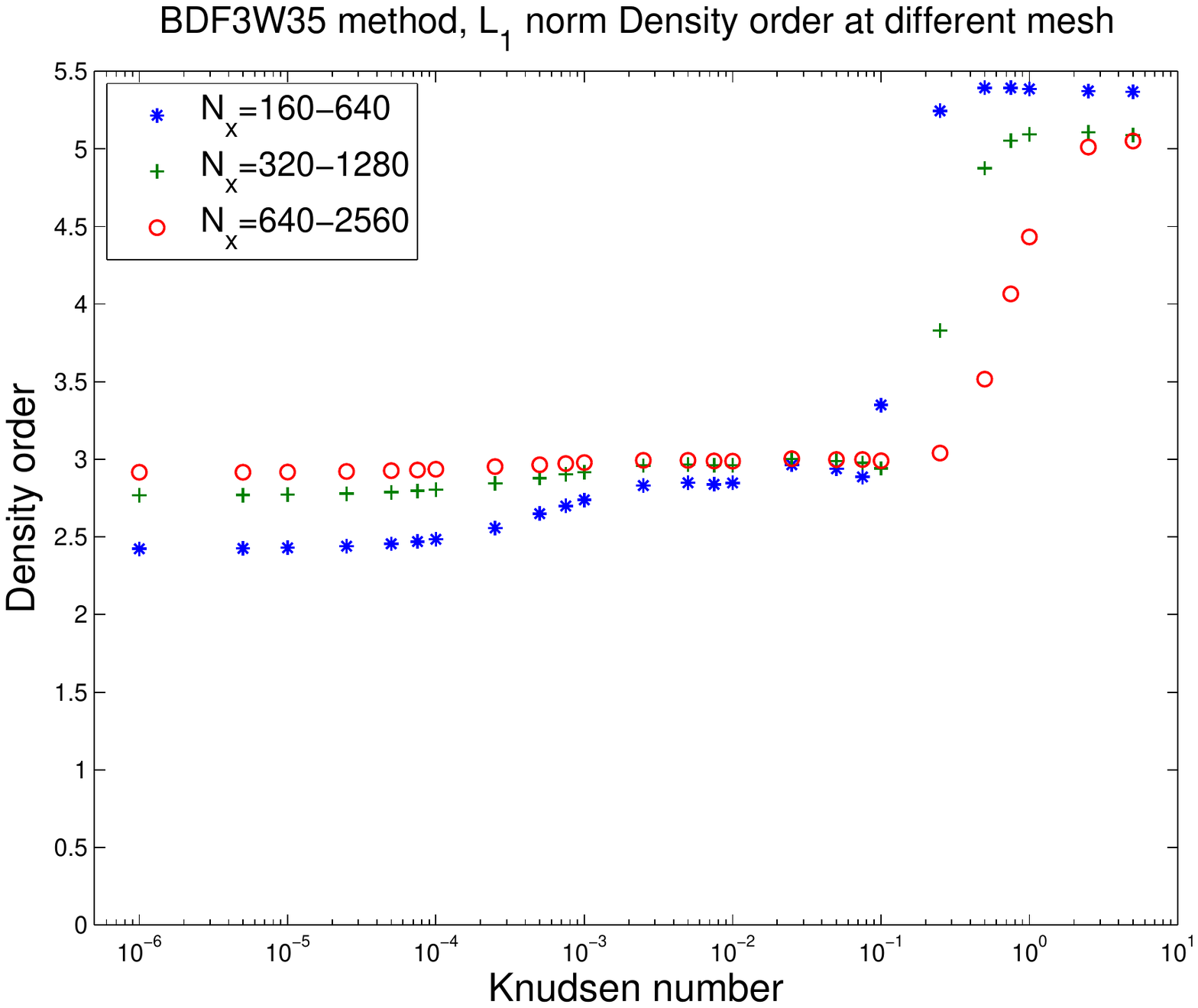}}
\caption{\footnotesize{$L_{1}$ error and accuracy order of BDF3W23 and BDF3W35, varying $\varepsilon$, using periodic boundary condition.}}\label{fig_BDF3W35_order_periodic}
\end{figure}

This test has been proposed in \cite{PP07}. Initial velocity profile is given by $$u=0.1\exp(-(10x-1)^{2})-2\exp(-(10x+3)^{2}), \;\;\;x\in[-1,1].$$
Initial density and temperature profiles are uniform, with constant value, $\rho=1$ and $T=1$.
\begin{figure}[h]
\subfigure
{\includegraphics[trim=3cm 7cm 2cm 8.5cm,scale=0.4]{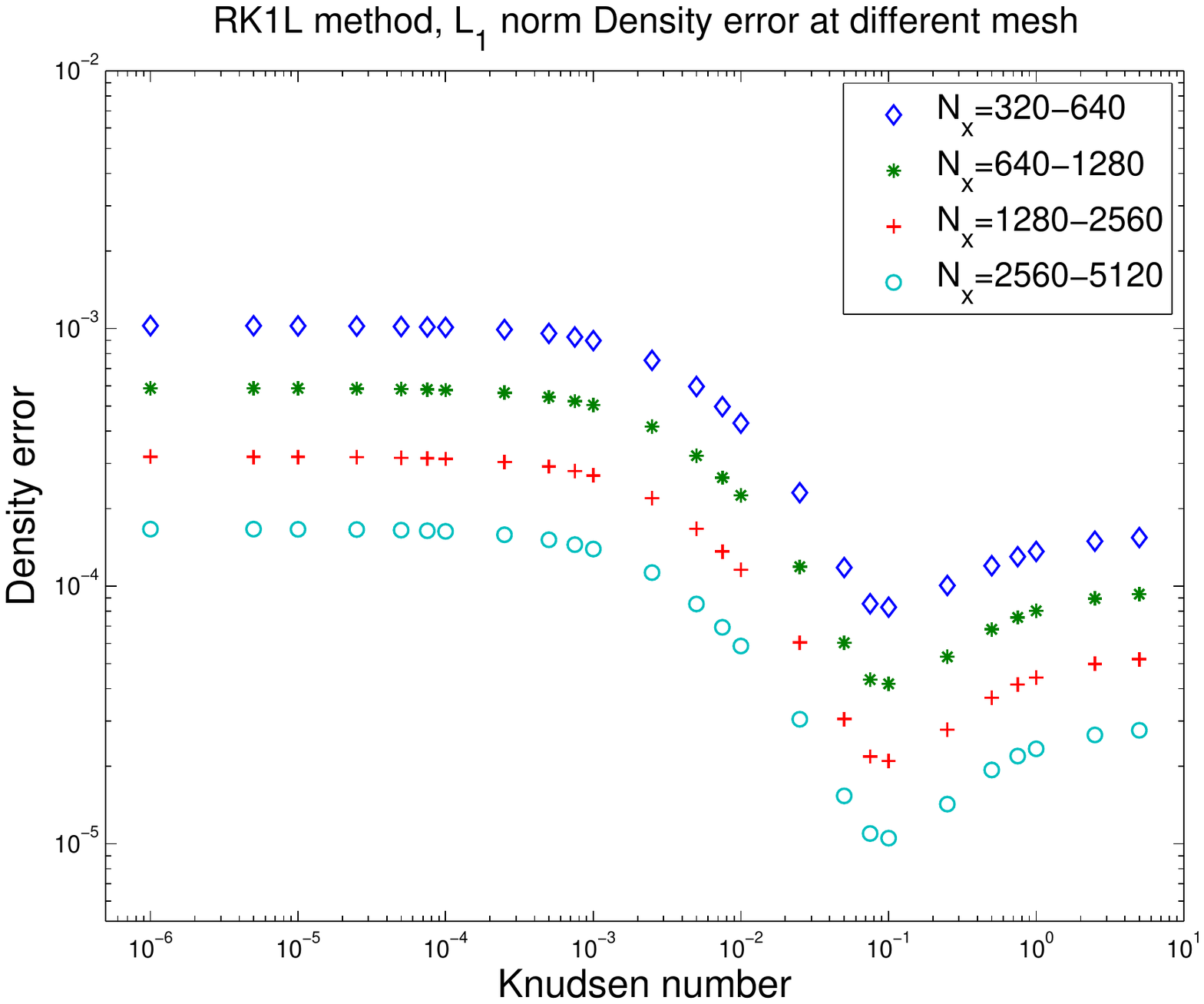}}
\subfigure
{\includegraphics[trim=1cm 7cm 2cm 8.5cm,scale=0.4]{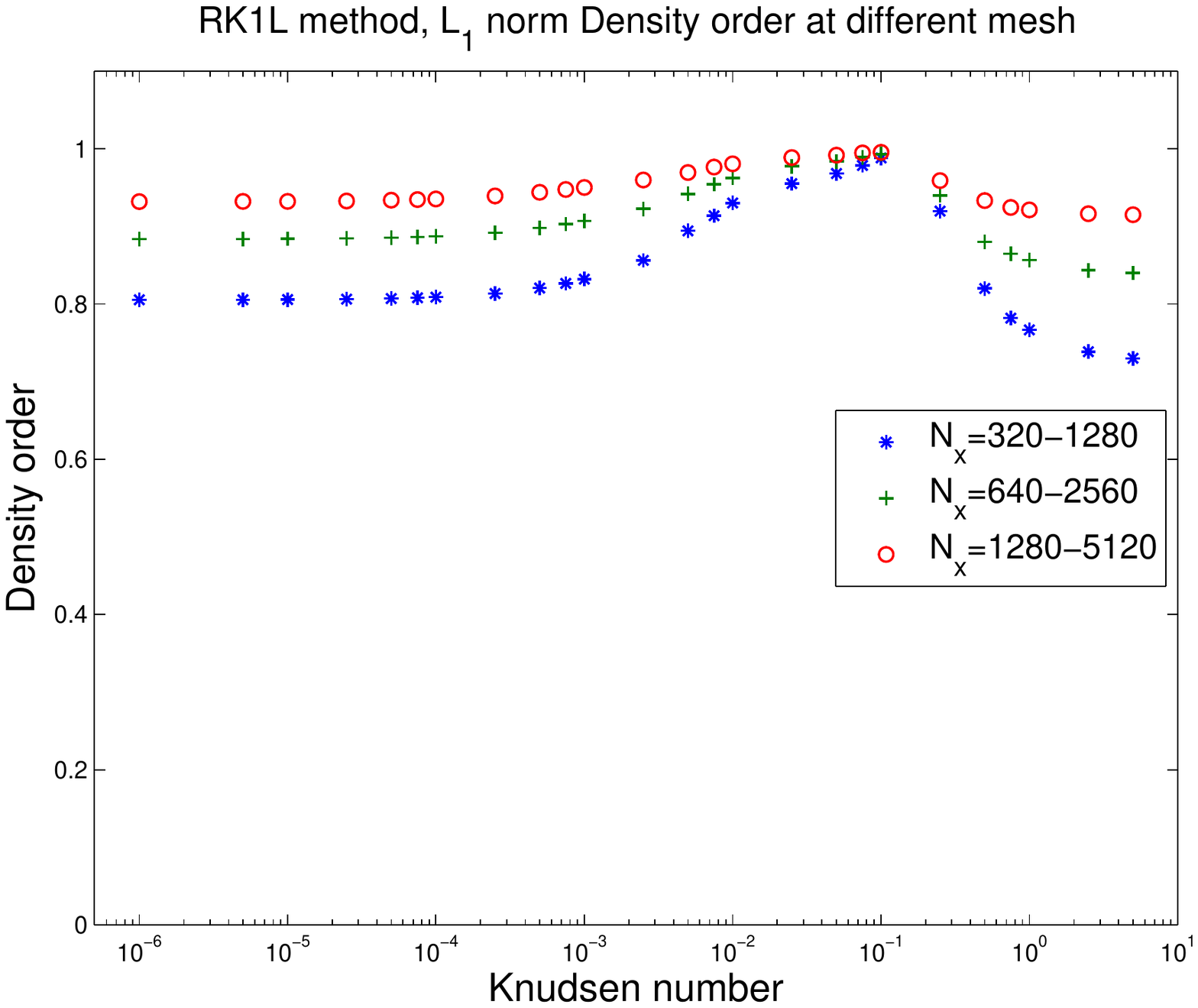}}
\caption{\footnotesize{$L_{1}$ error and accuracy order of implicit Euler methods coupled with linear interpolation, varying $\varepsilon$, using periodic boundary condition.}}\label{fig_RK1_order_periodic}
\end{figure}
The initial condition for the distribution function is the Maxwellian, computed by given macroscopic fields. To checked the accuracy order the solution must be smooth. Using periodic or reflective boundary conditions, we observe that some shocks appear in the solution around the time $t=0.35,$ so to test the accuracy order we use as final time $0.32,$ that is large enough to reach thermodynamic equilibrium.
\begin{figure}[h]
\subfigure
{\includegraphics[trim=3cm 7cm 2cm 8.5cm,scale=0.4]{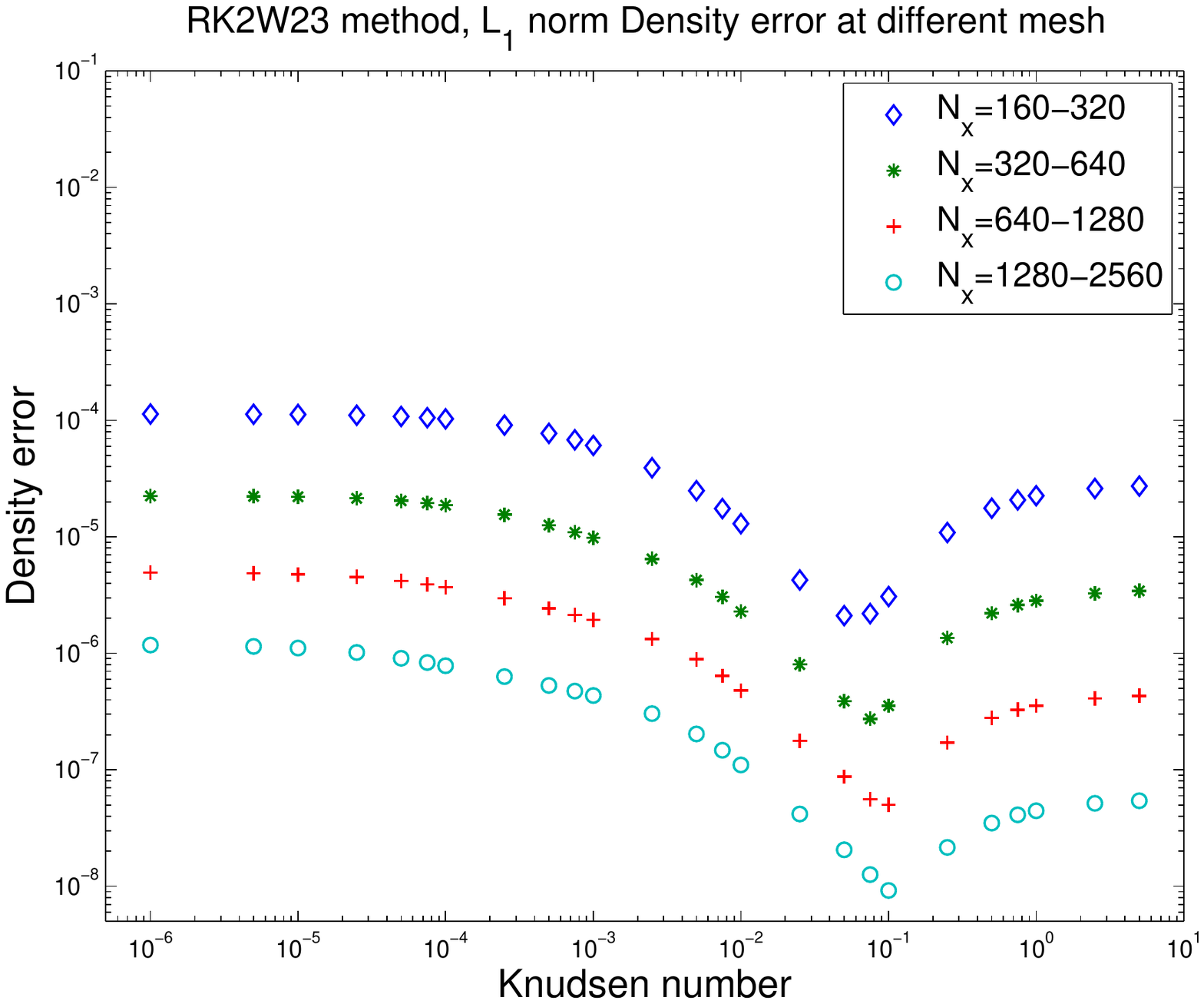}}
\subfigure
{\includegraphics[trim=1cm 7cm 2cm 8.5cm,scale=0.4]{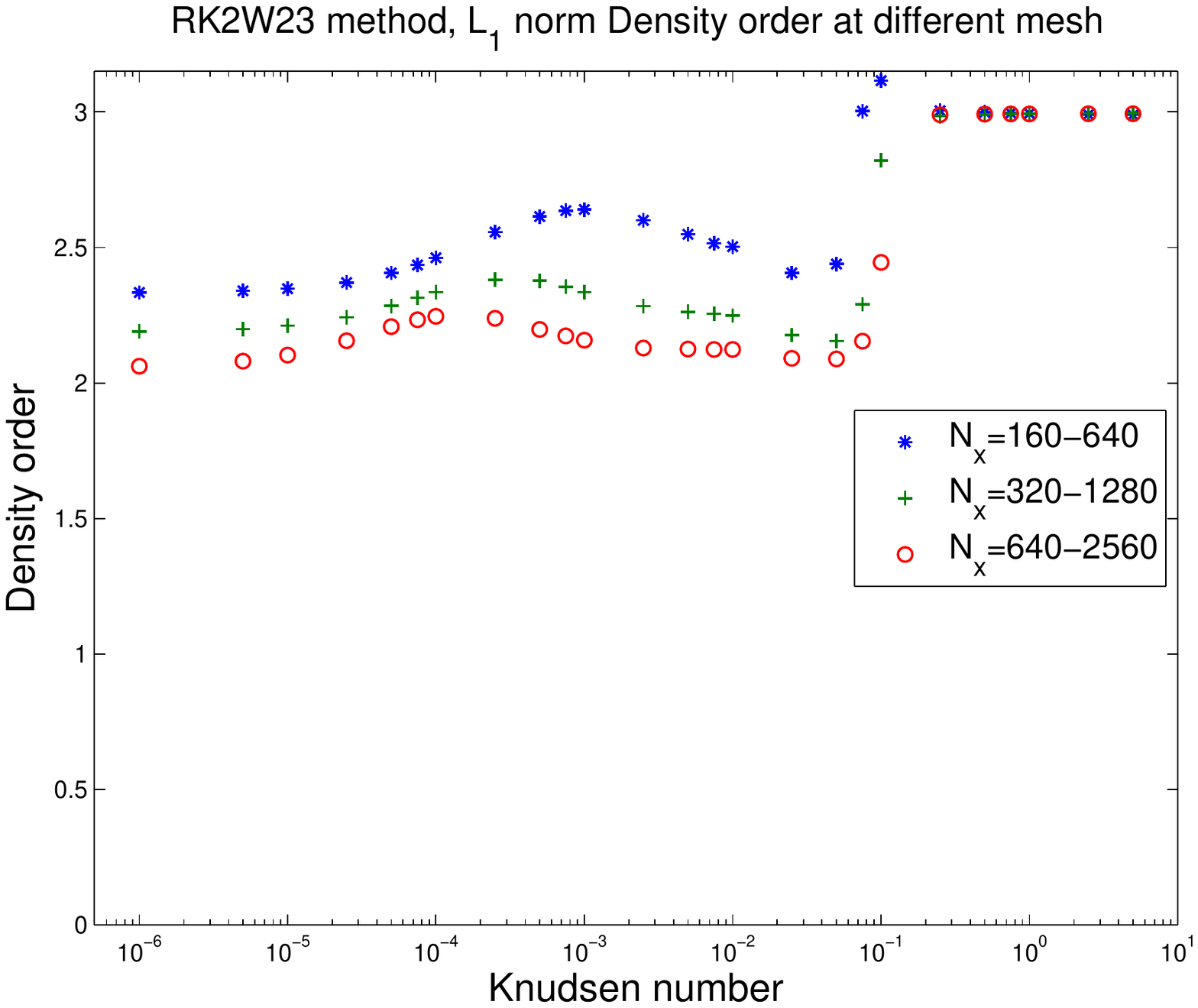}}
\subfigure
{\includegraphics[trim=3cm 7cm 2cm 7cm,scale=0.4]{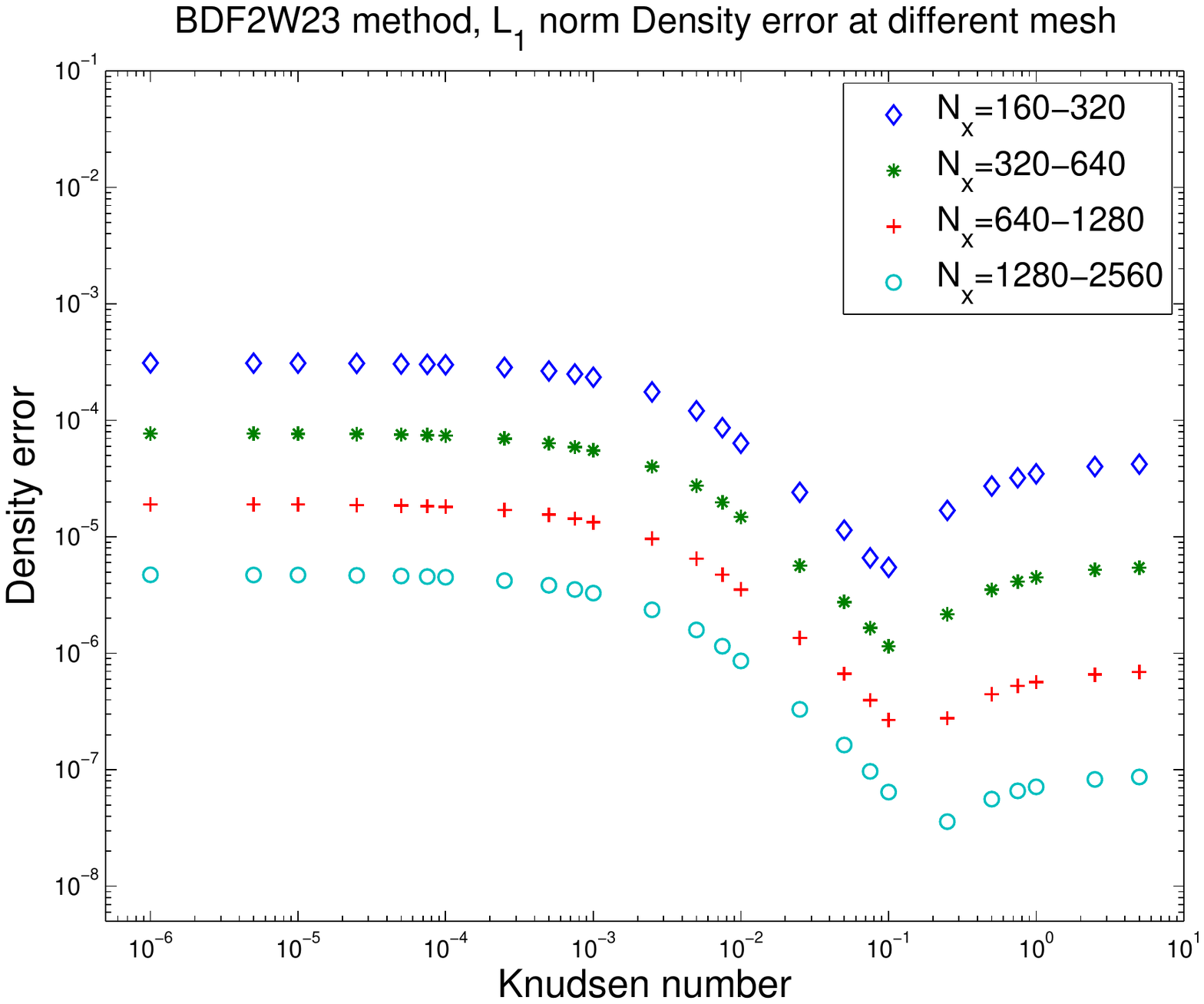}}
\subfigure
{\includegraphics[trim=1cm 7cm 2cm 7cm,scale=0.4]{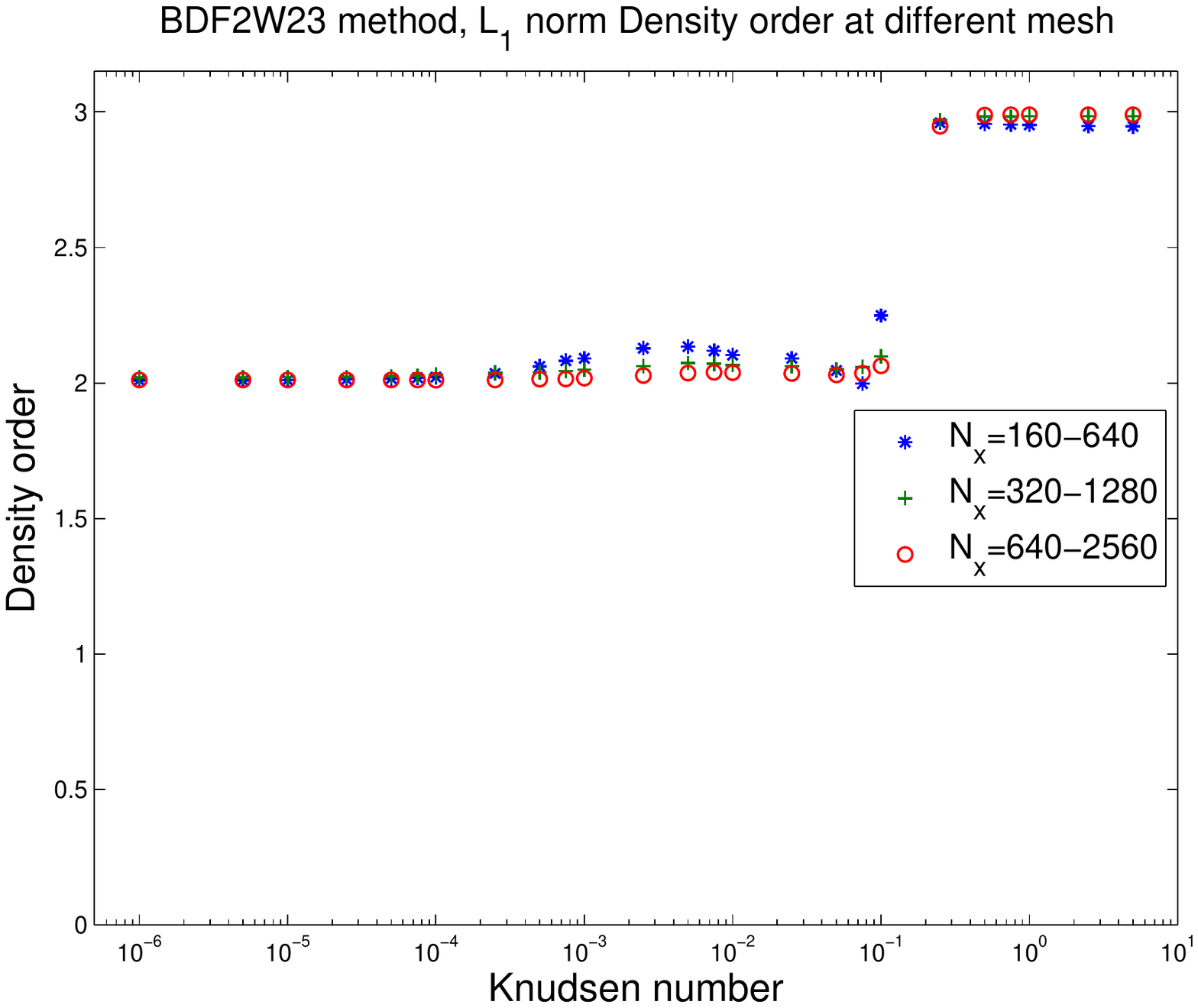}}
\caption{\footnotesize{$L_{1}$ error and accuracy order of RK2W23 and BDF2W23, varying $\varepsilon$, using periodic boundary condition.}}\label{fig_RK2W23_order_periodic}
\end{figure}
\noindent
In all tests  $N_{v}=20$ velocity points have been used, uniformly spaced in $[-10,10].$ For the time step, we set $\Delta t=$CFL$\,\Delta x/v_{max}$ and we have used CFL$=4.$ The spatial domain is $[-1,1].$\\
We have compared the following method:
\begin{itemize}
\item[-] First order implicit Euler coupled with linear interpolation;
\item[-] RK2W23 and BDF2W23 as second order methods;
\item[-] RK3W23, RK3W35, BDF3W23 and BDF3W35 as third order methods;
\end{itemize}
The Figures \ref{fig_RK1_order_periodic}-\ref{fig_BDF3W35_order_periodic} show the $L_{1}$ error and the rate of convergence related to macroscopic density of the scheme just mentioned using periodic boundary conditions. The same behaviour is observed when monitoring the error in mean velocity and in energy.
\subsubsection*{Remarks}
\begin{itemize}
\item[-] In most regimes the order of accuracy is the theoretical one. More precisely all schemes maintain the theoretical order of accuracy in the limit of small Knudsen number, except RK3-based scheme, whose order of accuracy degrades to 2, with both WENO23 and WENO35 interpolation. Some schemes (RK2W23, RK3W35, BDF2W23, BDF3W35) present a spuriously high order of accuracy for large Knudsen number. This is due to the fact that for such large Knudsen number and small final time most error is due to space discretization, which in such schemes is of order higher than time discretization.\\ The most uniform accuracy is obtained by the BDF3W23 scheme.
\item[-] Most tests have been conducted with periodic boundary condition. Similar results are obtained using reflecting boundary conditions (see Fig. \ref{fig_third_order_3D}).
\end{itemize}

\subsection{Optimal CFL}

The semilagrangian nature of the scheme allows us to avoid the classical CFL stability restriction. In this way, one can use large CFL numbers in order to obtain larger time step thus lowering the computational cost. How much can we increase the CFL number without degrading the accuracy?\\ Consistency analysis of semilagrangian schemes \cite{FalconeFerretti} shows that the error is composed by two part: one depending by the time integration and one depending on the interpolation. Therefore, if we use a small CFL number, the time step will be small and the error will be mainly due to the interpolation.

\noindent
On the other hand, if we use a big CFL number, the error will be mainly due to the time integration. This argument leads us to think that there is an optimal value of the CFL number, that allows us to minimize the error. The following Figures \ref{fig_second_order_optimal}-\ref{fig_third_order_optimal} show this behavior. Each picture shows the $L_{2}$ error of the macroscopic density of the previous smooth initial data, varying the CFL number from 0.05 to 20. The grid of the CFL values is not uniform because we want to work with constant time step until the final time, that for this test is 0.3. The error is computed using two numerical solution, obtained with $N_{x}=160$ and $N_{x}=320.$ $\varepsilon$ is fixed to the value $10^{-4}.$\\
When using accuracy in space which is not much larger than in time, as in the case of RK2W23, RK3W23, BDF3W23, an evident optimal CFL number appears, when interpolation error and time discretization error balance.\\ If space discretization is much more accurate than time discretization, the optimal CFL number decreases. Note that with the same formal order of accuracy, the optimal CFL number is larger for schemes based on RK than for schemes based on BDF, because RK have a smaller error constant.

\begin{figure}[htbp]
\subfigure
{\includegraphics[trim=3cm 7cm 2cm 8.5cm,scale=0.4]{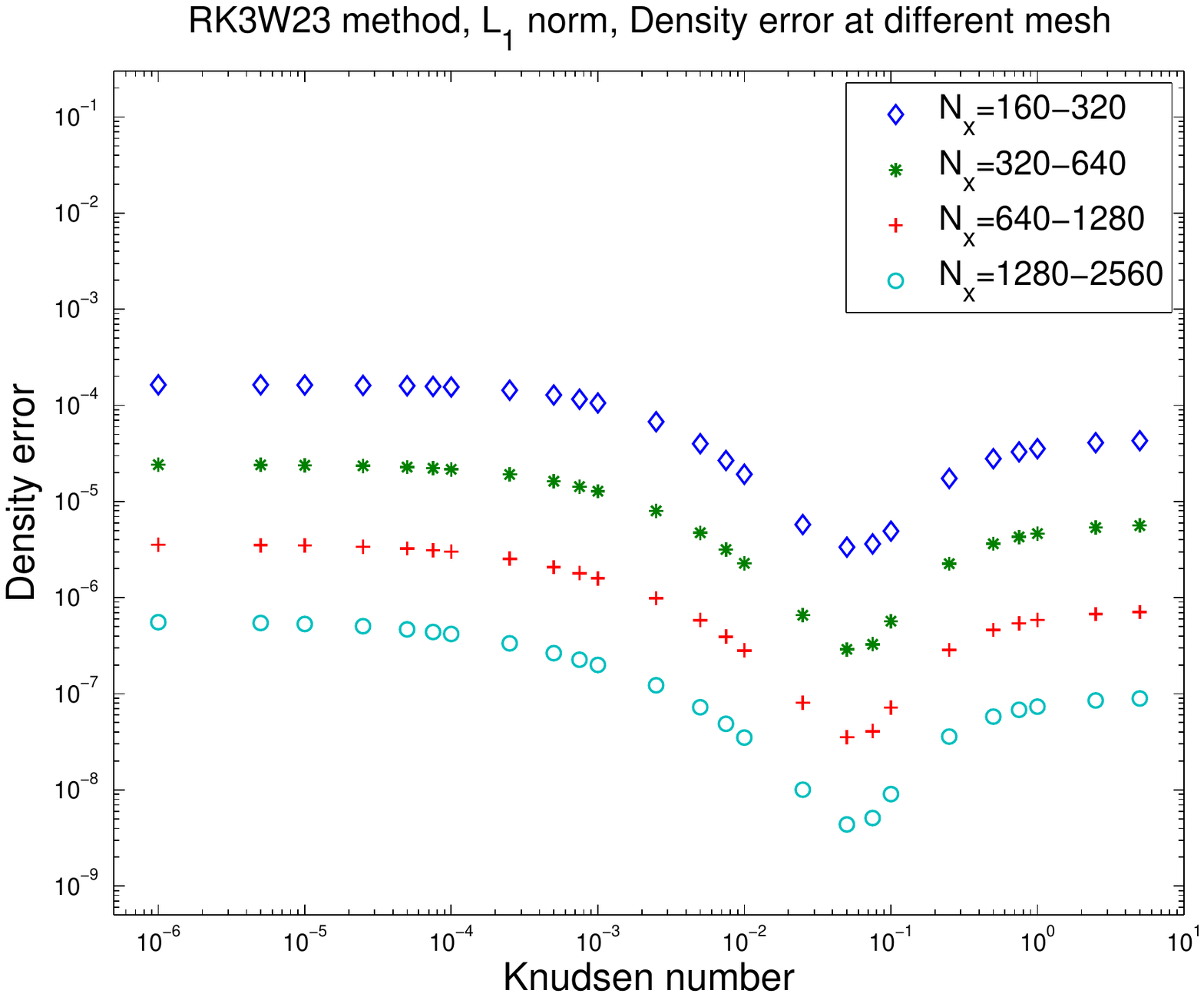}}
\subfigure
{\includegraphics[trim=1cm 7cm 2cm 8.5cm,scale=0.4]{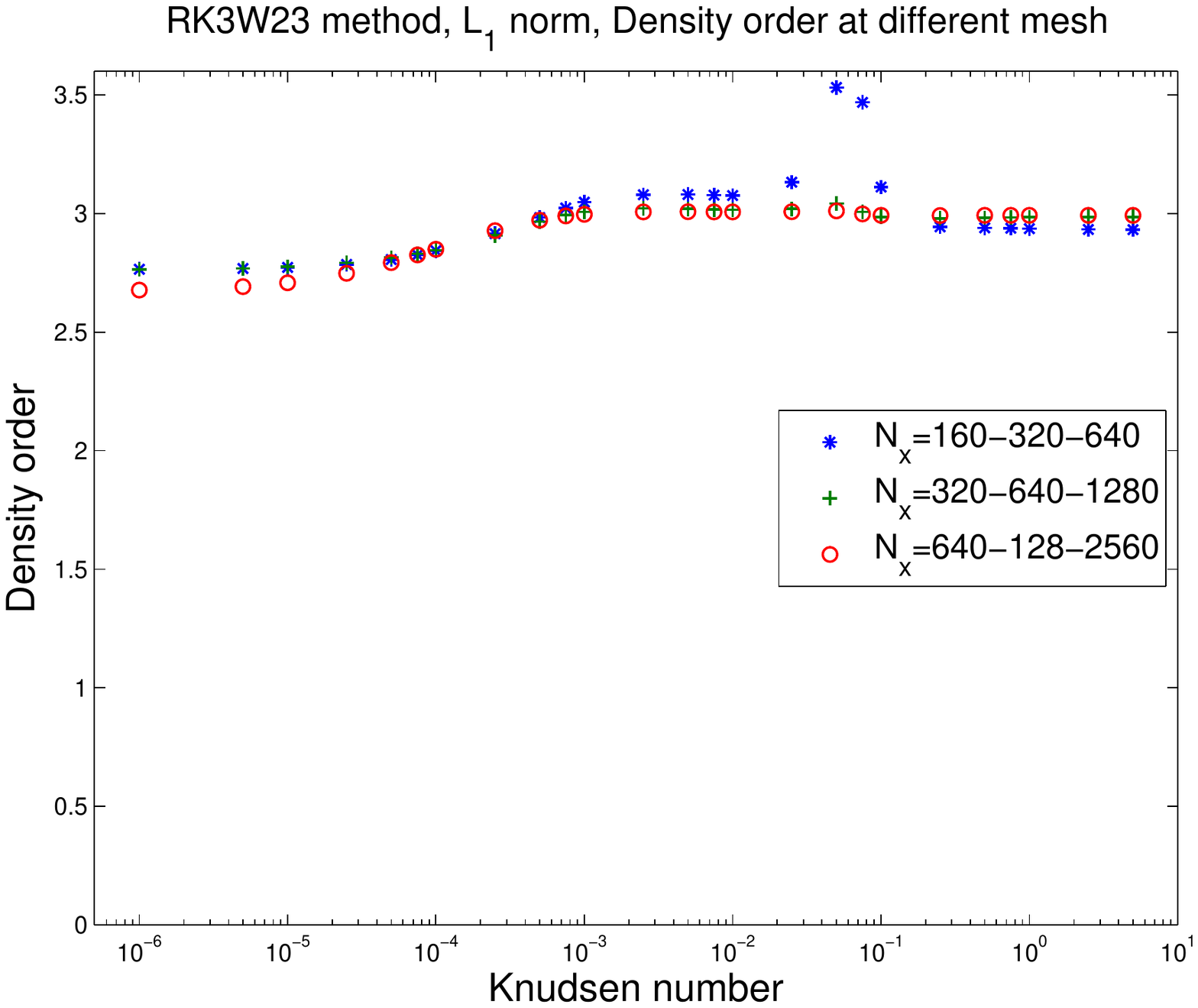}}
\subfigure
{\includegraphics[trim=3cm 7cm 2cm 7cm,scale=0.4]{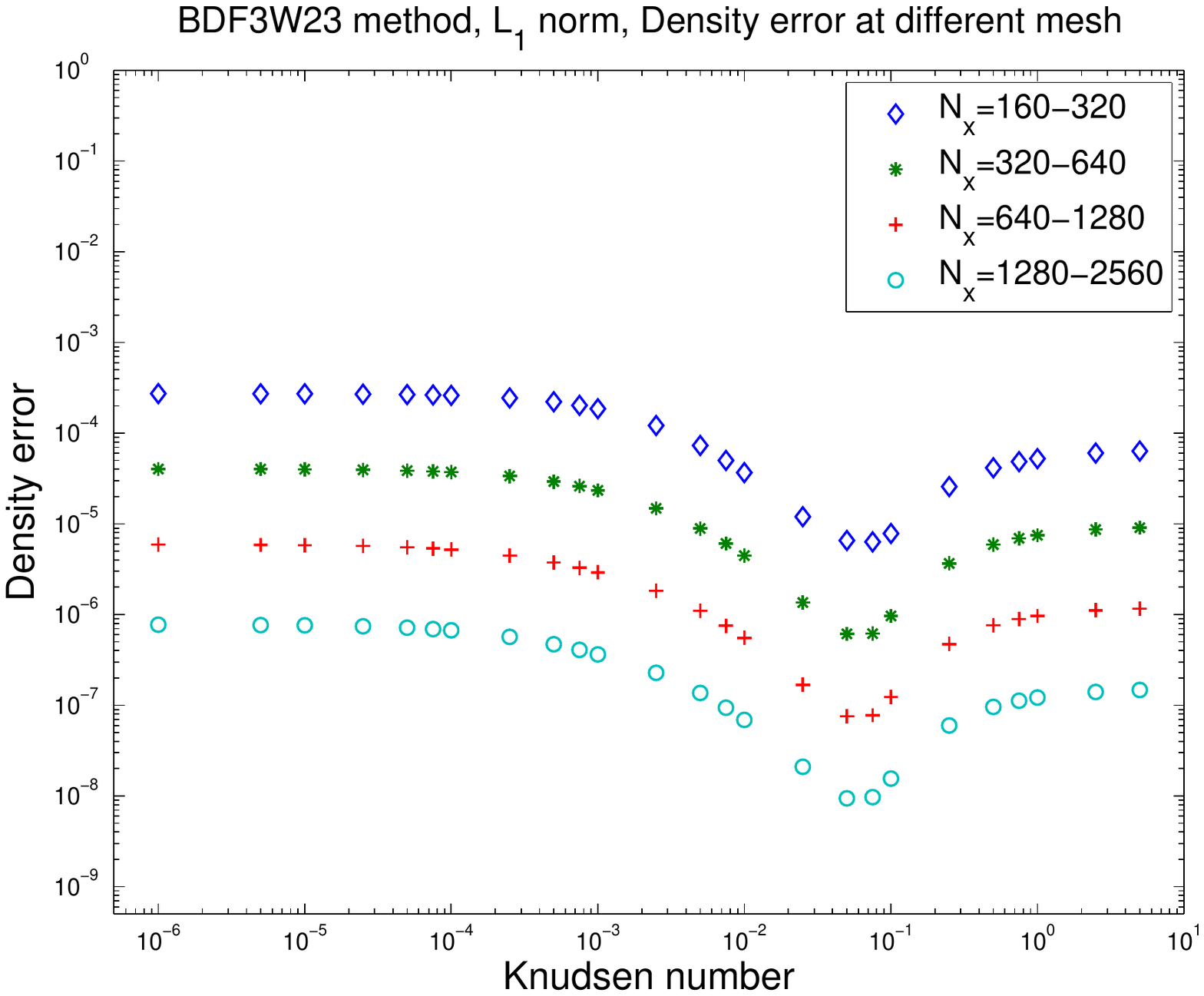}}
\subfigure
{\includegraphics[trim=1cm 7cm 2cm 7cm,scale=0.4]{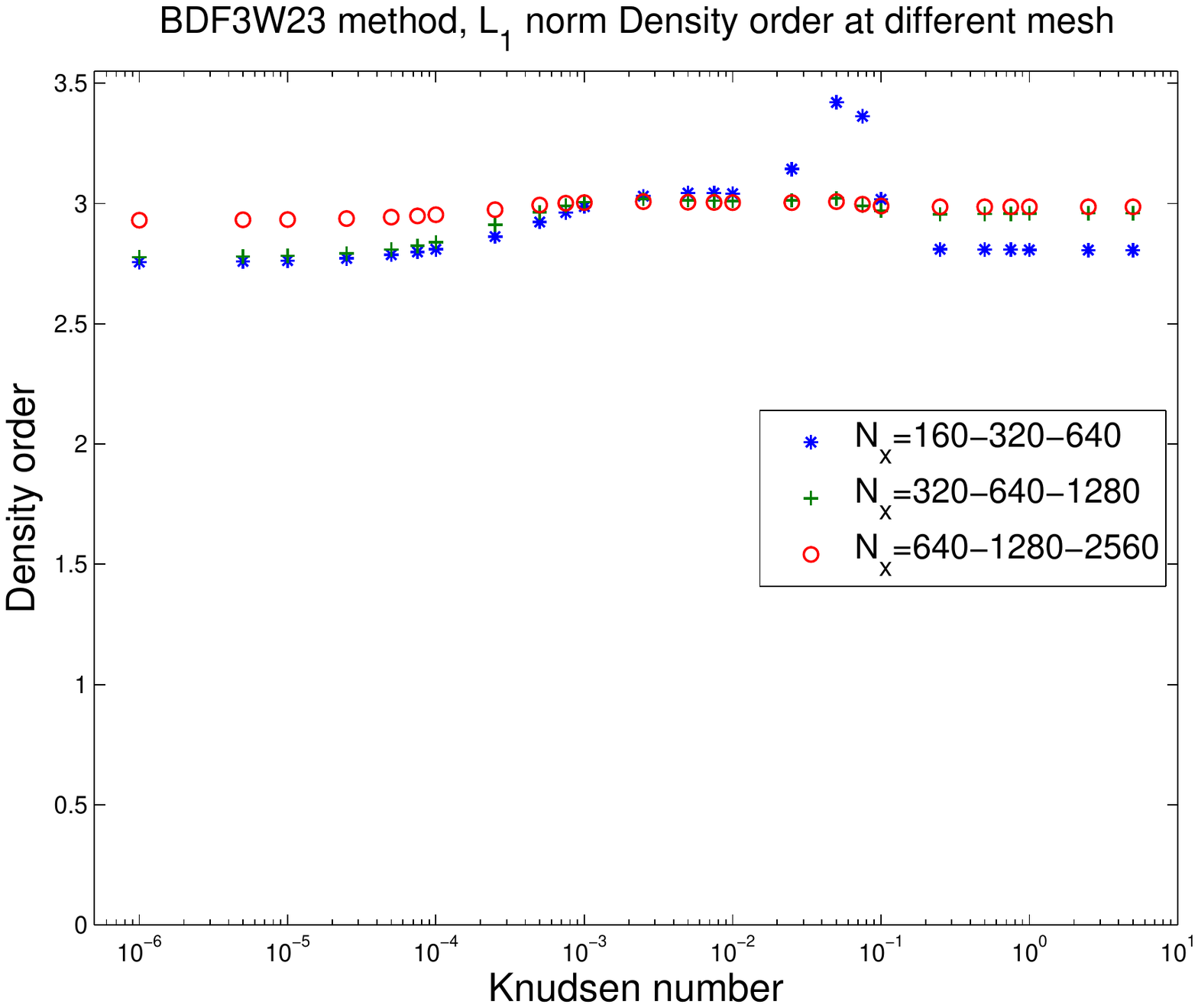}}
\caption{\footnotesize{$L_{1}$ error and accuracy order of RK3 and BDF3 methods coupled with WENO23, varying $\varepsilon$, using reflective boundary condition.}}\label{fig_third_order_reflective}
\end{figure}

\begin{figure}[h]
\subfigure
{\includegraphics[trim=3cm 7cm 2cm 7cm,scale=0.4]{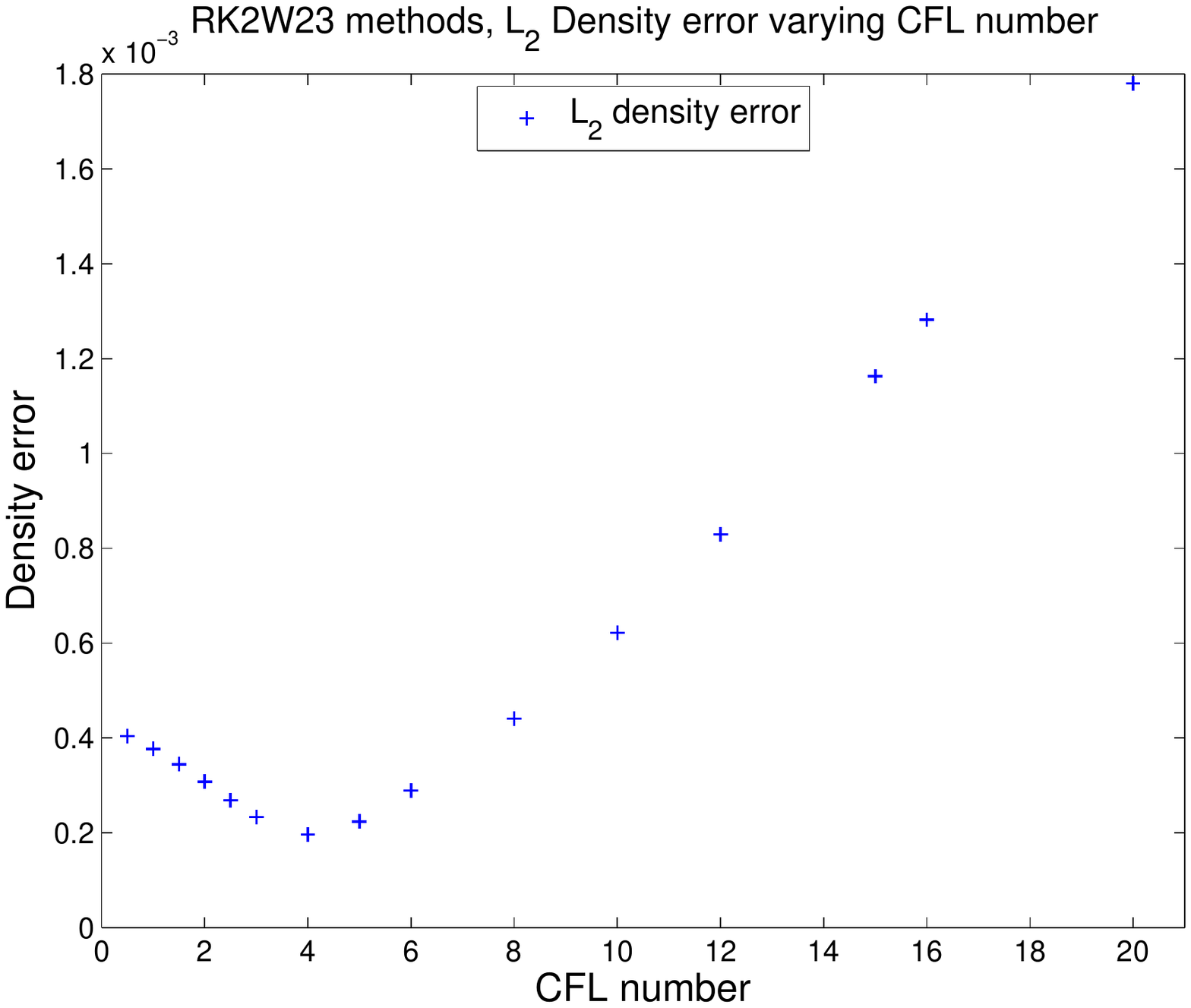}}
\subfigure
{\includegraphics[trim=1cm 7cm 2cm 7cm,scale=0.4]{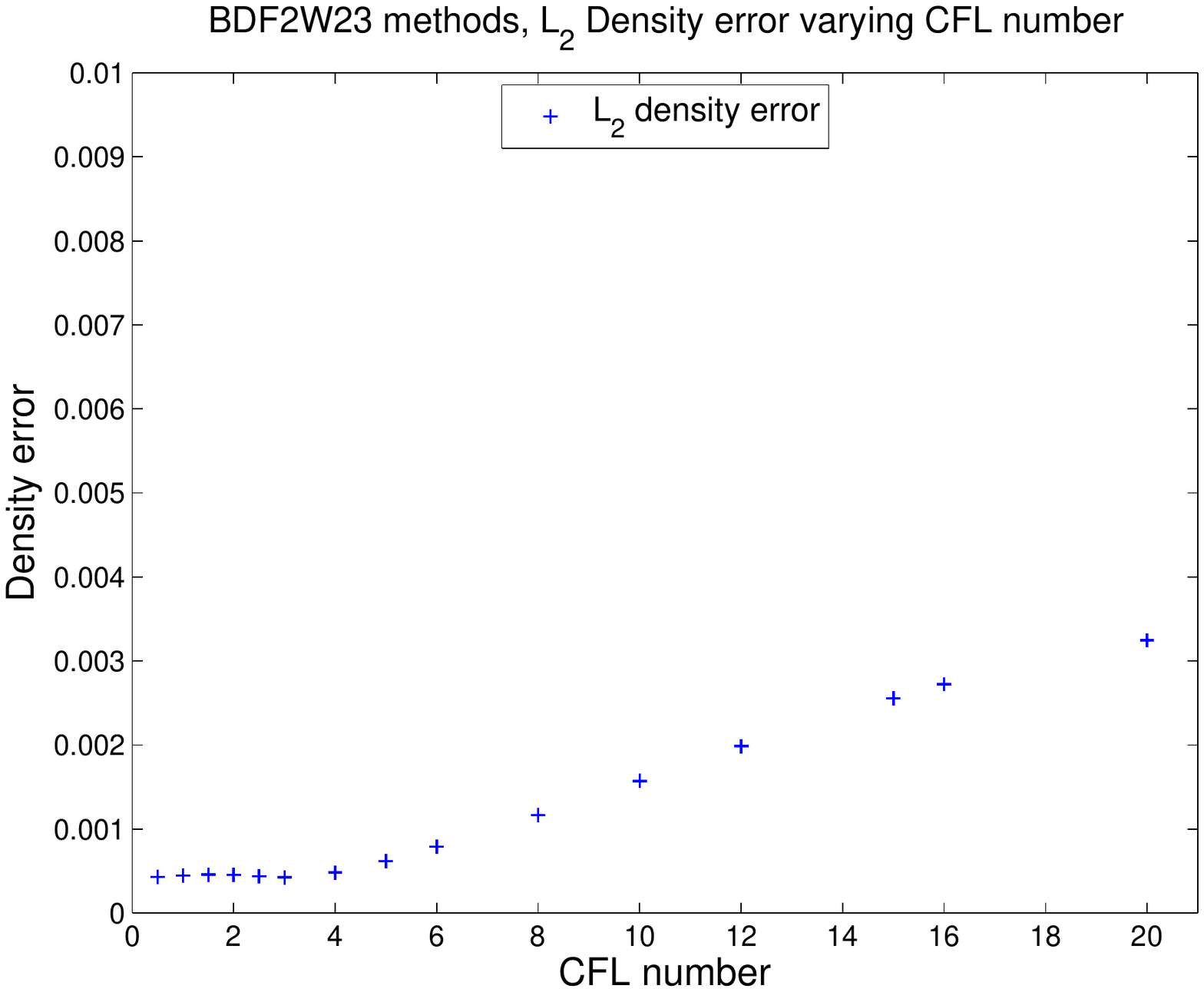}}
\subfigure
{\includegraphics[trim=3cm 7cm 2cm 7cm,scale=0.4]{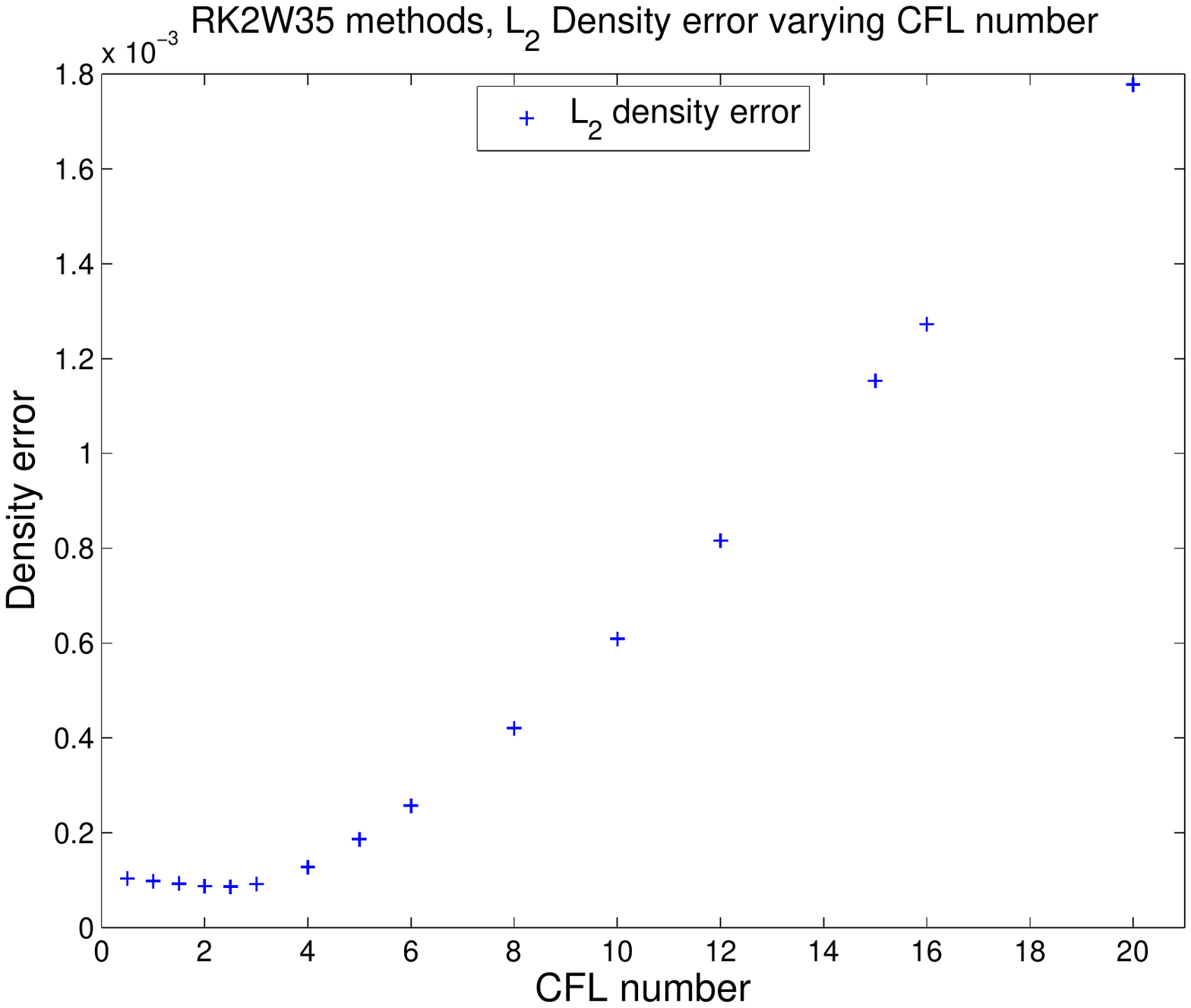}}
\subfigure
{\includegraphics[trim=1cm 7cm 2cm 7cm,scale=0.4]{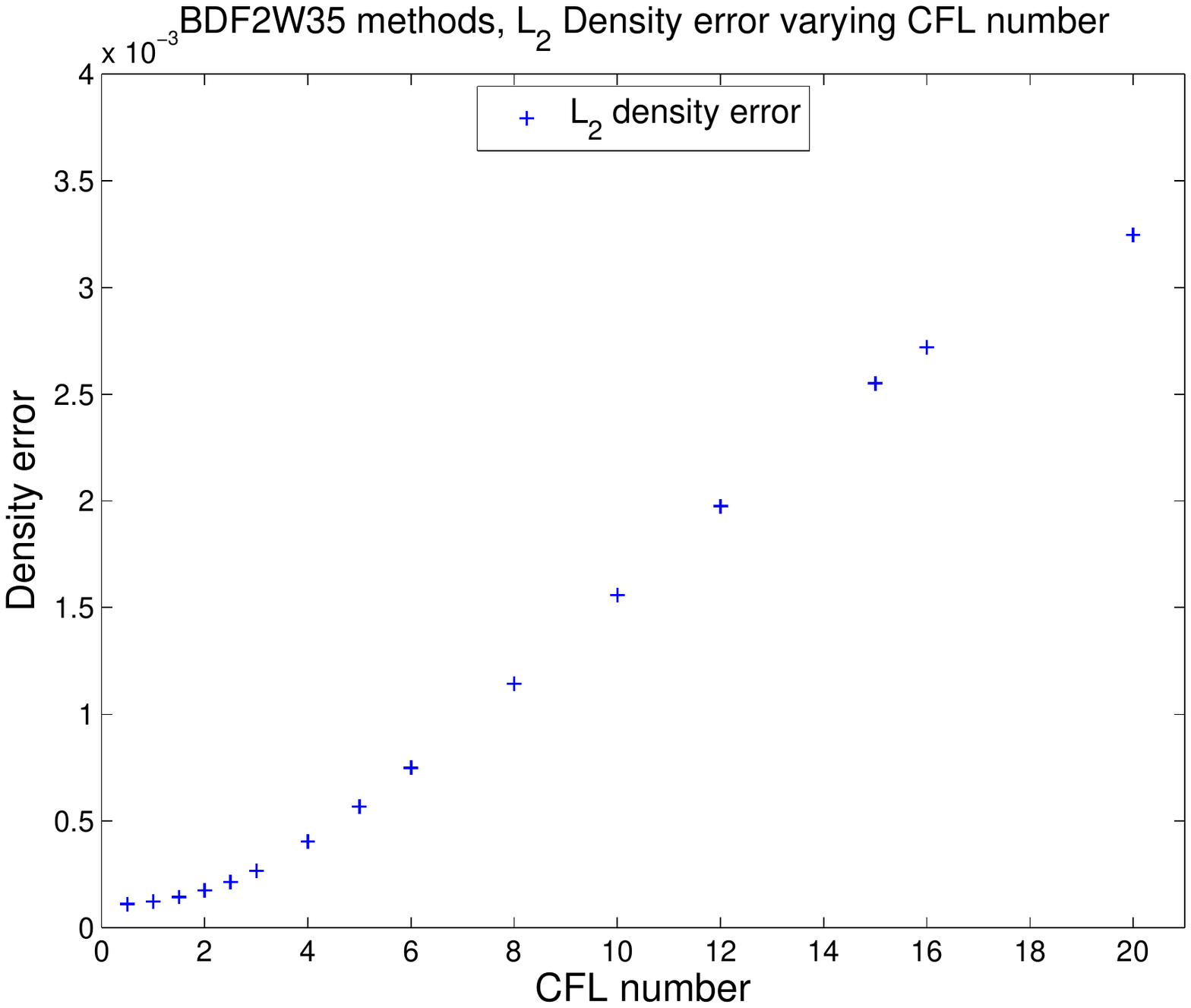}}
\caption{\footnotesize{Optimal CFL number. Left RK2, right BDF2. From top to bottom: WENO23, WENO35.}}\label{fig_second_order_optimal}
\end{figure}

\begin{figure}[h]
\subfigure
{\includegraphics[trim=3cm 7cm 2cm 7cm,scale=0.4]{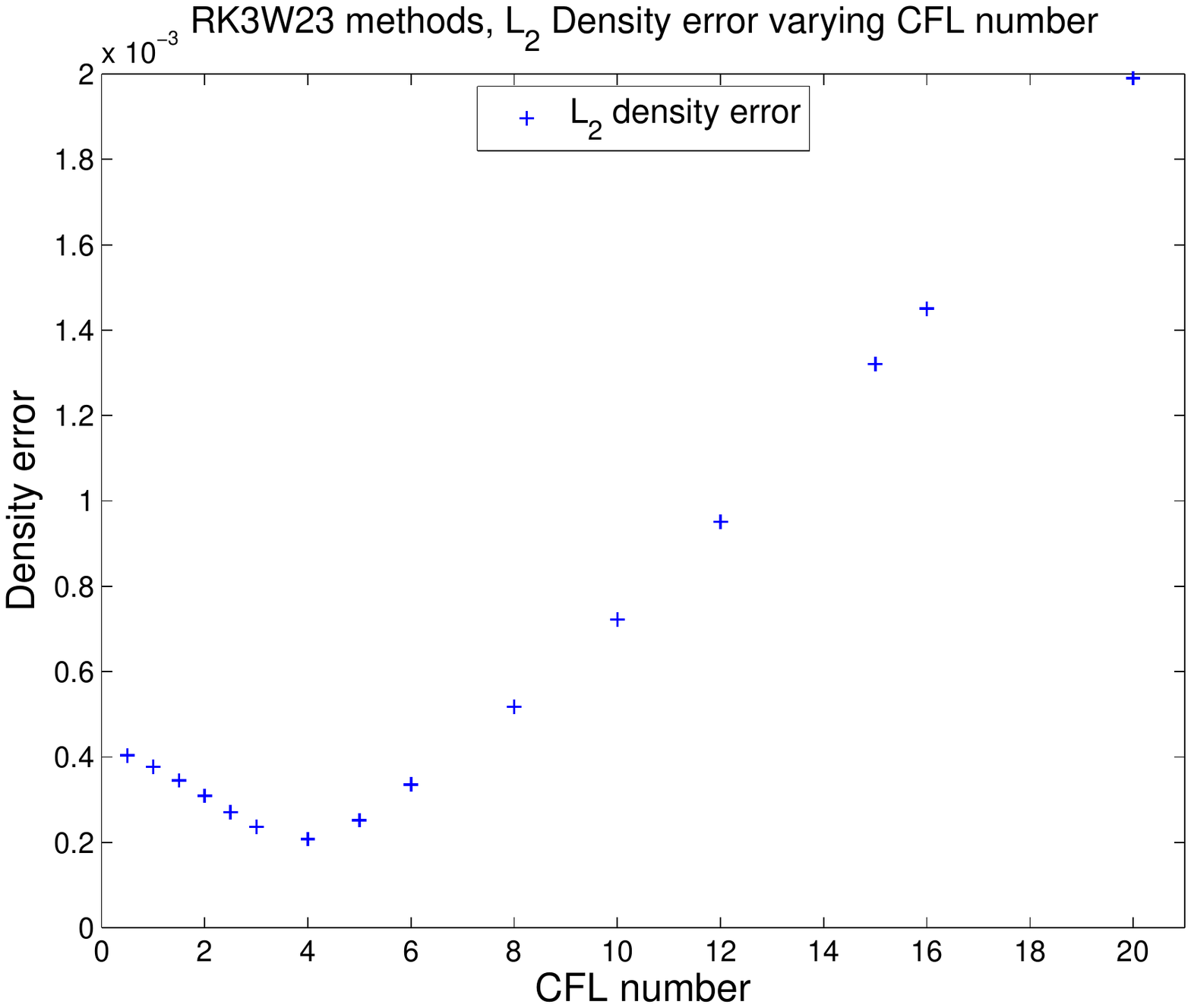}}
\subfigure
{\includegraphics[trim=1cm 7cm 2cm 7cm,scale=0.4]{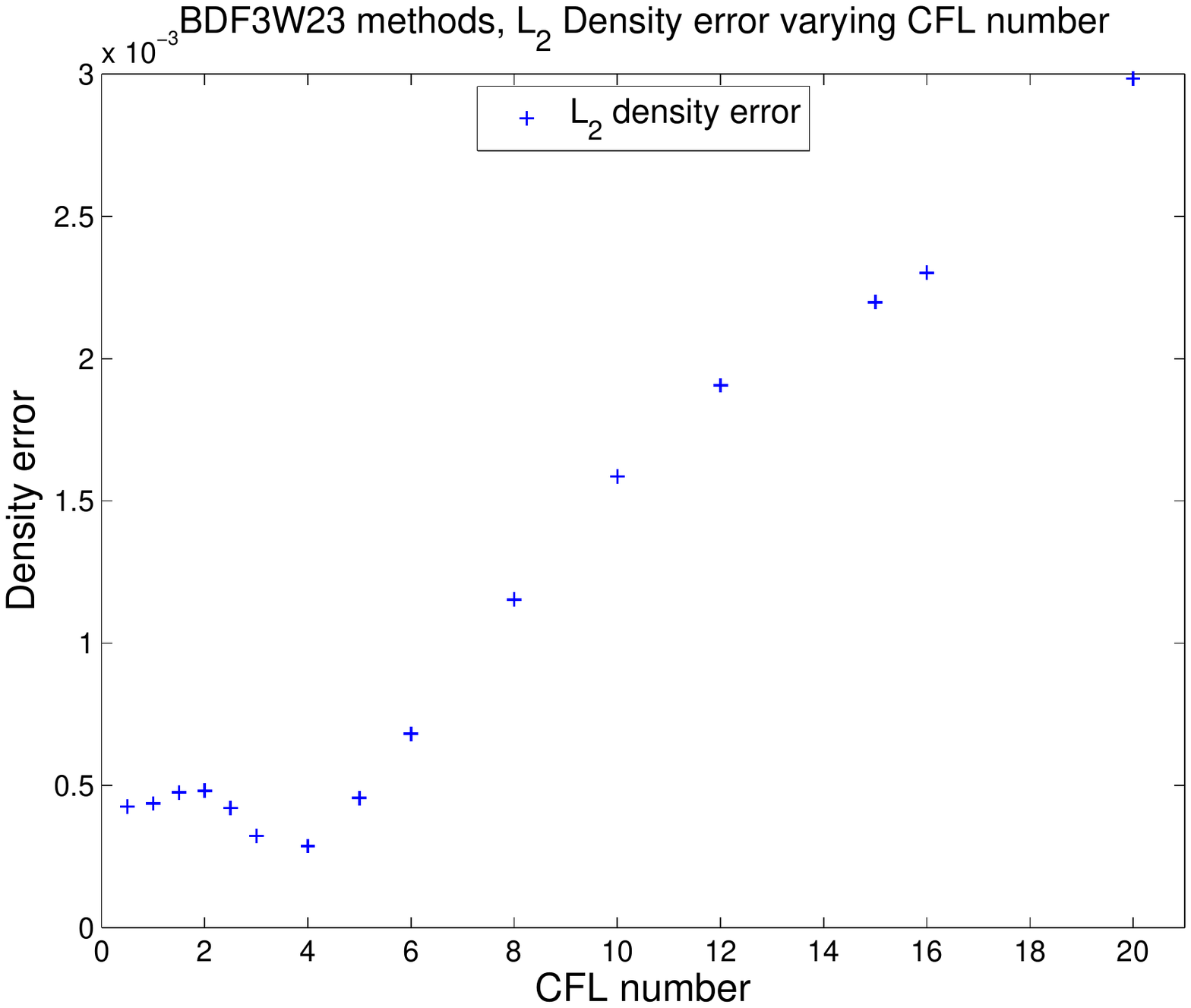}}
\subfigure
{\includegraphics[trim=3cm 7cm 2cm 7cm,scale=0.4]{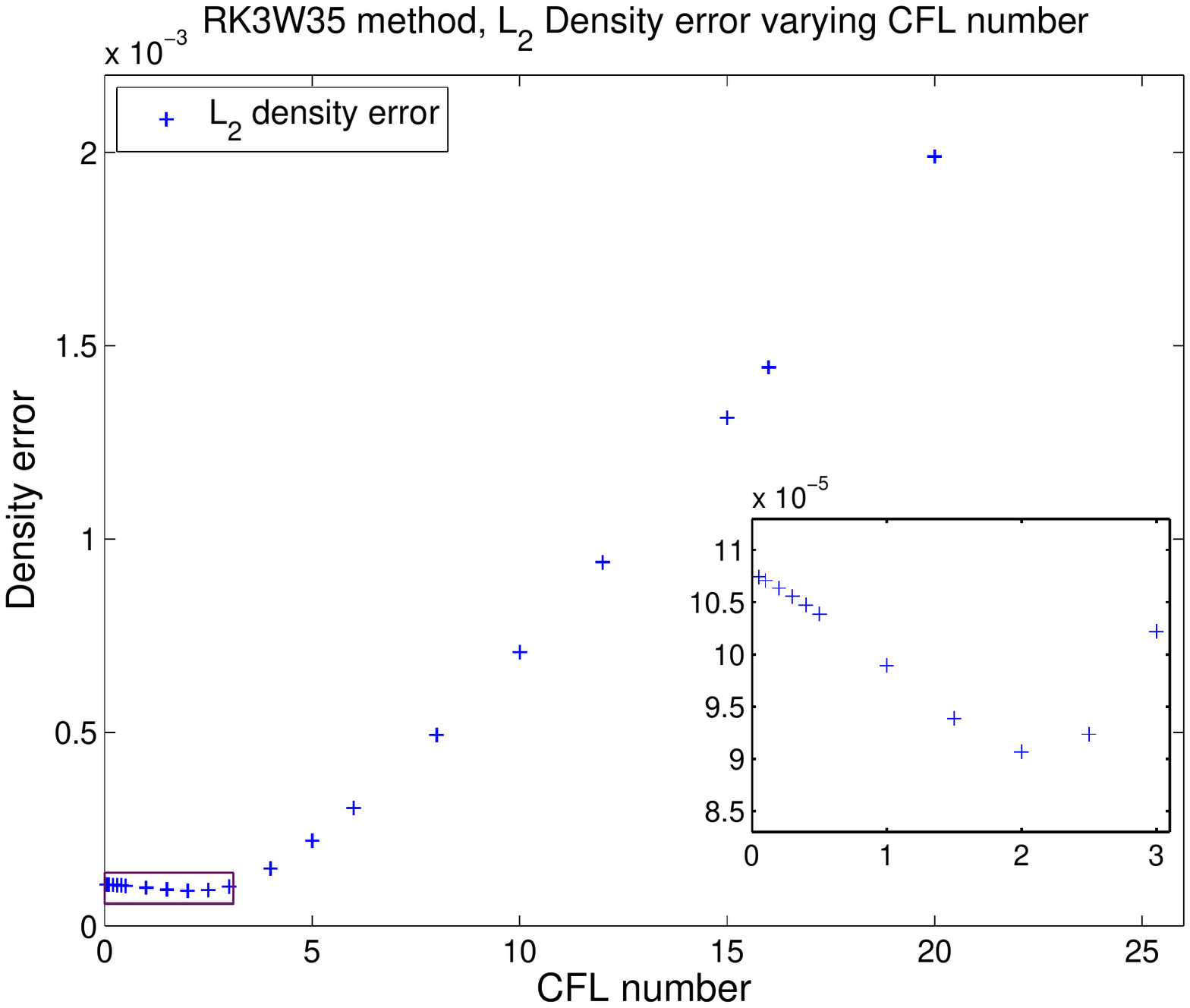}}
\subfigure
{\includegraphics[trim=1cm 7cm 2cm 7cm,scale=0.4]{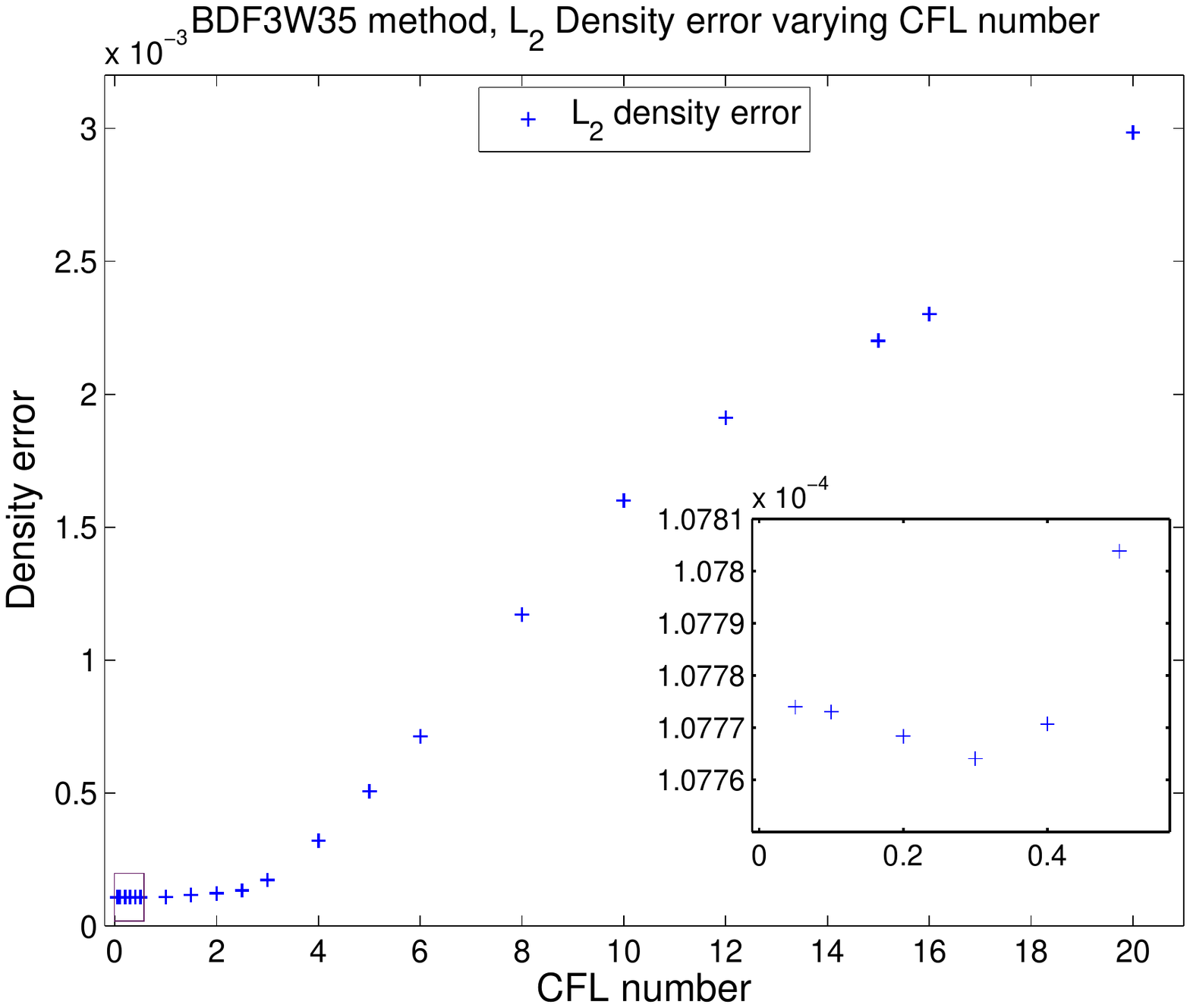}}
\caption{\footnotesize{Optimal CFL number. Left RK3, right BDF3. From top to bottom: WENO23, WENO35.}}\label{fig_third_order_optimal}
\end{figure}

\subsection{Riemann problem}
\begin{figure}[htbp]
\subfigure
{\includegraphics[trim=5cm 7cm 2cm 12cm,scale=0.45]{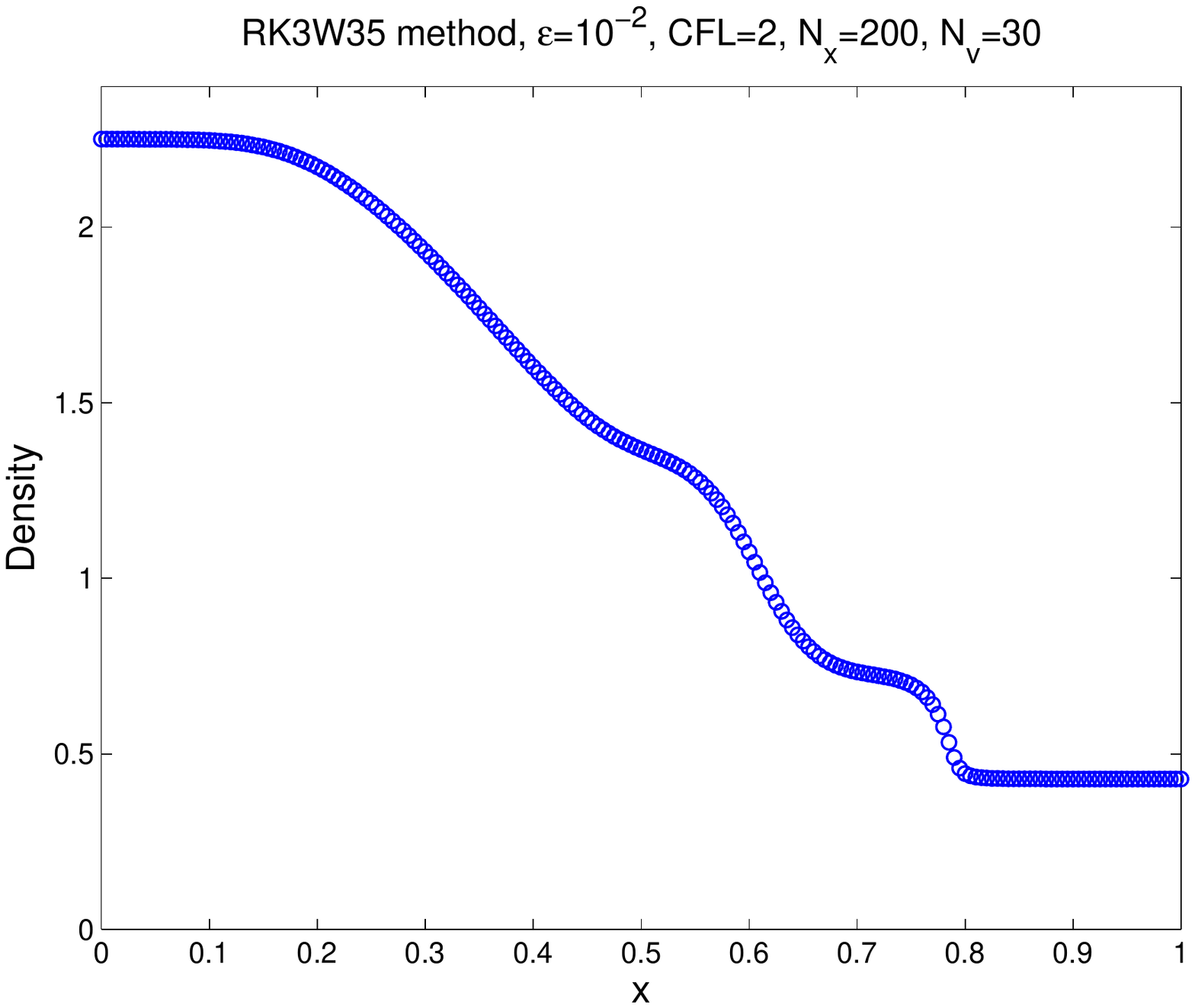}}
\subfigure
{\includegraphics[trim=1cm 7cm 2cm 12cm,scale=0.45]{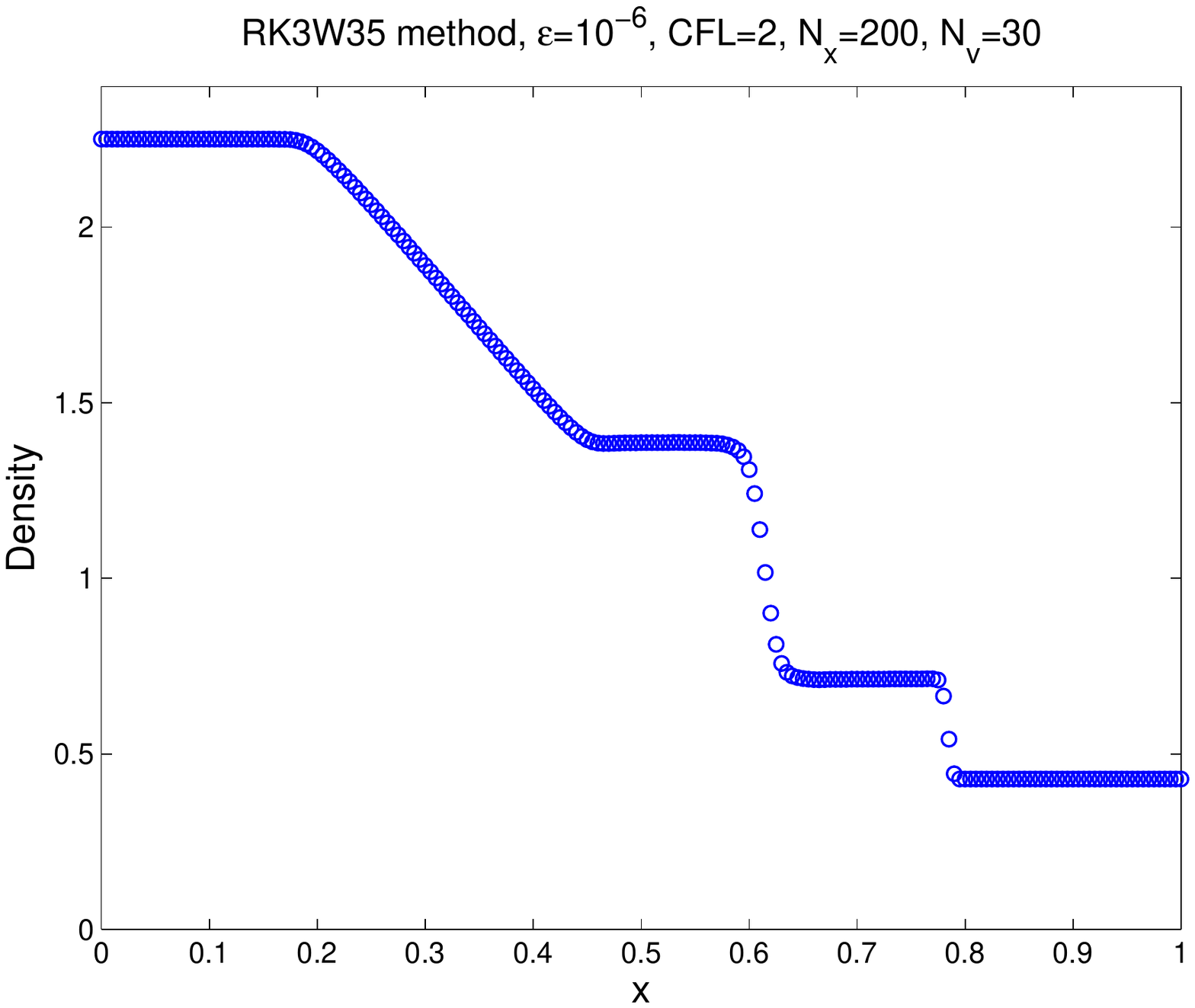}}
\subfigure
{\includegraphics[trim=5cm 7cm 2cm 7cm,scale=0.45]{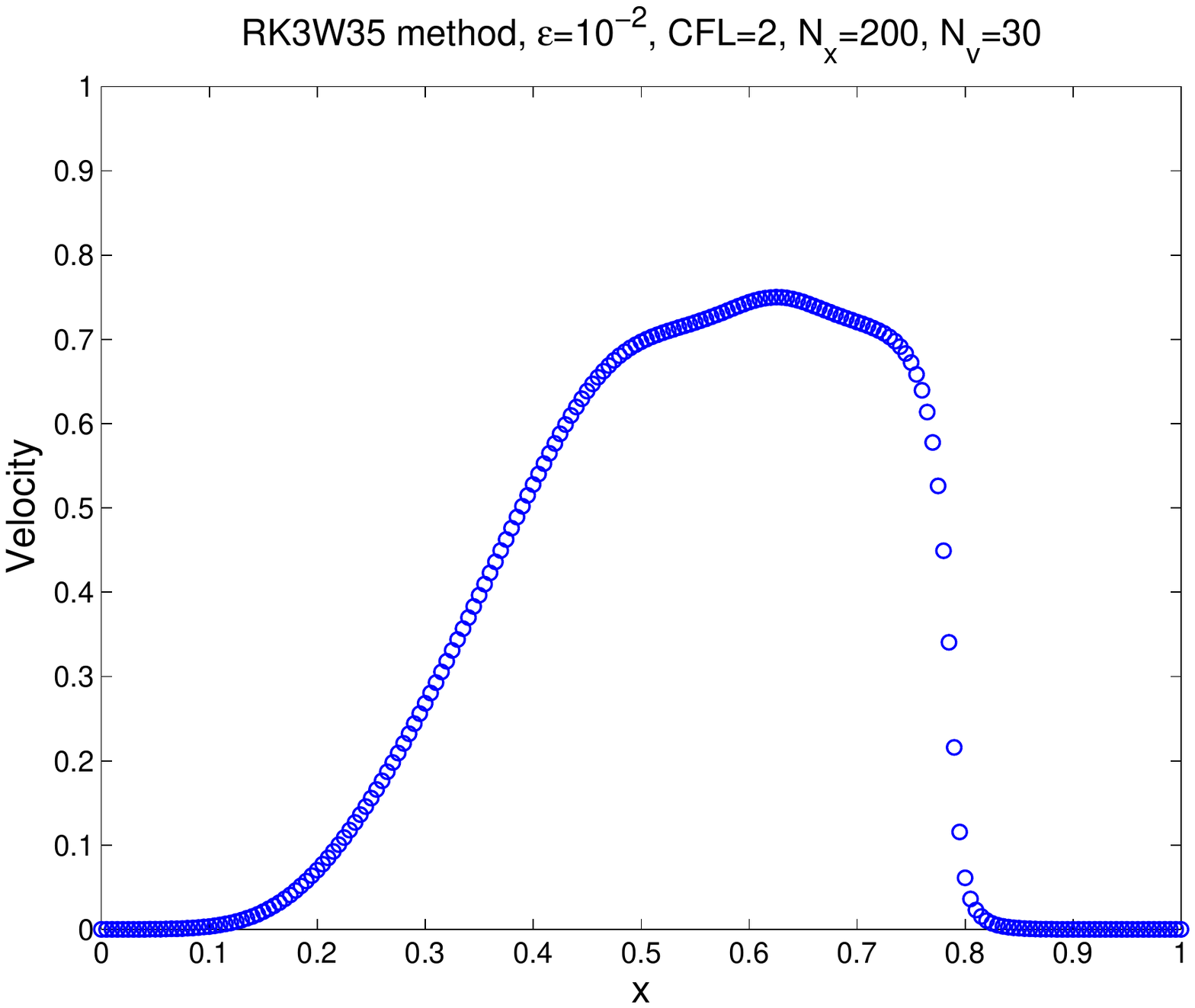}}
\subfigure
{\includegraphics[trim=1cm 7cm 2cm 7cm,scale=0.45]{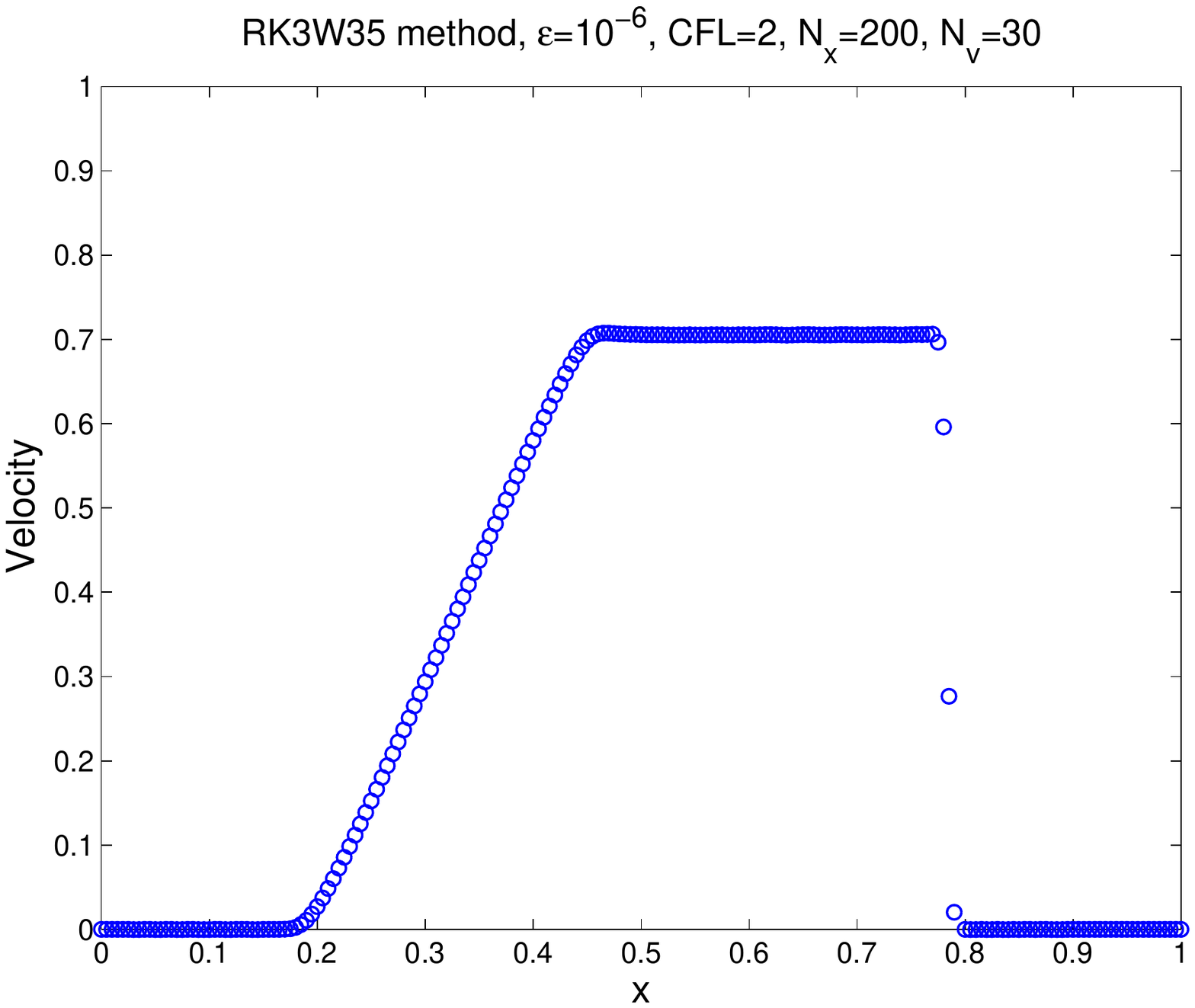}}
\subfigure
{\includegraphics[trim=5cm 7cm 2cm 7cm,scale=0.45]{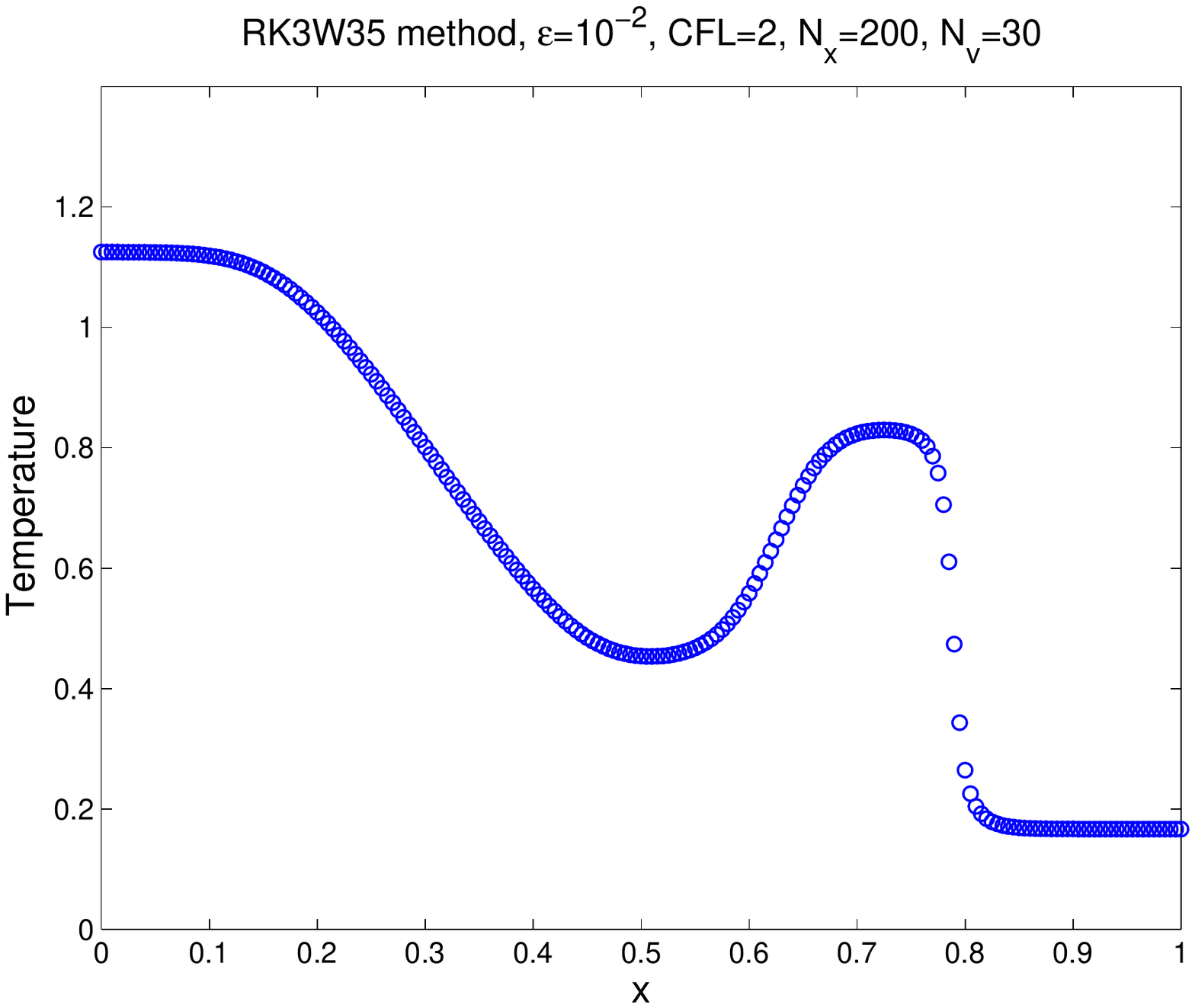}}
\subfigure
{\includegraphics[trim=1cm 7cm 2cm 7cm,scale=0.45]{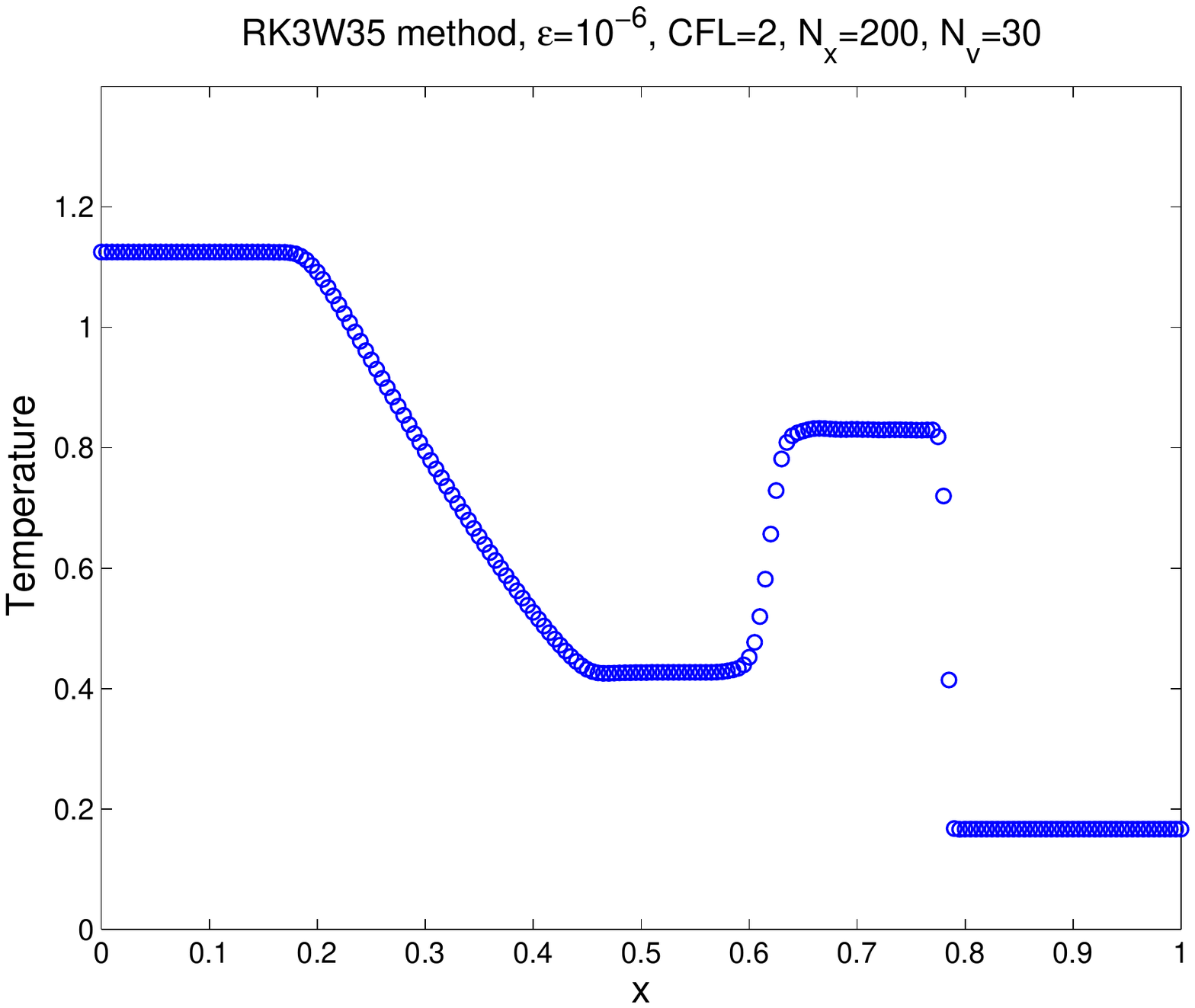}}
\caption{\footnotesize{RK3W35 scheme. Riemann problem in 1D space and velocity case. Left $\varepsilon=10^{-2}$; Right $\varepsilon=10^{-6}$. From top to bottom: Density, Velocity and Temperature.}}\label{fig_RK3_riemann}
\end{figure}

\begin{figure}[htbp]
\subfigure
{\includegraphics[trim=5cm 7cm 2cm 12cm,scale=0.45]{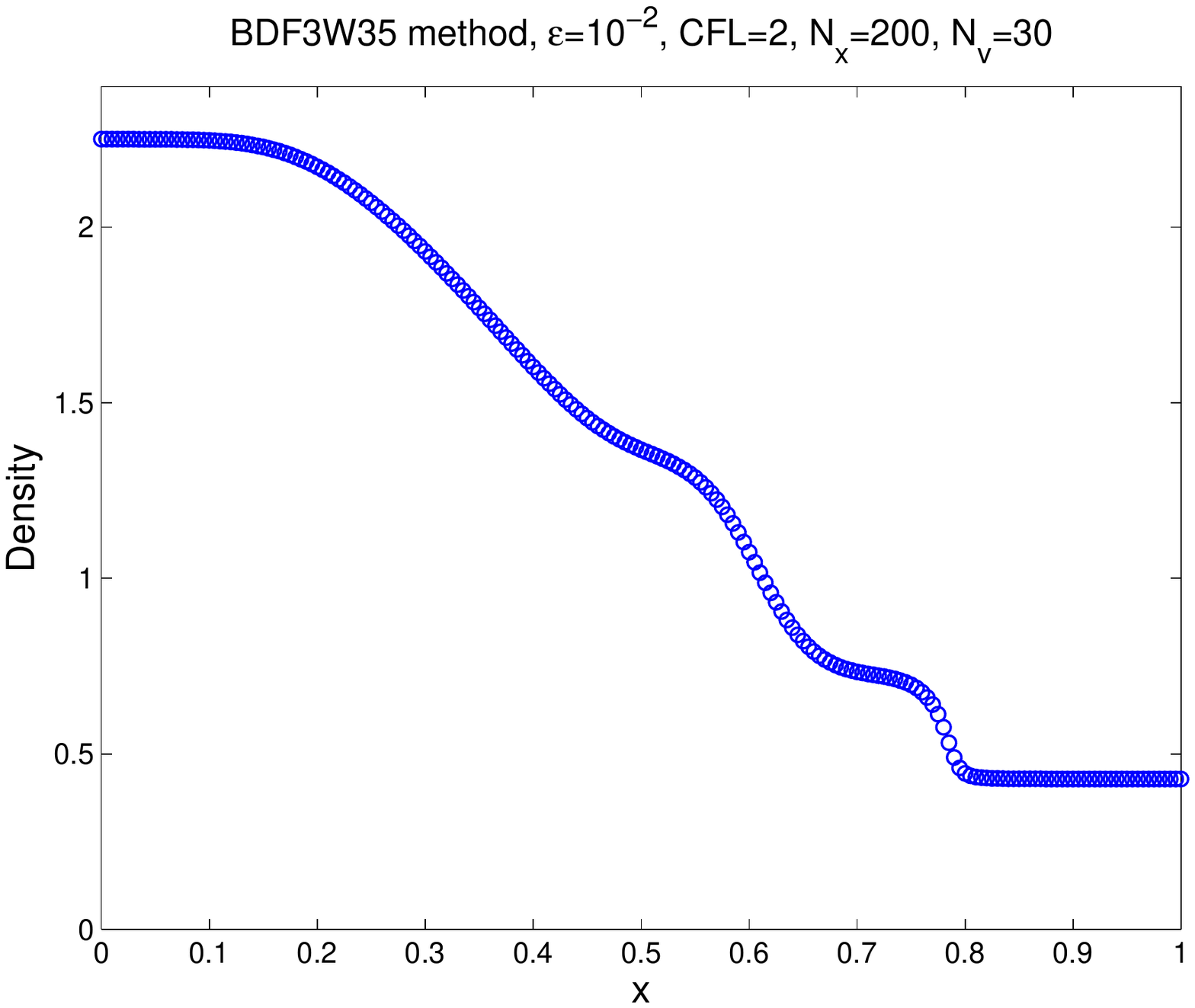}}
\subfigure
{\includegraphics[trim=1cm 7cm 2cm 12cm,scale=0.45]{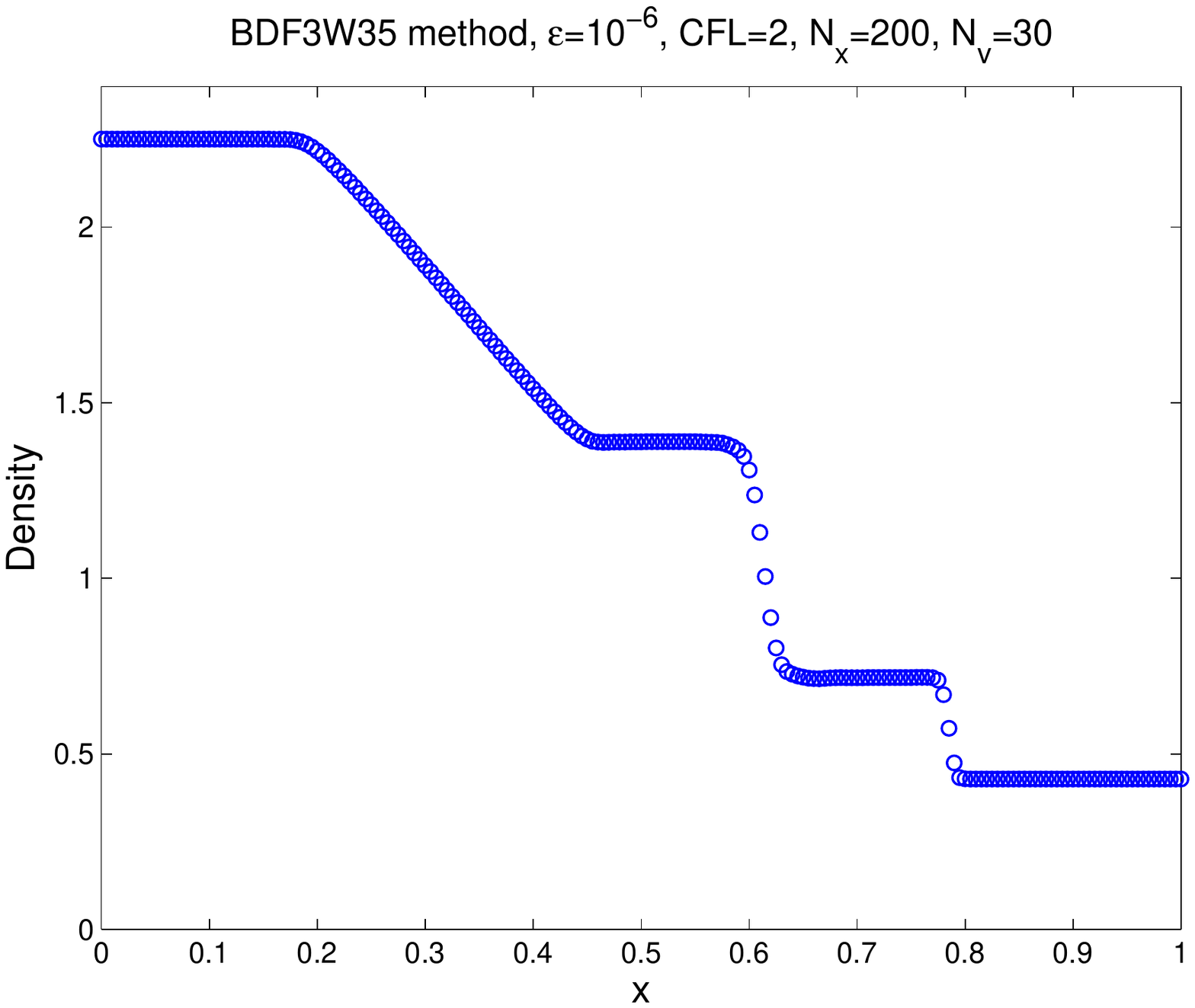}}
\subfigure
{\includegraphics[trim=5cm 7cm 2cm 7cm,scale=0.45]{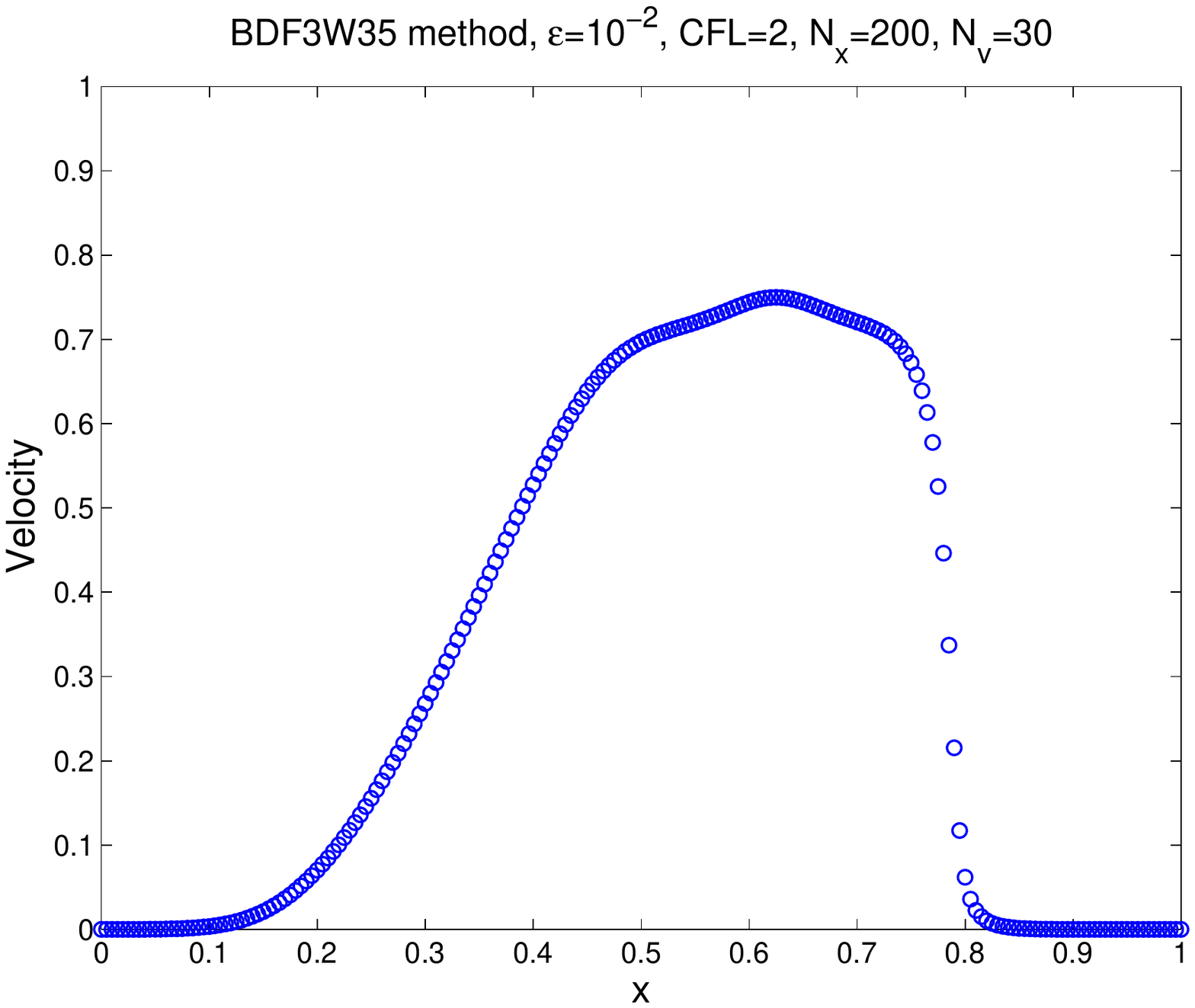}}
\subfigure
{\includegraphics[trim=1cm 7cm 2cm 7cm,scale=0.45]{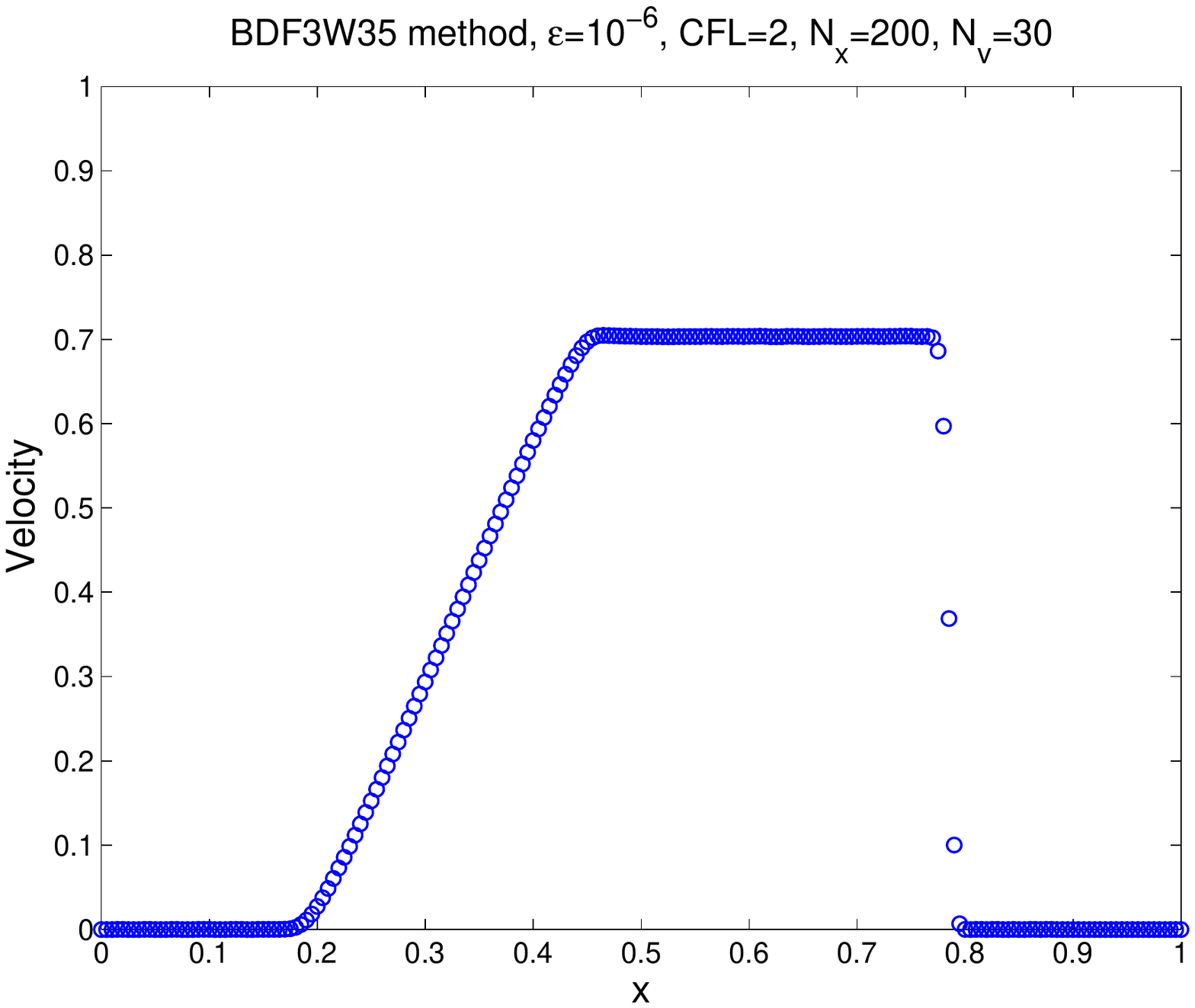}}
\subfigure
{\includegraphics[trim=5cm 7cm 2cm 7cm,scale=0.45]{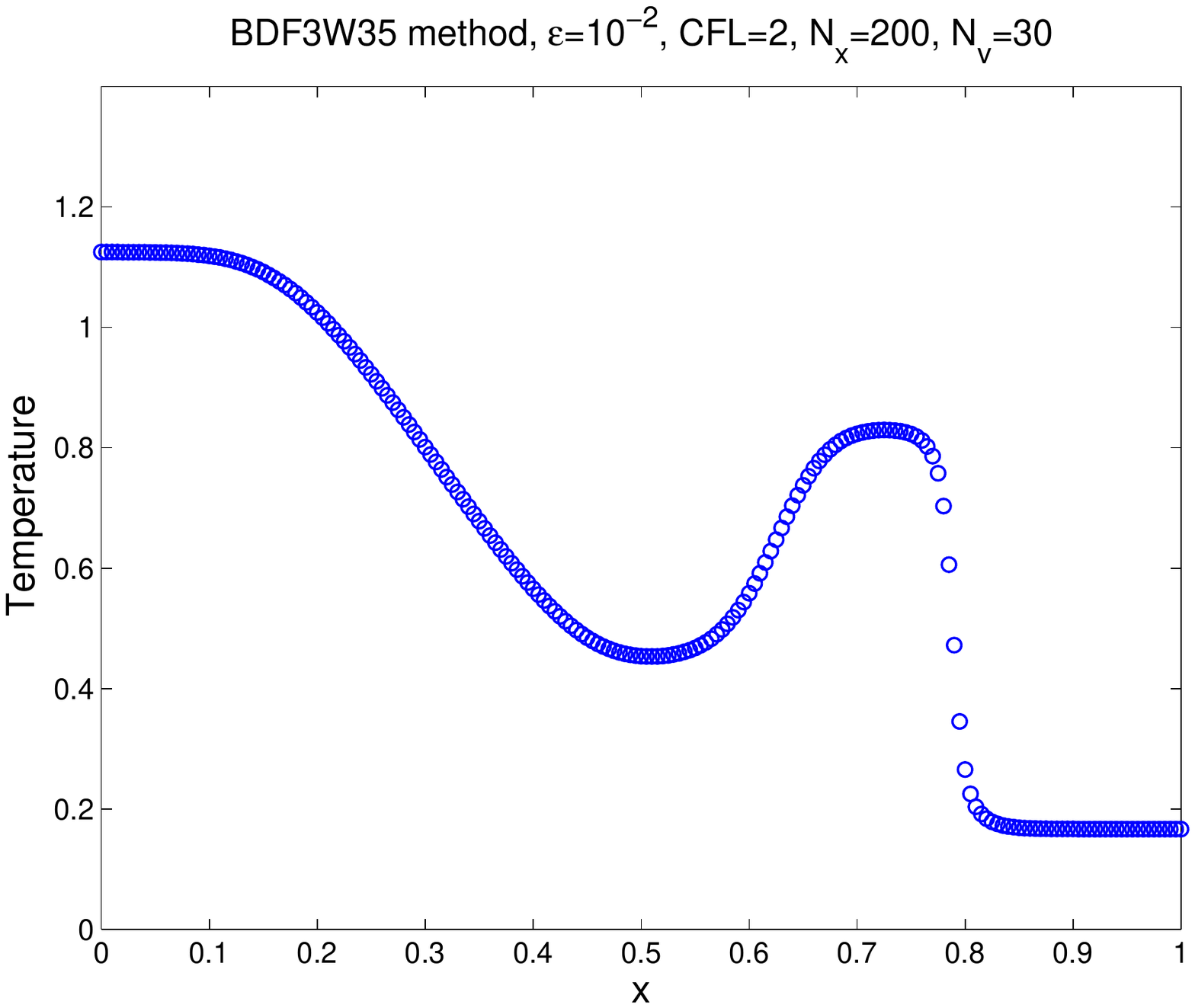}}
\subfigure
{\includegraphics[trim=1cm 7cm 2cm 7cm,scale=0.45]{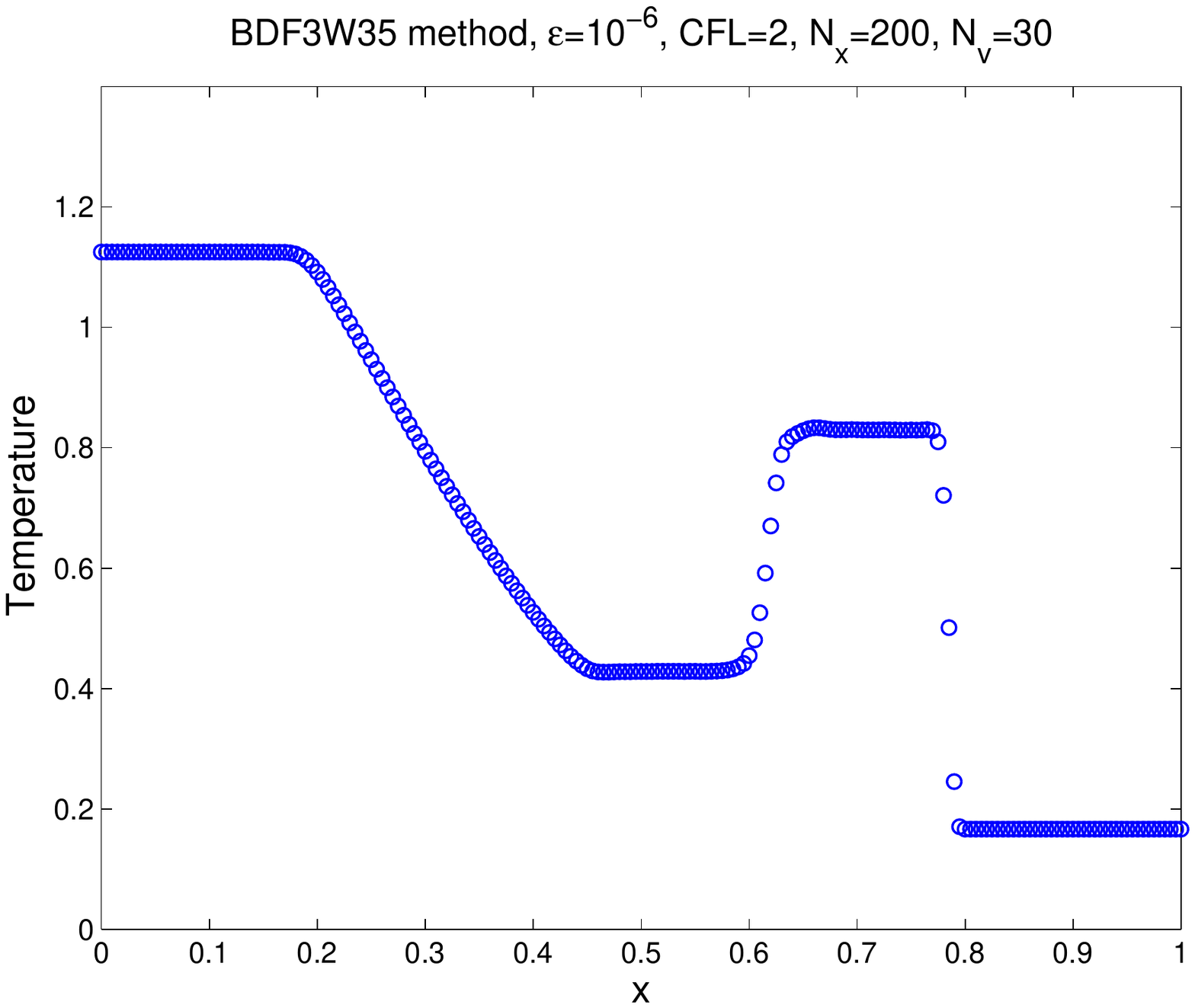}}
\caption{\footnotesize{BDF3W35 Riemann problem in 1D space and velocity case. Left $\varepsilon=10^{-2}$; Right $\varepsilon=10^{-6}$. From top to bottom: Density, Velocity and Temperature.}}\label{fig_BDF3_riemann}
\end{figure}

This test allows us to evaluate the capability of our class of schemes in capturing shocks and contact discontinuities. In particular, we are interested in the behavior of the schemes in the fluid regime. Here we illustrate the results obtained for moments, i.e density, mean velocity and temperature profiles, for $\varepsilon=10^{-2}$ and $\varepsilon=10^{-6}$ (see Fig. \ref{fig_RK3_riemann}-\ref{fig_BDF3_riemann}). The spatial domain chosen is $[0,1]$ and the discontinuity is taken at $x=0.5$. The initial condition for the distribution function is the Maxwellian computed with the following moments: $(\rho_{L},u_{L},T_{L})=(2.25,0,1.125)$, $(\rho_{R},u_{R},T_{R})=(3/7,0,1/6).$ Free-flow boundary conditions are assumed. The final time is 0.16. These tests have been performed using $N_{v}=30$ velocity nodes, uniformly spaced in $[-10,10].$\\
As it appears from Fig. \ref{fig_RK3_riemann} and \ref{fig_BDF3_riemann}, the schemes are able to capture the fluid dynamic limit for very small values of the relaxation time, where the evolution of the moments is governed by the Euler equations.

\subsection{Semi-lagrangian schemes without interpolation}
\begin{figure}[h]
\subfigure
{\includegraphics[trim=0cm 7cm 4cm 8cm,scale=0.65]{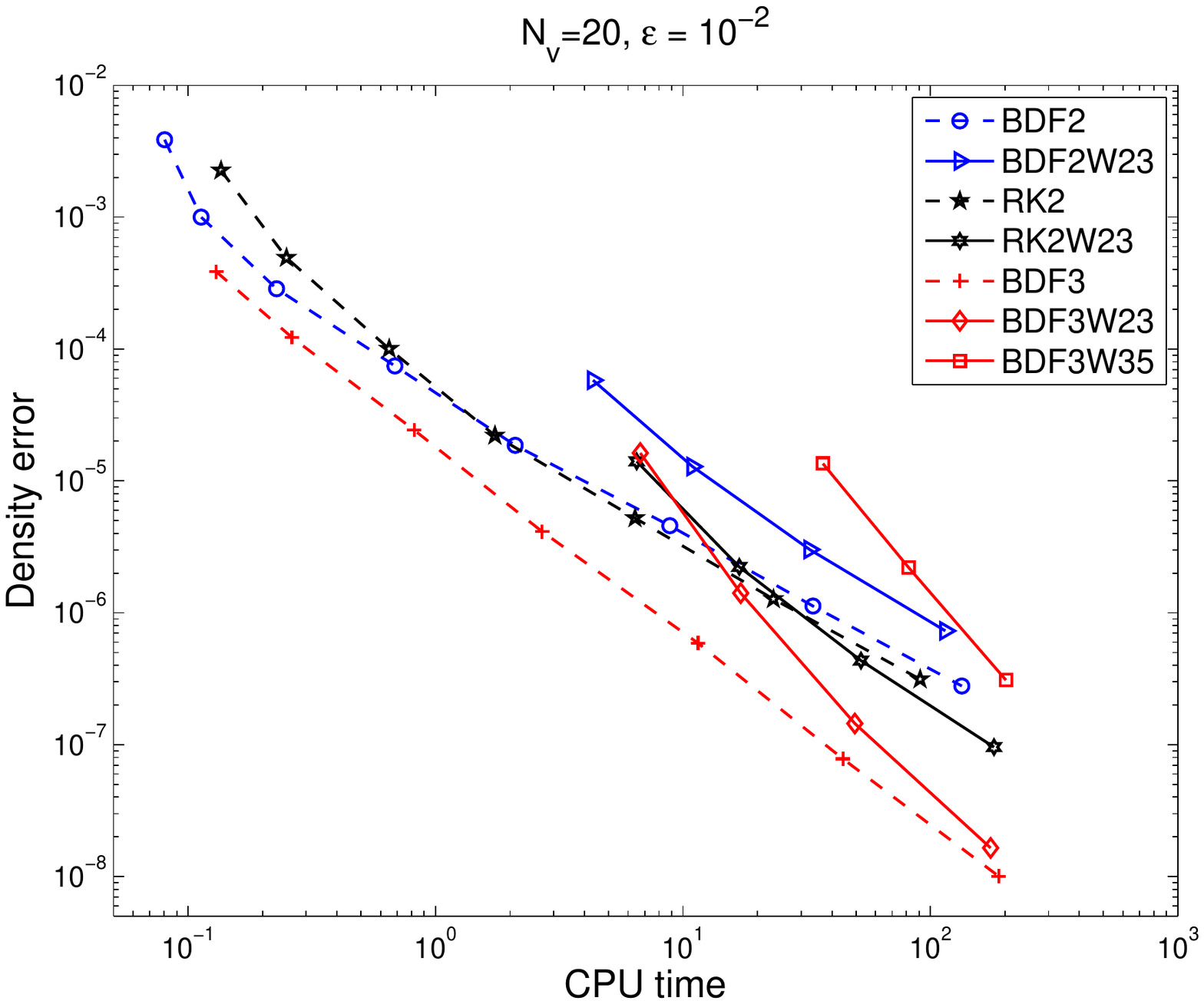}}
\caption{\footnotesize{Cpu time and $L_{1}$ error varying $N_{x}$.}}\label{fig_without_int}
\end{figure}

These schemes are very advantageous from a computational point of view. In Figure \ref{fig_without_int} we compare the cpu time and the $L_{1}$ error of the schemes with and without interpolation. At the third order of accuracy the relation between cpu time and error is better for the scheme without interpolation using $N_{v}=20$. The relative effectiveness of such schemes with respect to the ones that require interpolation decreases when increasing the number of velocities. However, these results are just indicative, as the schemes should be implemented efficiently.

\subsection{Numerical results - Chu reduction}
Also for the problem $3D$ in velocity we have considered two numerical test, that are aimed at verifying the accuracy and the shock capturing properties of the schemes. Different values of the Knudsen number have been investigated in order to observe the behavior of the methods varying from the rarefied to the fluid regime.\\ The initial data for test 1 are the same of the corresponding  test problem $1D$ in velocity, whereas the data for the second one are different. For the Riemann problem in this case, the initial condition for the distribution function is again a Maxwellian, having now the following initial macroscopic moments: $(\rho_{L},u_{L},T_{L})=(1,0,5/3)$, $(\rho_{R},u_{R},T_{R})=(1/8,0,4/3).$ As in the previous cases, free-flow boundary conditions are imposed. The final time is 0.25. This test has been performed using $N_{v}=30$ velocity nodes uniformly spaced in $[-10,10].$\\

\begin{figure}[htbp]
\subfigure
{\includegraphics[trim=3cm 7cm 2cm 8.5cm,scale=0.4]{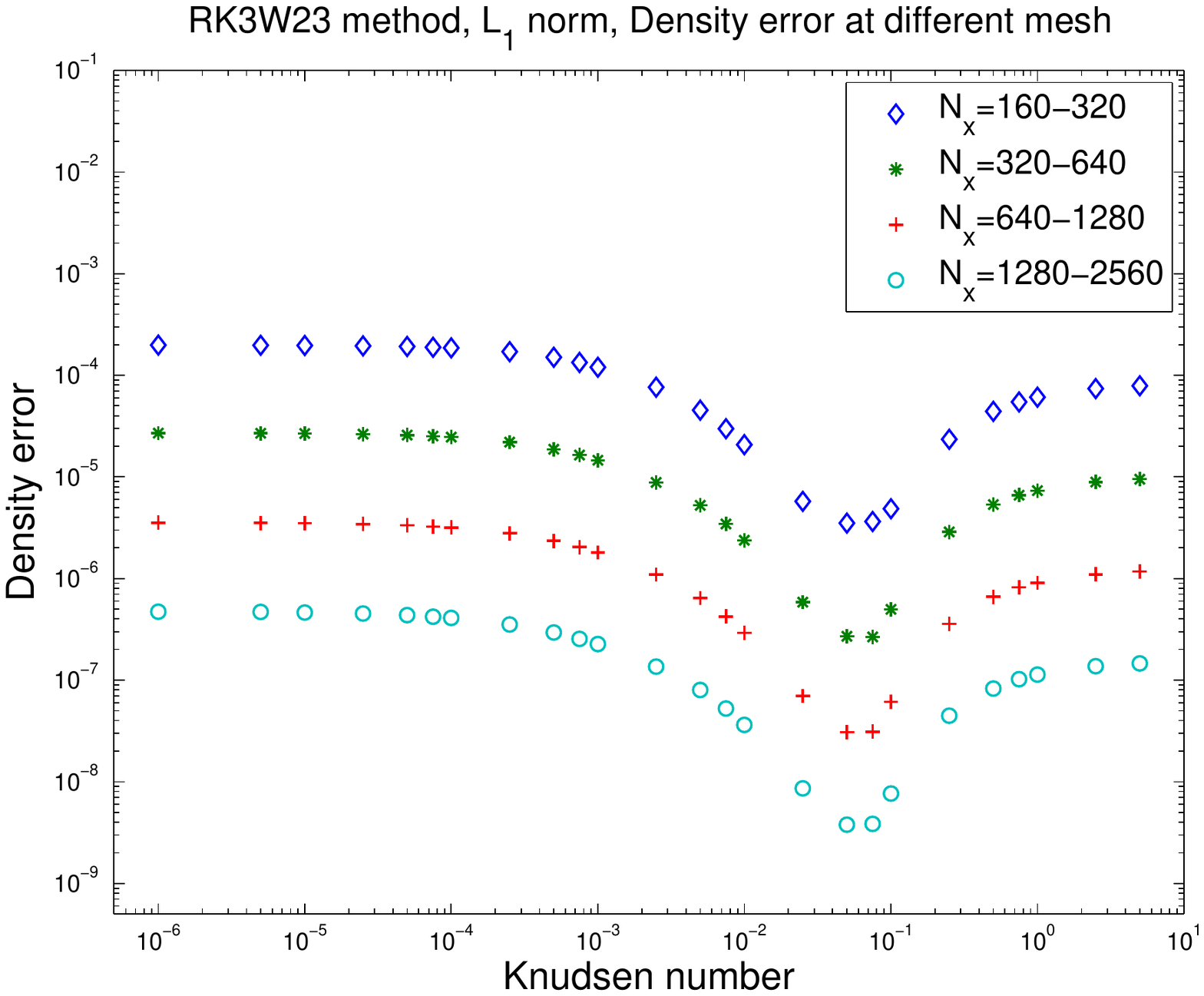}}
\subfigure
{\includegraphics[trim=1cm 7cm 2cm 8.5cm,scale=0.4]{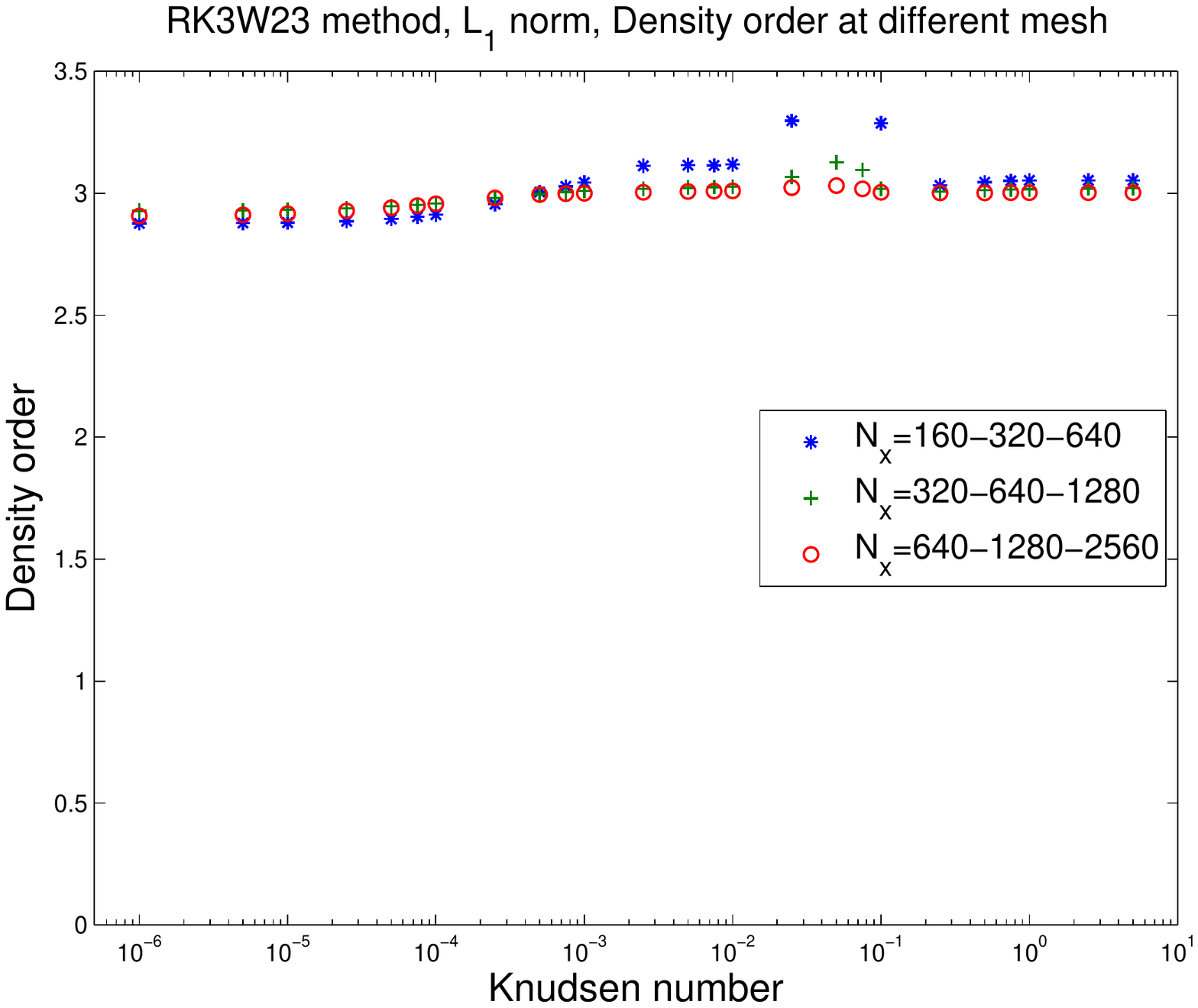}}
\subfigure
{\includegraphics[trim=3cm 7cm 2cm 7cm,scale=0.4]{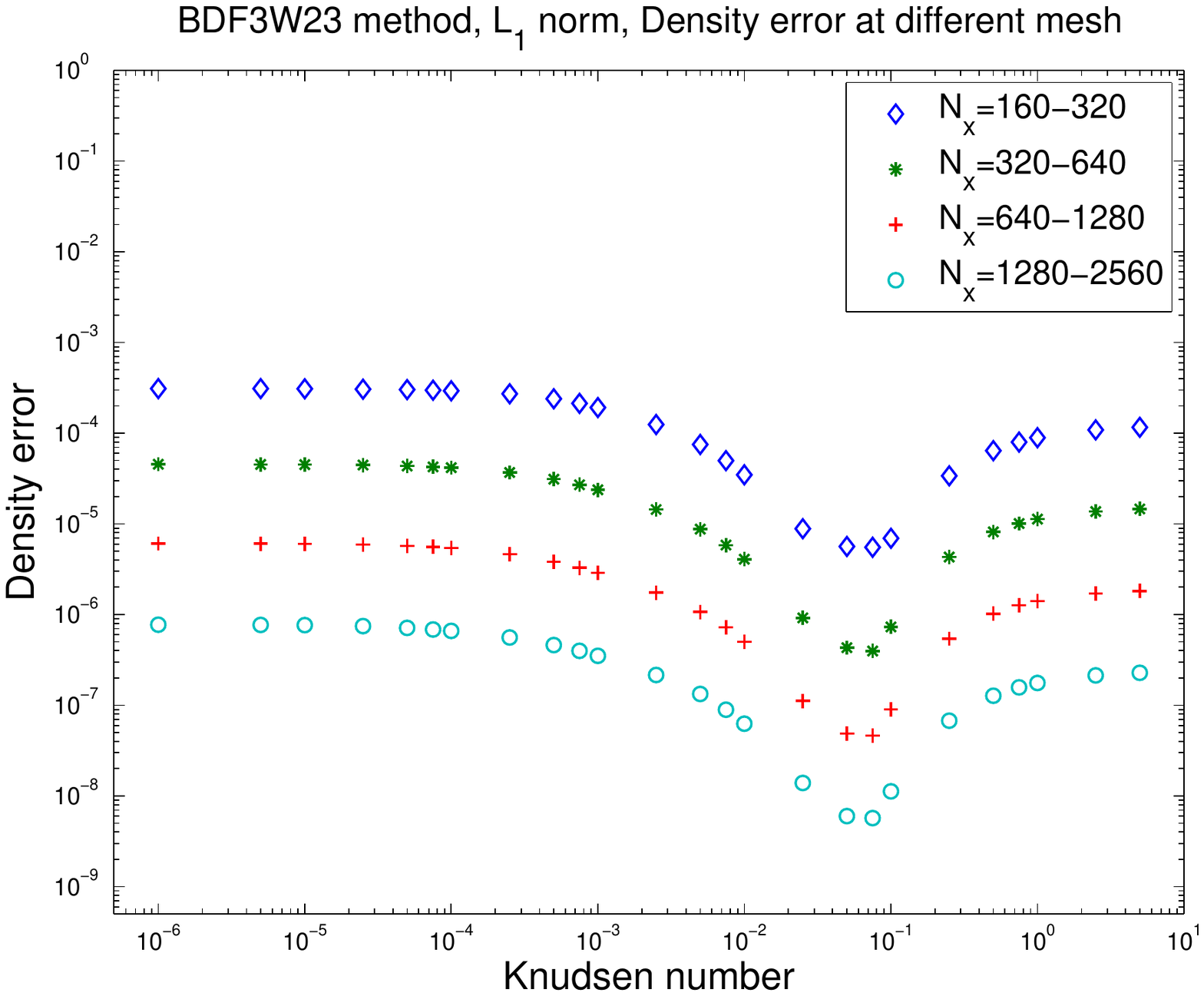}}
\subfigure
{\includegraphics[trim=1cm 7cm 2cm 7cm,scale=0.4]{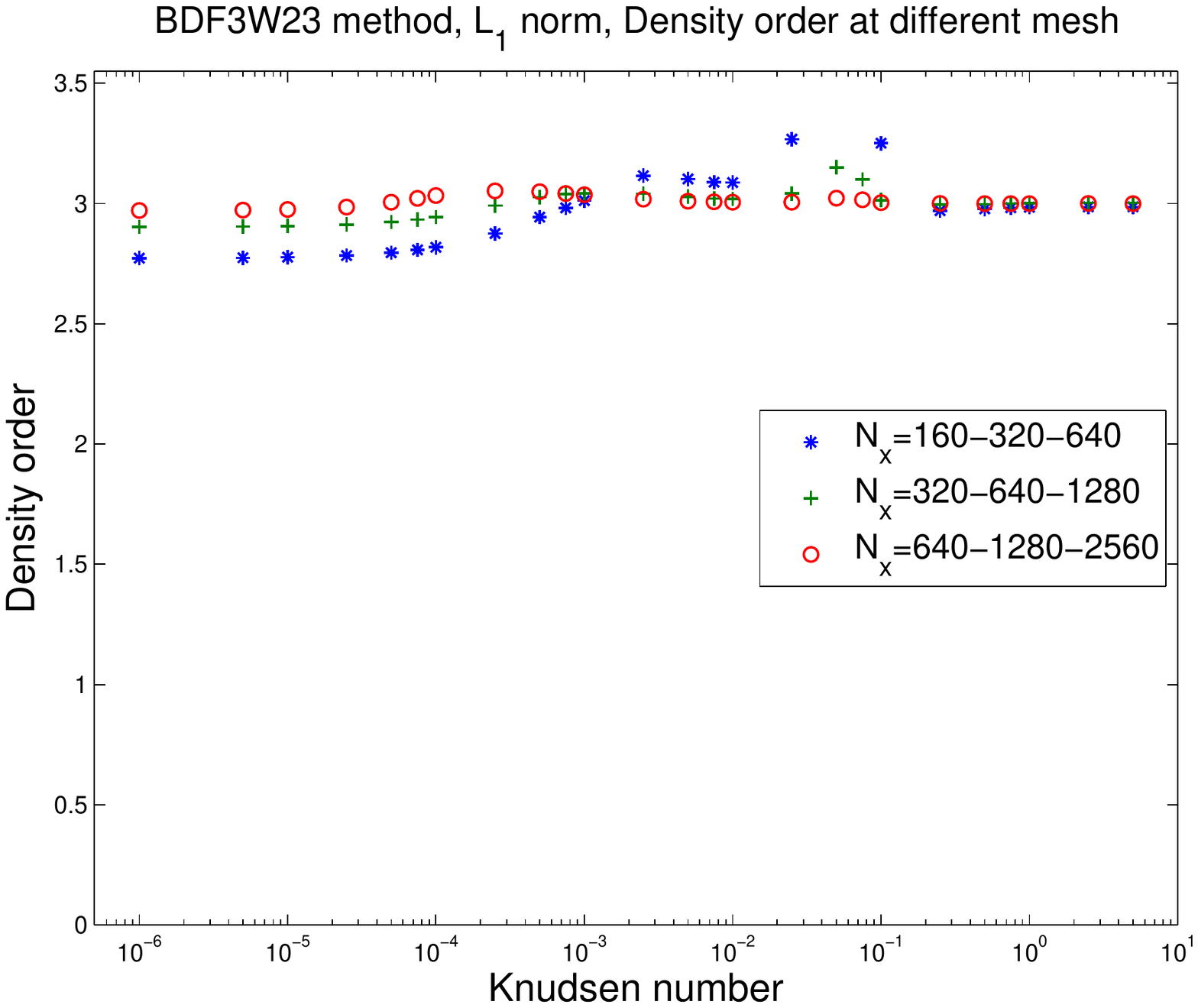}}
\caption{\footnotesize{$L_{1}$ error and accuracy order of RK3 and BDF3 methods coupled with WENO23, varying $\varepsilon$, using reflective boundary condition related to the 3D problem.}}\label{fig_third_order_3D}
\end{figure}
\begin{figure}[htbp]
\subfigure
{\includegraphics[trim=5cm 7cm 2cm 12cm,scale=0.45]{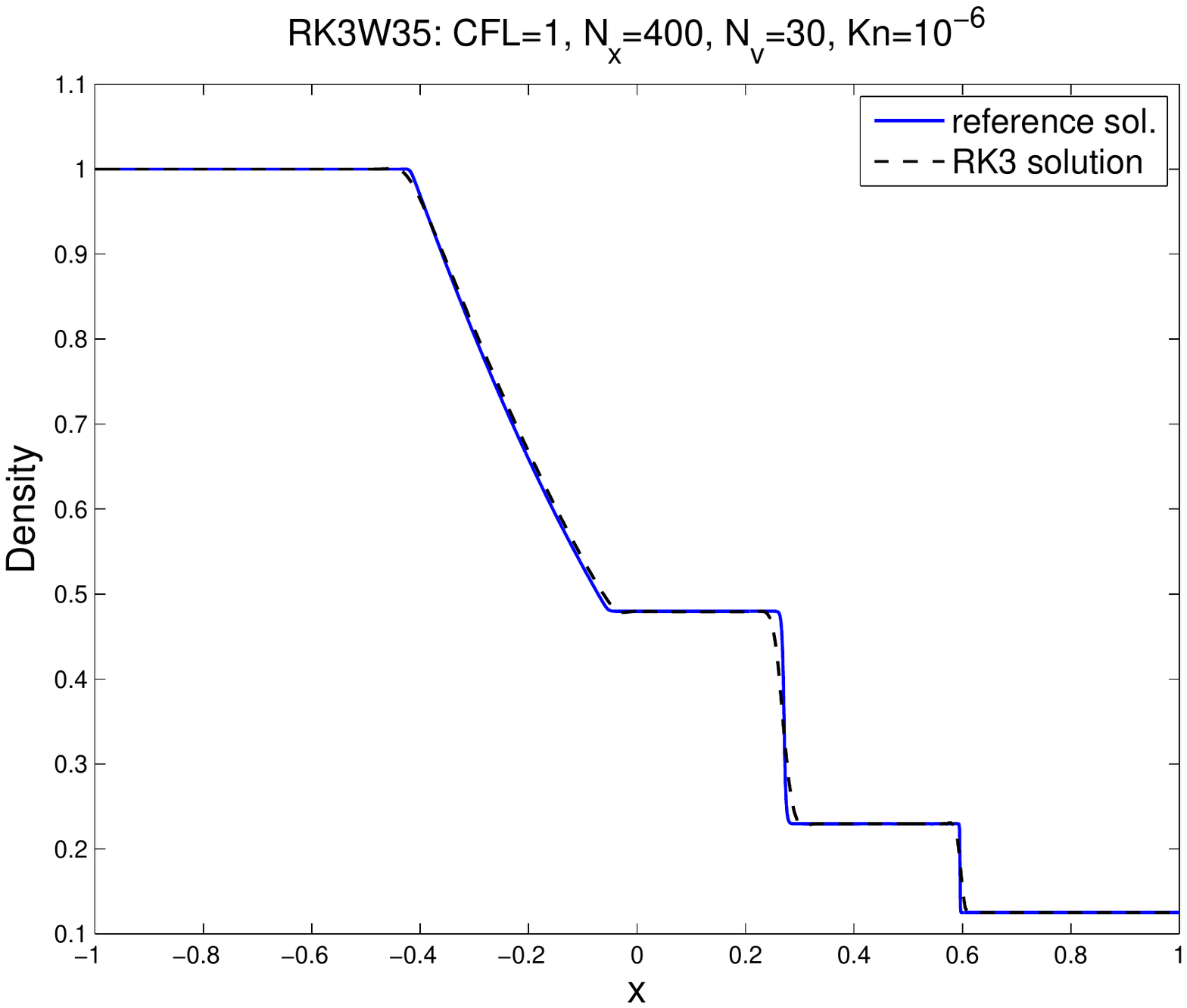}}
\subfigure
{\includegraphics[trim=1cm 7cm 2cm 12cm,scale=0.45]{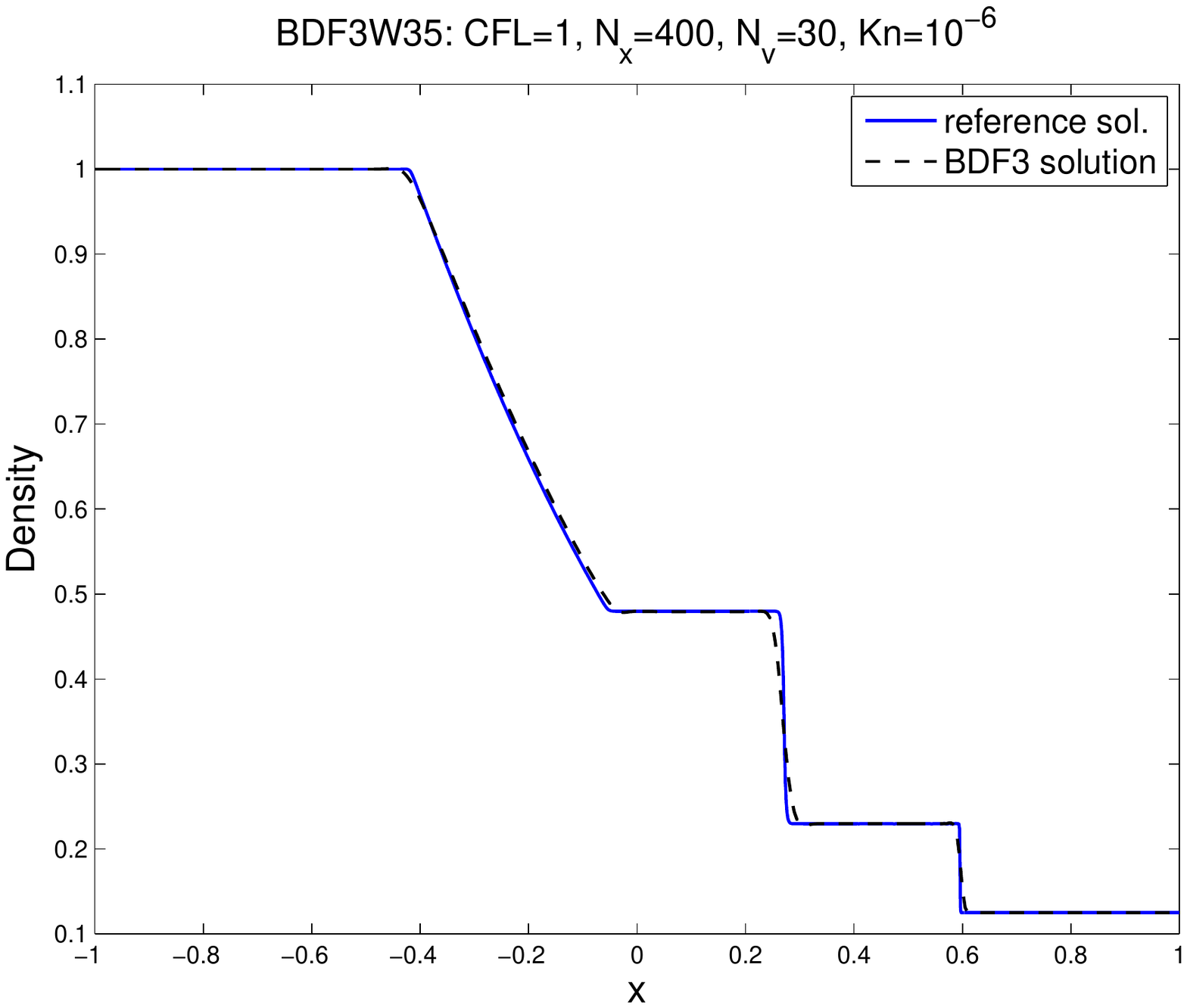}}
\subfigure
{\includegraphics[trim=5cm 7cm 2cm 7cm,scale=0.45]{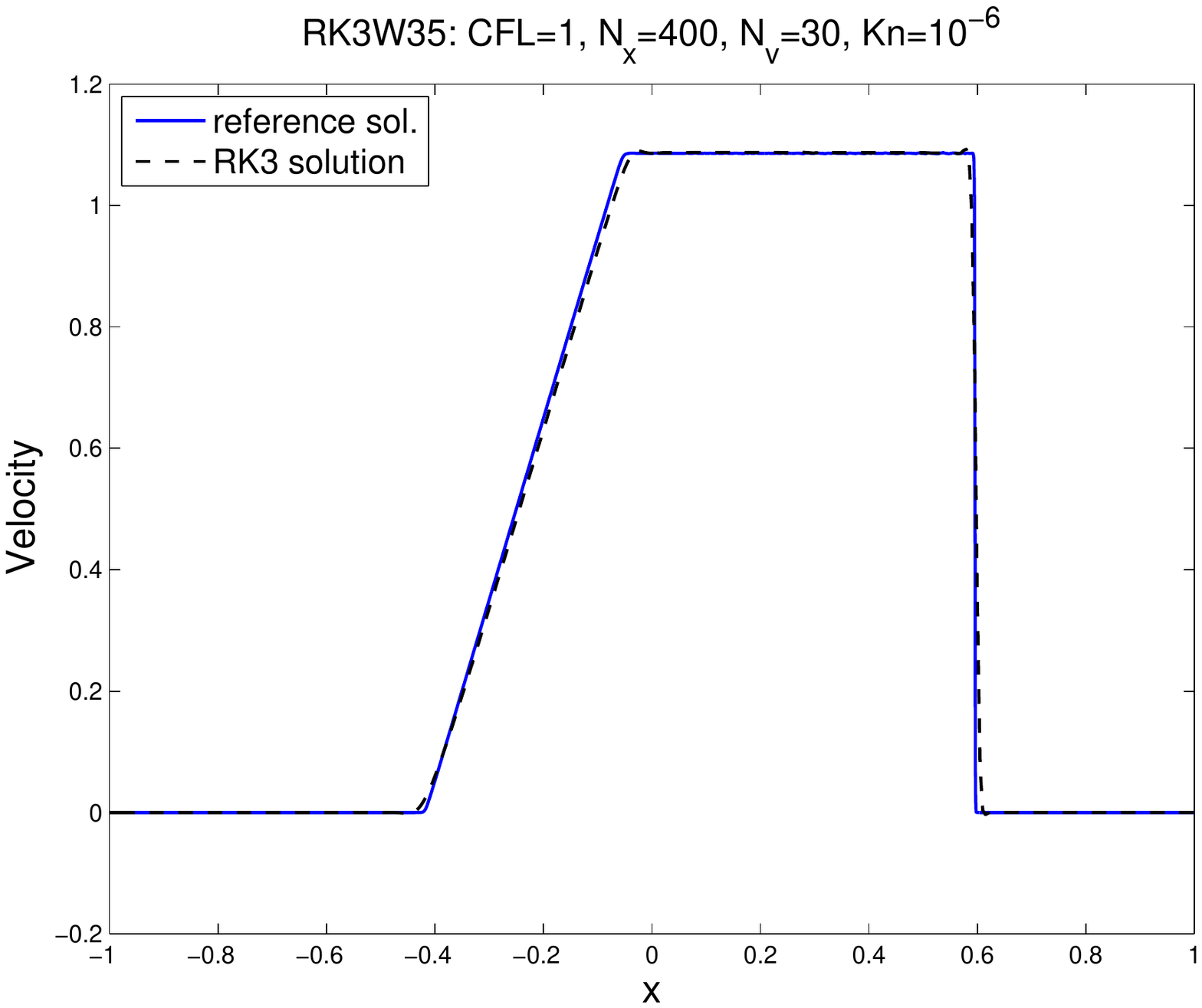}}
\subfigure
{\includegraphics[trim=1cm 7cm 2cm 7cm,scale=0.45]{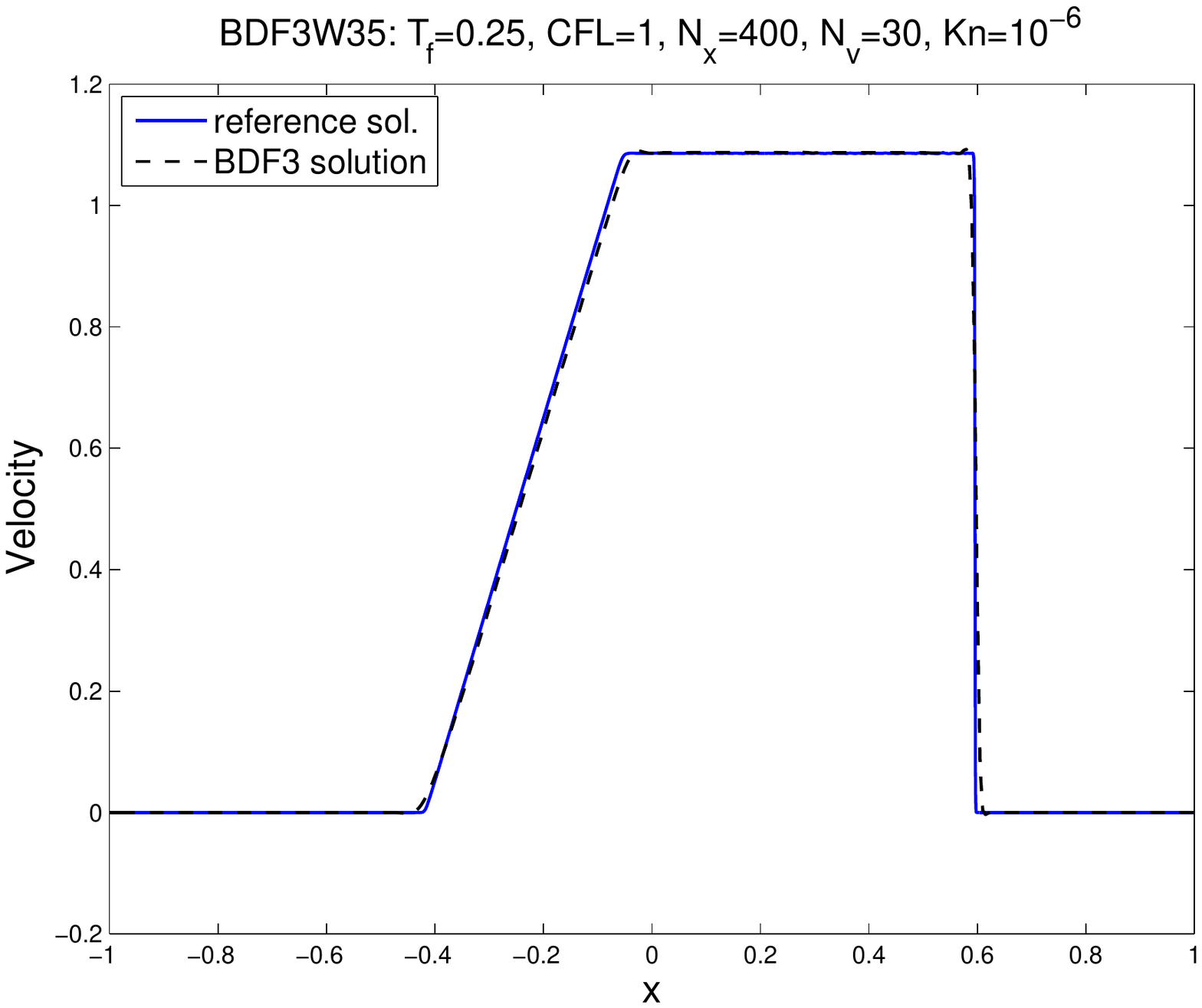}}
\subfigure
{\includegraphics[trim=5cm 7cm 2cm 7cm,scale=0.45]{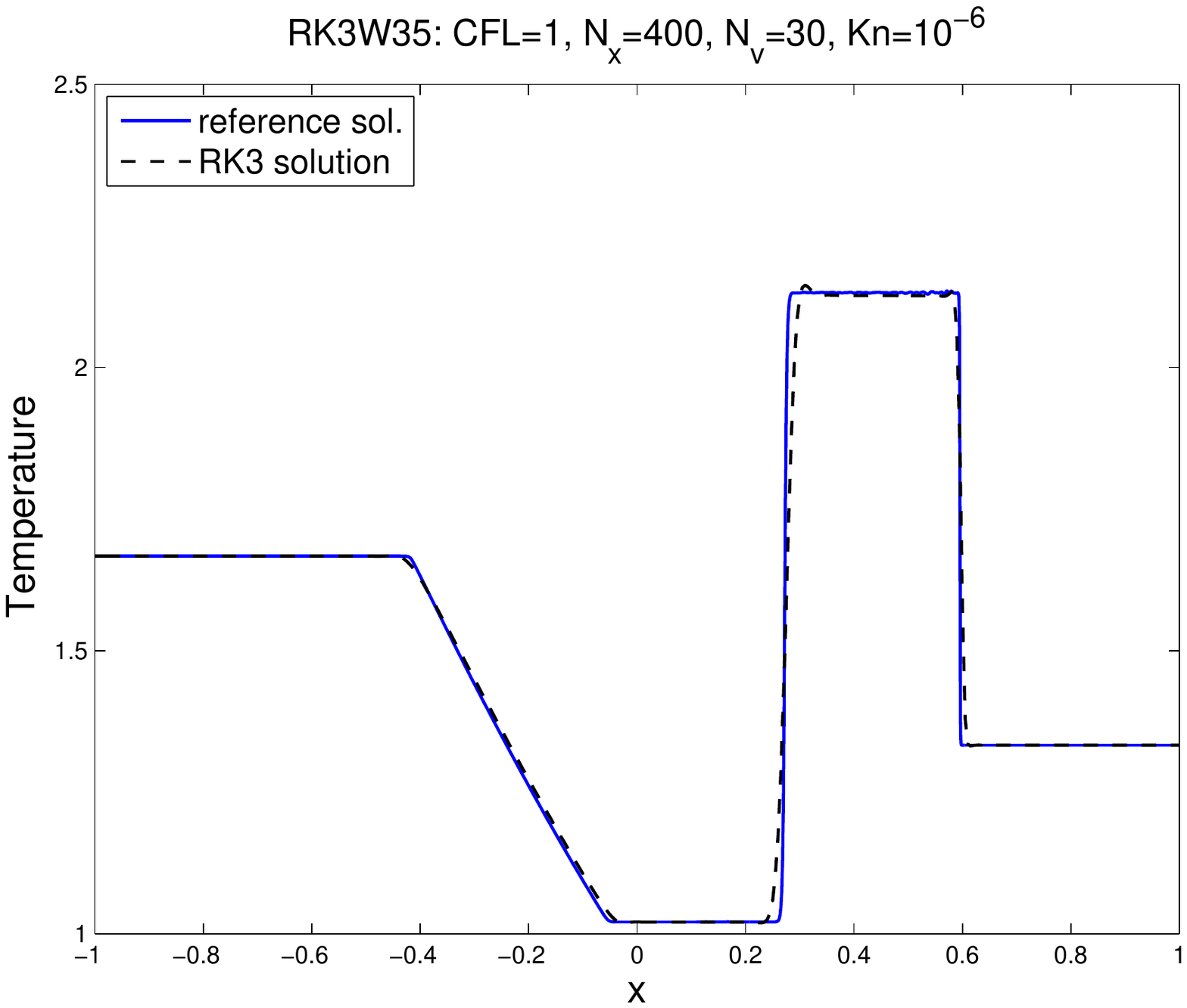}}
\subfigure
{\includegraphics[trim=1cm 7cm 2cm 7cm,scale=0.45]{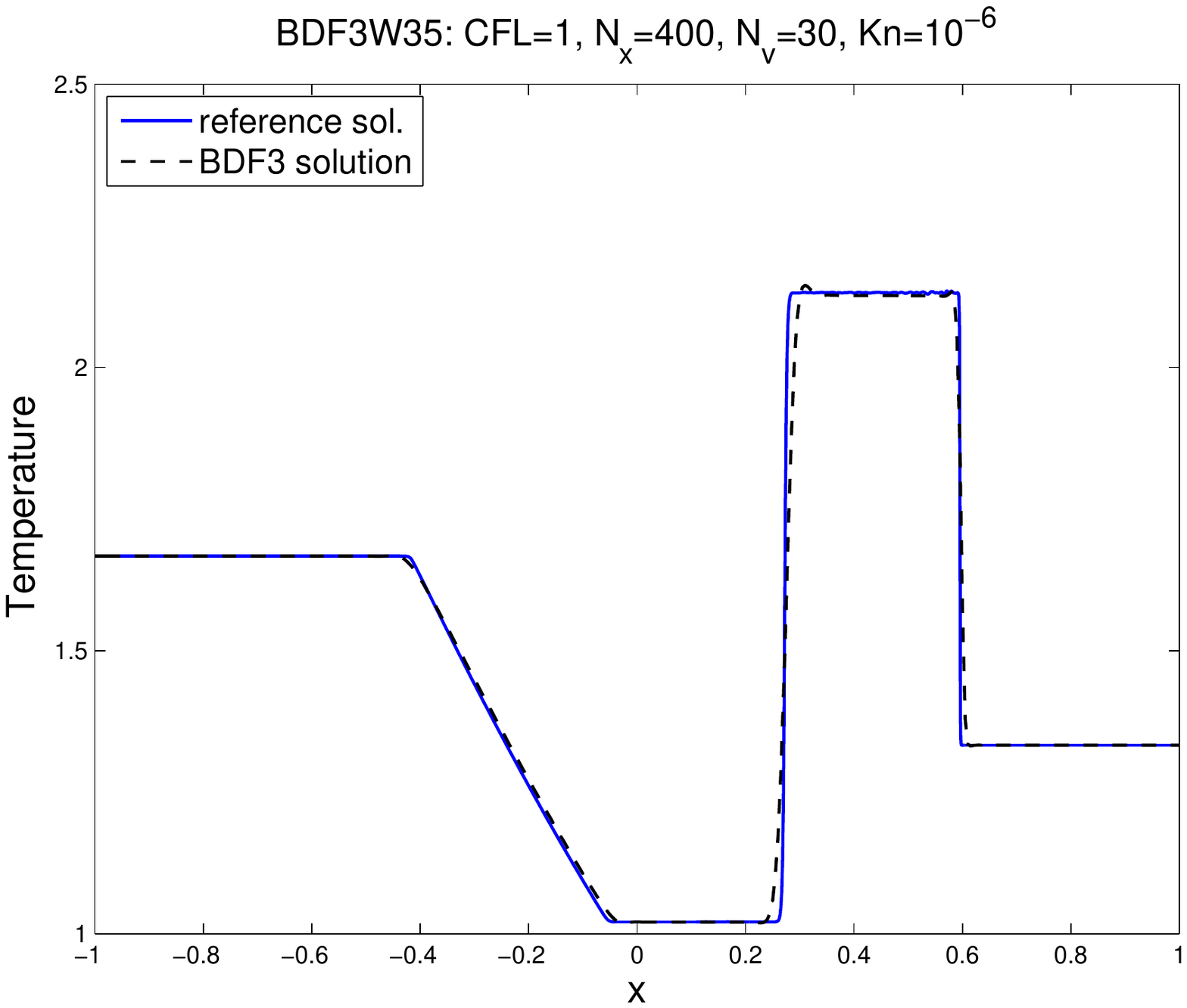}}
\caption{\footnotesize{Comparison with the gas dynamics. Left RK3W35; Right BDF3W35. From top to bottom: Density, Velocity and Temperature.}}\label{fig_comparison_gas_dynamics}
\end{figure}

\noindent
We will show only the order of accuracy related to the schemes RK3W23 and BDF3W23 (Fig. \ref{fig_third_order_3D}) using reflective boundary condition, in order to not be repetitive, as we get the same results of the 1D problem. In this test CFL$=2$ and the final time is 0.4. Regarding the Riemann problem we will show a comparison with the solution of the gas dynamics, for $\varepsilon=10^{-6},$ see Fig. \ref{fig_comparison_gas_dynamics}. As it appears from the results, also in this case the scheme is able to capture the fluid dynamic limit for very small values of the relaxation time, where the evolution of the moments is governed by the Euler equations.

\section{Appendix}
In order to obtain high order accuracy and to ensure the shock capturing properties of the proposed schemes near the fluid regime,  a suitable nonlinear reconstruction technique for the computation of $\tilde{f}^{n}_{ij}$ is required. ENO (essentially non oscillator) and WENO (weighted ENO) methods \cite{Shu} provide the desired high accuracy and non oscillatory properties. Both methods are based on the reconstruction of piecewise smooth functions by choosing the interpolation points on smooth side of the function. In ENO methods these points are chosen according to the magnitude of the divided differences evaluated by two candidate stencils. In WENO methods the different polynomials defined on the stencils are weighted in such a way that the information about the function from both sides can be used. Here we focus on WENO reconstruction \cite{WENO} by introducing the general framework for the implementation.

\subsection{Second-third order WENO interpolation (WENO23)}
To construct a third order interpolation we start from two polynomials of degree two, so that $$I[V^{n}](x)=\omega_{L}P_{L}(x)+\omega_{R}P_{R}(x),$$ where $P_{L}(x)$ and $P_{R}(x)$ are second order polynomials relevant to nodes $x_{j-1},x_{j},x_{j+1}$ and  $x_{j},x_{j+1},x_{j+2}$, respectively. The two linear weights $C_{L}$ and $C_{R}$ are first degree polynomials in $x$, and according to the general theory outlined so far, they read as $$C_{L}=\frac{x_{j+2}-x}{3\Delta x}, \,\,\,\,\, C_{R}=\frac{x-x_{j-1}}{3\Delta x};$$
the expressions of $\alpha_{L}$, $\alpha_{R}$, $\omega_{L}$ and $\omega_{R}$ may be easily recovered from the general form.

The smoothness indicators have the following explicit expressions
$$\beta_{L}=\frac{13}{12}v_{j-1}^{2} +\frac{16}{3}v_{j}^{2}+\frac{25}{12}v_{j+1}^{2} -\frac{13}{3}v_{j-1}v_{j}+\frac{13}{6}v_{j-1}v_{j+1} - \frac{19}{3}v_{j}v_{j+1},$$
$$\beta_{R}=\frac{13}{12}v_{j+2}^{2} +\frac{16}{3}v_{j+1}^{2}+\frac{25}{12}v_{j}^{2} -\frac{13}{3}v_{j+2}v_{j+1}+\frac{13}{6}v_{j+2}v_{j} - \frac{19}{3}v_{j}v_{j+1},$$
where
\begin{equation}\alpha_{k}(x)=\frac{C_{k}(x)}{(\beta_{k}+\epsilon)^{2}} \label{coef_alpha}
\end{equation}
(with $\epsilon$ a properly small parameter, usually of the order of $10^{-6}$), and then the nonlinear weights as  \begin{equation}
\omega_{k}=\frac{\alpha_{k}(x)}{\sum_{l}\alpha_{l}(x)}.\label{coef_omega}
\end{equation}

\subsection{Third-fifth order WENO interpolation (WENO35)}
To construct a fifth order interpolation we start from three polynomials of third degree: $$I[V^{n}](x) = \omega_{L}P_{L}(x) + \omega_{C}P_{C}(x) + \omega_{R}P_{R}(x),$$ where the third order polynomials $P_{L}(x)$, $P_{C}(x)$ and $P_{R}(x)$ are constructed, respectively, on $x_{j-2}$, $x_{j-1}$, $x_{j}$, $x_{j+1}$, on $x_{j-1}$, $x_{j}$, $x_{j+1}$, $x_{j+2}$, and on $x_{j}$, $x_{j+1}$, $x_{j+2}$, $x_{j+3}$. The weights $C_{L}$, $C_{C}$ and $C_{R}$ are second degree polynomials in $x$, and have the form
$$C_{L}=\frac{(x-x_{j+2})(x-x_{j+3})}{20\Delta x^{2}},\quad C_{C}=-\frac{(x-x_{j-2})(x-x_{j+3})}{10\Delta x^{2}},$$$$
C_{R}=\frac{(x-x_{j-2})(x-x_{j-1})}{20\Delta x^{2}},$$

while the smoothness indicators $\beta_{C}$ and $\beta_{R}$ have the expressions
$$\beta_{C}=\frac{61}{45}v_{j-1}^{2} + \frac{331}{30}v_{j}^{2} + \frac{331}{30}v_{j+1}^{2} + \frac{61}{45}v_{j+2}^{2} - \frac{141}{20}v_{j-1}v_{j} + \frac{179}{30}v_{j-1}v_{j+1} $$$$ - \frac{293}{180}v_{j-1}v_{j+2} - \frac{1259}{60}v_{j}v_{j+1} + \frac{179}{30}v_{j}v_{j+2} - \frac{141}{20}v_{j+1}v_{j+2},$$
$$\beta_{R}=\frac{407}{90}v_{j}^{2} + \frac{721}{30}v_{j+1}^{2} + \frac{248}{15}v_{j+2}^{2} + \frac{61}{45}v_{j+3}^{2} - \frac{1193}{60}v_{j}v_{j+3} + \frac{439}{30}v_{j}v_{j+2} $$$$ - \frac{683}{180}v_{j}v_{j+3} - \frac{2309}{60}v_{j+1}v_{j+2} + \frac{309}{30}v_{j+1}v_{j+3} - \frac{553}{60}v_{j+2}v_{j+3},$$
and $\beta_{L}$ can be obtained using the same set of coefficients of $\beta_{R}$ in a symmetric way (that is, replacing the indices $j-2,\cdots,j+3$ with $j+3,\cdots,j-2$) and $\alpha_{k}$ and $\omega_{k}$ are computed as in (\ref{coef_alpha}) and in (\ref{coef_omega}).

\section{Conclusions}
This paper presents high order shock capturing  semilagrangian methods for the numerical solutions of BGK-type equations.\\ The methods are based on L-stable schemes for solution of the BGK equations along the characteristics, and are asymptotic preserving, in the sense that are able to solve the equations also in the fluid dynamic limit.\\ Two families of schemes are presented, which differ for the choice of the time integrator: Runge-Kutta or BDF. A further distinction concerns space discretization: some schemes are based on high order reconstruction, while other are constructed on the lattice in phase space, thus requiring no space interpolation.\\ Numerical experiments show that schemes without interpolation can be cost-effective, especially for problems that do not require a fine mesh in velocity. In particular, BDF3 without interpolation appears to have the best performance in most tests.\\ Future plans includes to extend such schemes to problem in several space dimension and treat more general boundary condition.

\end{document}